\documentclass[]{amsbook}
\usepackage{amssymb}
\usepackage{amsmath}
\usepackage{amsfonts}
\usepackage{amsthm}
\usepackage{amscd}

\theoremstyle{plain}
\newtheorem{thm}{Theorem}[subsection]
\newtheorem{lem}[thm]{Lemma}
\newtheorem{cor}[thm]{Corollary}
\newtheorem{prop}[thm]{Proposition}
\theoremstyle{definition}
\newtheorem{defn}[thm]{Definition}

\newtheorem{rem}[thm]{Remark}

\newcommand{\bZ}{{\mathbb Z}}
\newcommand{\bQ}{{\mathbb Q}}
\newcommand{\bR}{{\mathbb R}}
\newcommand{\bI}{{\mathbb I}}
\newcommand{\bJ}{{\mathbb J}}
\newcommand{\bK}{{\mathbb K}}
\newcommand{\bH}{{\mathbb H}}
\newcommand{\fm}{{\mathfrak m}}
\newcommand{\fD}{{\mathfrak D}}
\newcommand{\fR}{{\mathfrak R}}
\newcommand{\fB}{{\mathfrak B}}
\newcommand{\cI}{{\mathcal I}}
\newcommand{\cH}{{\mathcal H}}
\newcommand{\cO}{{\mathcal O}}
\newcommand{\gpure}[1]{F^{#1}(G_1)}
\newcommand{\on}[1]{\mathop{\mathrm{#1}}\nolimits}
\newcommand{\dif}{\on{Diff}}

\renewcommand{\mid}{\hskip.02in;\hskip.02in}

\renewcommand{\labelenumi}{(\arabic{enumi})}

\makeatletter
\newcounter{step}
\newenvironment{step}[1]{
\refstepcounter{step}\par
{\bf Step}\ $\mathbf\thestep$.\quad{\it #1}
\par\noindent}{}
\newcounter{case}
\newenvironment{case}[1]{
\refstepcounter{case}\par
\underline{Case\ $\thecase$}.\quad{\it #1}
\par\noindent}{}
\renewenvironment{proof}[1][\proofname]{\par
\normalfont \topsep6\p@\@plus6\p@\relax
\trivlist\itemindent\normalparindent
\item[\hskip\labelsep\scshape 
#1\@addpunct{.}]\ignorespaces
\setcounter{step}0
\setcounter{case}0
}\@endtheorem
\makeatother

\newcounter{komoku}[subsection] 
\renewcommand{\thekomoku}{\thechapter.\arabic{section}.%
\arabic{subsection}.\arabic{komoku}} 
\newcommand{\komoku}[1]{\refstepcounter{komoku}%
\hskip0in\par\noindent\thekomoku\ {\bf#1}\ } 
\newcounter{subkomoku}[komoku] 
\renewcommand{\thesubkomoku}{\thechapter.\arabic{section}.%
\arabic{subsection}.\arabic{komoku}.\arabic{subkomoku}} 
\newcommand{\subkomoku}[1]{\refstepcounter{subkomoku}%
\hskip0in\par\noindent\thesubkomoku\ {\bf#1}\ } 

\addtocontents{toc}{{\noindent\bf Part I.\quad Foundation; %
the language of the idealistic filtration}}
\makeatletter
\renewcommand{\tocsubsection}[3]{%
\phantom{\S}{\@ifnotempty{#2}{\ignorespaces#1 #2.\quad}}#3}
\makeatother
\begin{document}
\title{
{\huge\sc Toward resolution of singularities\linebreak
over a field of\linebreak 
positive characteristic} 
\vskip.6in 
{\it Dedicated to Professor Heisuke Hironaka} 
\vskip1.2in 
{Part\ I.
\vskip0in
Foundation; the language of the idealistic filtration}
\vskip3.3in 
}
\author{\sc Hiraku\ Kawanoue}
\date{today}
\maketitle
\tableofcontents
\begin{chapter}{Introduction} 
\begin{section}{Goal of this series of papers.}\label{0.1} 
This is the first of the series of papers under the title 
\begin{center} 
``Toward resolution of singularities  over a field of positive 
characteristic'' 
\begin{tabular}{ll} 
Part I. &Foundation; the language of the idealistic filtration 
\\ 
Part II. &Basic invariants associated to the idealistic filtration 
\\ 
&and their properties 
\\ 
Part III. &Transformations and modifications of the idealistic filtration 
\\ 
Part IV.  &Algorithm in the framework of the idealistic filtration. 
\end{tabular} 
\end{center} 
Our goal is to 
present a program toward constructing an algorithm for 
resolution of singularities of an algebraic variety 
over a perfect field $k$ of positive characteristic 
$p = \on{char}(k) >0$.  We would like to emphasize, 
however, that the 
program is created in the spirit of developing a uniform point of view 
toward the problem of resolution of singularities in 
all characteristics, and hence that it is also valid in characteristic zero. 
\footnote{During the preparation of the 
manuscript for Part I, we were informed that Professor 
Hironaka announced a 
program of resolution of singularities in all 
characteristics $p > 0$ and in all dimensions at the summer school in 
Trieste 2006 (cf.~\cite{H-Trieste}).} 
 
In Part I, we establish the notion and some fundamental properties of an 
{\it idealistic filtration}, which 
is the main language to describe the program.  This part, therefore, forms 
the foundation of the program. 
 
In Part II, we study the basic invariants $\sigma$ and $\widetilde{\mu}$ 
associated to an idealistic 
filtration, which will become the building blocks toward constructing the 
strand of invariants used in our 
algorithm, and discuss their properties. 
 
In Part III, we analyze the behavior of an idealistic filtration under the 
two main operations in the process 
of our algorithm for resolution of singularities: 
\begin{itemize} 
\item 
transformations of an idealistic filtration under the operation of 
blowup, and 
\item 
modifications of an idealistic filtration under the operation of 
constructing the strand of invariants. 
\end{itemize} 
 
Part II and Part III should play the role of a bridge between the foundation 
in Part I and the 
presentation of our algorithm in Part IV.

In Part IV, we present our algorithm for resolution of singularities 
according to the program as a summary 
of the series.  In characteristic zero, the program leads to a complete 
algorithm (slightly different from 
the existing ones), which then serves as a prototype toward the case in 
positive characteristic.  In positive 
characteristic, all the ingredients of the program work nicely forming a 
perfect parallel to the case in 
characteristic zero, except for the problem of termination: we do 
{\it not} know at this point whether our algorithm 
terminates after finitely many steps or not.  Although we do know that the 
strand of invariants we construct 
strictly drops after each blowup, we can not exclude the possibility that 
the denominators of some 
invariants in the strand may indefinitely increase and hence that the 
descending chain condition may not be 
satisfied.  The problem of termination remains as the only 
missing piece toward completing our algorithm in positive characteristic. 
We hope, however, that we may be able to come up 
with a solution to the problem during the process of writing down all the 
details of the program in this series of papers. 
\end{section} 
\begin{section}{Overview of the program.}\label{0.2} 
Below we present an 
overview of the program, by first 
giving a crash course on the existing algorithm(s) in characteristic zero, 
then pinpointing the main source 
of troubles if we try to apply the same methods to the case in positive 
characteristic, and finally 
describing how our program attempts to overcome these troubles. 
\begin{subsection} 
{Crash course on the existing algorithm(s) in characteristic 
zero.}\label{0.2.1} 
\komoku{Standard reduction.}\label{0.2.1.1} 
By a standard argument 
free of characteristic, the problem 
of resolution of singularities of an abstract algebraic variety is reduced 
to, and reformulated as, the 
problem of transforming a given ideal 
$\cI \subset \cO_W$ on a nonsingular 
variety $W$ over $k$ into the one whose multiplicity (order) 
becomes lower than the aimed (or expected) multiplicity 
$a$ everywhere, through a sequence of blowups and through 
a certain transformation rule for the ideal.  We 
require that each center of blowup to be nonsingular and transversal to the 
boundary, which consists of the 
exceptional divisor and the strict transform of a simple normal crossing 
divisor $E$ on 
$W$ given at the beginning.  We call this reformulation the problem of 
resolution of 
singularities of the triplet $(W, (\cI,a), E)$, and call 
$\on{Sing}(\cI,a) = \{P \in W 
\mid \on{ord}_P(\cI) \geq a\}$ 
its singular locus or support. 
\komoku{Inductive scheme in characteristic zero.}\label{0.2.1.2} 
At the very core of all the existing algorithmic approaches in 
characteristic zero lies the common inductive 
scheme on dimension, that is, reduce the problem of 
resolution of singularities of $(W,(\cI,a),E)$ 
to that of 
$(H, ({\mathcal J}, b), D)$, where $H$ is a smooth hypersurface in $W$.  The 
hypersurface $H$ is called a 
hypersurface of maximal contact, since it contains ({\it contacts}) the 
singular locus $\on{Sing}(\cI,a)$ 
and since so do its strict transforms throughout 
any sequence of 
transformations.  The ideal ${\mathcal J}$ on $H$ is usually realized as 
${\mathcal 
J} = C(\cI)|_H$, where 
$C(\cI)$ is the so-called coefficient ideal of the original ideal 
$\cI$, which is larger 
than 
$\cI$.  (It is worthwhile noting that the mere restriction ${\mathcal 
I}|_H$ of the 
original ideal would fail to provide the inductive scheme in general, and it 
is necessary to take a larger 
ideal.)  In short, we decrease the dimension by converting the problem on 
$W$ into the one on the 
hypersurface of maximal contact $H$ with 
$\dim H = \dim W -1$. 
\komoku 
{ 
Algorithm: modifications and construction of 
the strand of invariants.}\label{0.2.1.3} 
The above description of the inductive scheme is, 
however, oversimplified. For 
an arbitrary triplet $(W,({\mathcal 
I},a),E)$, a hypersurface of maximal contact may not exist at all. 
In order to gurantee that a hypersurface 
of maximal contact $H$ exists, we have to take the ``companion 
modification'' 
associated to the weak-order 
``$w$''.  Furthermore, in order to guarantee that $H$ is transversal to $E$ 
and hence that we can take $D = 
E|_H$, we have to take the ``boundary modification'' associated to the 
invariant ``$s$''.  In other words, only 
after considering the pair of invariants $(w,s)$ and taking the 
corresponding companion modification and 
its boundary modification, we can find the triplet 
$(H, ({\mathcal J},b), D)$ of dimension one less as in 
\ref{0.2.1.2}, whose resolution of singularities corresponds 
to the decrease of the pair of invariants $(w,s)$. 
 
Therefore, the actual algorithm realizing the inductive scheme is carried 
out in such a way that we 
construct the strand of invariants 
$$\on{\it inv}_{\on{classical}} 
= (w,s)(w,s)(w,s) \dotsb$$ 
by repeating the operations of taking the companion modification, boundary 
modification, and taking the 
restriction to a hypersurface of maximal contact, and that at the end the 
maximum locus 
of the strand $\on{\it inv}_{\on{classical}}$ 
of invariants coincides with the last 
hypersurface of maximal 
contact, which is hence nonsingular and which we choose as the center of 
blowup.  After the blowup, we repeat the same process.  We can repeat the 
process only finitely many 
times, since after each blowup the value of the strand of invariants 
strictly drops and since the set of its 
values satisfies the descending chain condition, leading to the termination 
of the algorithm.  (See, 
e.g.,\cite{MR1440306}\cite{MR1748620}%
\cite{MR1949115}\cite{MR2163383}\cite{AG0103120} 
for details of the construction of the strand of 
invariants and the corresponding modifications 
in the classical setting.) 
\end{subsection} 
\begin{subsection} 
{Trouble in positive characteristic.}\label{0.2.2} 
In positive characteristic, 
however, the examples by R.~Narasimhan 
\cite{MR684627}\cite{MR715853} 
and others demonstrate that there is {\it no} hope of finding a 
hypersurface of maximal contact in general 
(even after companion or boundary 
modification), as long as we require it to 
contain the singular locus and to be nonsingular.  This lack of a 
hypersurface of maximal contact and hence 
of an apparent inductive scheme is the main source of troubles, which 
allowed the problem in positive 
characteristic to elude any systematic attempt to find an algorithm for its 
solution so far. 
\end{subsection} 
\begin{subsection} 
{Our program: a new approach in the framework 
of the idealistic filtration}\label{0.2.3} 
Our program offers a new approach to overcome the main 
source of troubles in the language 
of the 
{\it idealistic filtration}, which is a refined extension of 
such classical notions as the idealistic exponent by Hironaka, 
the presentation by Bierstone-Milman, the basic object by Villamayor, 
and the marked ideal by W{\l}odarczyk.  We devote Part I of the series 
of papers to introducing the notion of an idealistic 
filtration, and to establishing its fundamental properties. 
\komoku{ 
What is an idealistic filtration? 
}\label{0.2.3.1} 
In the classical setting, we consider the pair 
$(\cI,a)$ consisting of an ideal $\cI \subset \cO_W$ on a 
nonsingular variety $W$ and the 
aimed multiplicity 
$a \in \bZ_{> 0}$.  Stalkwise at a point $P \in W$, this is equivalent 
to considering the collection 
of pairs $\{(f,a)\mid f \in \cI_P\}$. 
 
Suppose we interpret the pair $(f,a)$ as a statement saying that 
``the multiplicity of $f$ is at least $a$''. 
In this interpretation, the problem of resolution of singularities 
(cf.~\ref{0.2.1.1}) is, after a sequence 
of blowups and through transformations and at every point of the ambient 
space, to negate at least one 
statement in the collection. 
 
Observe in this interpretation that the following conditions 
naturally hold: 
$$ 
\left\{\begin{array}{cll} 
(\on{o}) 
&(f,0) \quad \forall f \in \cO_{W,P}, 
&(0,a) \quad \forall a \in \bZ 
\\ 
(\on{i}) 
&(f,a), (g,a) 
&\Longrightarrow (f + g,a) 
\\ 
&r \in \cO_{W,P}, (f,a) 
&\Longrightarrow (rf,a) 
\\ 
(\on{ii}) 
&(f,a), (h,b) 
&\Longrightarrow (fh,a + b) 
\\ 
(\on{iii}) 
&(f,a), b \leq a 
&\Longrightarrow (f,b). 
\end{array}\right. 
$$ 
 
Observe also that the problem of resolution of singularities stays 
unchanged, 
even if we add the statements derived from the given collection using the 
above conditions (implications). 
For example, starting from the given collection 
$\{(f,a)\mid f \in \cI_P\}$, the problem 
stays unchanged even if we consider the new collection 
$\{(f,n)\mid f\in \cI_P^{\lceil{n}/{a}\rceil}, 
n\in \bZ_{\geq 0}\}$.  Our philosophy is that it should be 
theoretically more desirable to consider 
the larger or largest collection of statements toward the problem of 
resolution of singularities. 
 
Accordingly we define an idealistic filtration, at a point $P \in W$, to be 
a subset $\bI \subset {\cO}_{W,P} 
\times \bR$ satisfying the following conditions: 
$$ 
\left\{\begin{array}{clll} 
(\on{o}) 
&\lefteqn{(f,0) \in \bI 
\ \forall f \in {\cO}_{W,P}, 
\quad 
(0,a) \in \bI 
\ \forall a\in \bR 
}\\ 
(\on{i}) 
&(f,a), (g,a) \in \bI 
&\Longrightarrow 
&(f + g,a)\in \bI 
\\ 
&r \in \cO_{W,P}, (f,a) \in \bI 
&\Longrightarrow 
&(rf,a)\in \bI \\ 
(\on{ii}) 
&(f,a), (h,b) \in \bI 
&\Longrightarrow 
&(fh,a + b)\in \bI \\ 
(\on{iii}) 
&(f,a) \in \bI, b \leq a 
&\Longrightarrow 
&(f,b)\in \bI. 
\end{array}\right. 
$$ 
Note that, as a consequence of conditions (o) and (iii), we have 
$$(f,a) \in \bI \quad\text{for\ any}\quad 
f \in \cO_{W,P},\ a \in\bR_{\leq0}.$$ 
 
We say an element $(f,a) \in \bI$ is at level $a$.  Note that we let 
the level 
vary in $\bR$.  Starting from the level varying in $\bZ$, we are 
naturally led to 
the situation where we let the level varying in the fractions ${\mathbb Q}$ 
when we start considering 
the condition (cf.~$\fR$-saturation) 
$$(\on{radical}) \quad 
(f^n,na) \in \bI,\ n \in \bZ_{> 0} 
\Longrightarrow (f,a) \in \bI,$$ 
and then to the situation where we let the level varying in $\bR$ when 
we start considering the 
condition of continuity 
$$(\on{continuity}) \quad 
(f,a_l) \in \bI\text{\ for\ a\ sequence\ } 
\{a_l\} \text{\ with\ } \lim_{l\rightarrow \infty} a_l = a 
\Longrightarrow (f,a) \in \bI.$$ 
 
Note that there is one more natural condition to consider related to the 
differential operators 
$$(\on{differential}) \quad 
(f,a) \in \bI, d \text{\ a\ differential\ operator of\ degree\ }t 
\Longrightarrow (d(f),a - t) \in \bI.$$ 
 
We remark that we do not include condition (radical), (continuity) or 
(differential) in the definition 
of an idealistic filtration, even though these conditions play crucial roles 
when we consider the 
radical and differential saturations of an idealistic filtration 
(cf.~\ref{0.2.3.2.3}).  We also introduce the notion of an 
idealistic filtration of r.f.g.~type (cf.~\ref{0.8}). 
\komoku{Distinguished features.}\label{0.2.3.2} 
Being framed in an extension of the classical notions, 
our program in the language of the idealistic filtration 
shares some common spirit with the existing 
approaches.  However, the following four features 
distinguish our program from them in a decisive way: 
\subkomoku{ 
Leading generator system as a collective 
substitute for a hypersurface of maximal contact.}\label{0.2.3.2.1} 
Given an idealistic filtration 
$\bI\subset \cO_{W,P}\times \bR$ 
at a point $P \in W$, we look at the graded ring 
of its leading terms $L(\bI) := 
\bigoplus_{n \in \bZ_{\geq 0}}L(\bI)_n$ where 
$L(\bI)_n = \{f \bmod{m_{W,P}^{n+1}} 
\mid (f,n) \in \bI, f \in\fm_{W,P}^n\}$. 
If we fix a regular system of parameters 
$(x_1, \dotsc, x_d)$ at $P$ and if we fix a natural isomorphism of 
$G =\bigoplus_{n \in \bZ_{\geq0}}m_{W,P}^n/m_{W,P}^{n+1}$ 
with the polynomial ring $k[x_1, \dotsc, x_d]$, the 
graded ring $L(\bI)$ can be 
considered as a graded $k$-subalgebra of 
$G = k[x_1, \dotsc, x_d]$. 
 
Now the fundamental observation is that (if 
the idealistic filtration is 
differentially saturated 
(cf.~$\fD$-saturation in \ref{0.2.3.2.3})) 
for a suitably chosen regular system 
of parameters, we can choose 
the generators of 
$L(\bI)$, as a graded $k$-subalgebra of $k[x_1, \dotsc, x_d]$, 
to be of the form 
$$\{x_i^{p^{e_i}}\mid e_i \in \bZ_{\geq 0}\}_{i \in 
I} \quad\text{for\ some}\quad I \subset \{1, 
\dotsc, d\}$$ 
when we are in positive characteristic $\on{char}(k) = p> 0$. 
We define a leading generator system of 
the idealistic filtration to be a set of elements 
$\{({h_i},p^{e_i})\}_{i\in I}\subset \bI$ 
whose leading terms give rise to the set of generators as above, i.e., 
$h_i \bmod{m_{W,P}^{p^{e_i}+1}} = x_i^{p^{e_i}}$ for $i \in I$.  We 
emphasize that the leading terms of the 
elements in the leading generator system lie in degrees 
$p^0, p^1, p^2, p^3, \dotsc$, and hence that the leading generator 
system may not form (a part of) a regular 
system of parameters when we are in positive characteristic 
$\on{char}(k)= p > 0$.  In the example by 
R.~Narasimhan, where there is no nonsingular hypersurface of maximal 
contact, there is no leading term of 
degree one in any leading generator system. 
When we are in characteristic zero $\on{char}(k) = 0$, 
in contrast, we can choose the generators of $L(\bI)$ to 
be concentrated all in degree one, i.e., of the form 
$$\{x_i\}_{i \in I} \quad\text{for\ some}\quad 
I \subset \{1, \dotsc,d\}.$$ 
Accordingly, we can take a leading generator system 
to be a set of elements 
\linebreak 
$\{(h_i,1)\}_{i \in I}\subset\bI$ with 
$h_i\bmod{m_{W,P}^2} = x_i$ for $i \in I$. 
If we look at the classical algorithm(s), then a 
hypersurface of maximal contact (locally at $P$) is 
given by $\{h_i = 0\}$ (for some $i\in I$). 
Since the leading term of $h_i$ is linear, it is guaranteed to 
define a {\it nonsingular} hypersurface. 
 
However, the case in positive characteristic and the case in characteristic 
zero should not be considered 
as two separate entities.  Rather, the case in characteristic zero should be 
considered as a special case of 
the uniform phenomenon: Traditionally we define the characteristic 
$\on{char}(k)$ to be the 
(non-negative) generator of the set of the annihilators of the unit ``$1$'' 
in the field $k$.  However, for 
the purpose of considering the problem of resolution of singularities, it is 
more natural to adopt the 
convention that the ``characteristic'' $p$ attached to the field $k$ is 
defined by 
$$p = \inf\{n \in \bZ_{> 0}\mid n \cdot 1 = 0 \in k\}.$$ 
In other words, we expect the behavior in characteristic zero to be similar 
to the one in positive 
characteristic with large $p$, and ultimately to lie at the limit when 
$p\rightarrow \infty$.  In this 
regard with the above convention, in characteristic zero, the (virtual) 
leading terms of the leading 
generator system in degrees $p^1 = p^2 = \dotsb = \infty$ are 
invisible (non-existent), while the 
actual leading terms are concentrated all in degree 
$\lim_{p \rightarrow\infty}p^0 = 1$. 
 
That is to say, we consider the notion of a hypersurface of maximal contact 
in characteristic zero to be a 
special case of the notion of a leading generator system, which is valid in 
all characteristics. 
Accordingly, we use the notion of a leading generator system as a collective 
substitute in positive 
characteristic for the notion of a hypersurface of maximal contact in 
characteristic zero in the process of 
constructing an algorithm according to our program. 
\subkomoku{Enlargement vs. restriction. 
{\it (Construction of the strand of invariants 
only through enlargements (modifications) of an 
idealistic filtration, and without using restriction to a 
hypersurface of maximal contact.)}}\label{0.2.3.2.2} 
At first sight, the introduction of the 
notion of a leading generator 
system does not seem to contribute toward overcoming the main source of 
troubles at all.  Recall (cf.~\ref{0.2.1.3}) that in the classical setting 
in characteristic zero the strand of 
invariants is constructed in 
such a way that a unit 
$(w,s)$ is added to the strand constructed so far every time we decrease the 
dimension by one, and then 
continue the construction by restricting ourselves to a hypersurface of 
maximal contact.  Nonsingularity of a 
hypersurface of maximal contact is absolutely crucial in order to continue 
the construction by restriction. 
Therefore, in the new setting in positive characteristic where we use a 
leading generator system, we seem to 
fail to construct the strand of invariants if any of the elements in the 
leading generator system defines a 
singular hypersurface.  However, in the construction of the strand of 
invariants in 
the new setting, we do not use any restriction but only use enlargements 
(modifications) of the idealistic 
filtration.  In fact, starting from a given idealistic filtration on a 
nonsingular variety $W$, we construct 
the triplet of invariants $(\sigma,\widetilde{\mu},s)$, where $\sigma$ 
reflects the degrees of the leading 
terms of a leading generator system, and $\widetilde{\mu}$ and $s$ are the 
weak-order (with respect to a 
leading generator system) and the invariant determined by the boundary, 
respectively, corresponding to the 
invariants $w$ and $s$ as before.  In the classical setting, after taking 
the corresponding companion 
modification and boundary modification, we take a hypersurface of maximal 
contact at this point and 
continue the process by taking the restriction to it.  In the new setting, 
however, after taking the 
companion modification and boundary modification, we consider a leading 
generator system of the newly 
modified idealistic filtration and continue the process.  In other words, in 
the new setting, we construct 
the strand of invariants in the following form 
$$\on{\it inv}_{\on{new}} = 
(\sigma,\widetilde{\mu},s)(\sigma,\widetilde{\mu},s) 
(\sigma,\widetilde{\mu},s)\dotsb,$$ 
and the construction is done only through enlargement 
keeing the ambient space $W$ intact, 
and hence the crucial nonsingularity intact. 
 
It is worthwhile noting that $\widetilde{\mu}$ is independent of the choice 
of a leading 
generator system, which is a priori needed for its definition, and hence is 
an invariant canonically 
attached to the idealistic filtration (if it is appropriately saturated 
(See \ref{0.2.3.2.3} below.)).  This 
implies that the strand of invariants 
$\on{\it inv}_{\on{new}}$ 
is also 
canonically determined globally. 
Therefore, we see that the center of each blowup in our algorithm, which is 
the maximum locus of the strand 
of invariants, is also canonically and globally defined, without the 
so-called Hironaka's trick needed in the classical setting 
(cf.~\ref{0.2.3.2.3} and \cite{MR2163383}).  

In Part II, we will define the two basic invariants denoted by 
$\sigma$ and $\widetilde{\mu}$ in the cotext of 
an idealistic filtration as above.  They form the building blocks for 
constructing the strand of invariants (together with invariant $s$ related 
to the boundary).  Some of their properties which 
are straightforward in characterisic zero, e.g., the upper semi-continuity, 
become highly non-trivial in positive 
characteristic and are also discussed in Part II.  
 
Discussion of the modifications is one of the main themes of Part III, where 
the classical notion of 
the companion modification and that of the boundary modification find their 
perfect analogs in the context 
of the enlargements of an idealistic filtration with respect to a leading 
generator system. 
\subkomoku{Saturations.}\label{0.2.3.2.3} 
It is important in our program to 
make a given idealistic 
filtration ``larger'' without changing the associated problem of resolution 
of singularities.  Ultimately, we 
would like to find the largest of all such (with respect to a certain fixed 
kind of operations ``$X$''), 
leading to the notion of the ($X$-)saturation.  Dealing with the saturated 
idealistic filtration, we expect 
to extract more intrinsic information toward a solution of the problem of 
resolution of singularities (e.g.~invariants 
which are independent of the choice of a leading generator system 
in the new setting, or the choice of a hypersurface of maximal 
contact in the classical setting).  The two key saturations in our 
program are the differential saturation (called the 
$\fD$-saturation for short, with respect to the operation 
of taking differentiations) and the radical saturation (called the 
$\fR$-saturation for short, with respect to the operation 
of taking the $n$-th roots (radicals)), the 
latter being equivalent to taking the integral closure (for an 
idealistic filtration of r.f.g.~type).  (The operation of 
taking the coefficient ideal and the operation of taking the 
``homogenization'' in the sense of \cite{MR2163383} share 
the same spirit with $\fD$-saturation.  In fact, we 
can obtain new formulas for the coefficient ideal and the 
homogenization as byproducts of the notion of the $\fD$-saturation 
of an idealistic filtration.  See \cite{AG0103120} for details. 
We also invite the reader to look at 
\cite{AG0508332}, which discusses several 
extensions of the idea of homogenization.)  At the center of our program 
sits the analysis of the interaction of these two 
saturations, leading to the notion of the bi-saturation (called the 
${\mathfrak B}$-saturation) and its explicit description as the 
$\fR\fD$-saturation.  Note that the notion of a leading 
generator system in \ref{0.2.3.2.1} is defined only through 
$\fD$-saturation, and the new nonsingularity principle 
in \ref{0.2.3.2.4} only through ${\mathfrak B}$-saturation. 
\subkomoku{New nonsingularity principle.}\label{0.2.3.2.4} 
There is another problem which comes along with 
using a leading generator system as a collective substitute for a 
hypersurface of maximal contact.  In the classical setting in 
characteristic zero, what 
guarantees the nonsingularity of the 
center is the nonsingularity of a hypersurface of maximal contact 
(cf.~\ref{0.2.1.3}).  In our new setting in positive 
characteristic, we no longer have this guarantee.  In fact, at the 
intermediate stage of the construction of 
the strand of invariants, the leading generator system may not be 
(a part of) a regular system of parameters and hence may define 
a singular subscheme.  We observe, however, that at the end of the 
construction of the strand of invariants the enlarged idealistic 
filtration takes such a special form that 
guaratees the corresponding leading generator system to be 
(a part of) a regular system of parameters. 
The maximum locus of the strand of the invariants, which we 
choose as the center, is defined by this leading 
generator system, and hence is nonsingular. 
We call this observation the new nonsingularity principle 
of the center. 
\komoku{Uniformity of our program in all characteristics.}\label{0.2.3.3} 
It should be emphasized 
that our program is not designed to come up with an esoteric strategy 
peculiar to the situation in positive 
characteristic, but rather intended to develop a uniform point of view 
toward the problem of resolution of 
singularities valid in all characteristics.  Part IV is devoted to letting 
this point of view manifest 
itself in the form of an algorithm, summarizing all the ingredients of the 
program. 
\end{subsection} 
\end{section} 
\begin{section}{Algorithm constructed according to the program.}\label{0.3} 
\begin{subsection}{Algorithm in characteristic zero.}\label{0.3.1} 
Aiming at uniformity, our program makes 
perfect sense and works just as well in characteristic zero, leading to a 
new algorithm slightly different 
from the existing ones.  We will demonstrate how the distinguished features 
of our program described in \ref{0.2.3.2} work in the new algorithm. 
\end{subsection} 
\begin{subsection}{Algorithm in positive characteristic; 
the remaining problem of termination.}\label{0.3.2} 
The algorithm in characteristic zero, now through uniformity, serves as a 
prototype toward establishing 
an algorithm in positive characteristic.  In fact, we can carry out almost 
all the procedures in positive 
characteristic, forming a perfect parallel to the case in characteristic 
zero, except for the problem of termination. 
\komoku{Termination.}\label{0.3.2.1} 
It is easy to see that in 
characteristic zero the invariants 
constituting the strand, constructed according to the program, have bounded 
denominators, and hence that the 
strand takes its value in the set satisfying the descending chain condition. 
Since the value of the strand 
strictly drops after each blowup, we conclude that the algorithm terminates 
after finitely many steps.  However, in positive 
characteristic, we can not exclude the 
possibility that the denominators may increase indefinitely as we carry out 
the processes (blowups) of the 
algorithm.  (In the 
unit $(\sigma,\widetilde{\mu},s)$ for the strand, the values of invariant 
$\sigma$ and $s$ are integral by 
definition.  Therefore, more specifically, 
the only issue is the boundedness 
of the denominators for the values of 
$\widetilde{\mu}$, which are fractional.) 
Therefore, we do not know at the 
moment if the algorithm 
terminates after finitely many steps. 
 
The problem of termination remains as the 
only missing piece in our quest to complete an algorithm for resolution of 
singularities in positive characteristic according 
to the program. 
\end{subsection} 
\end{section} 
\begin{section}{Assumption on the base field.}\label{0.4} 
We carry out our entire program assuming that the base field $k$ is 
algebraically closed field of characteristic 
$\on{char}(k) = p \geq 0$. 
 
Our definition of a leading generator system, 
the key notion of the program, at a closed point $P \in W$ 
where $W$ is a variety of 
dimension $d$ smooth over $k$, needs the assumption of the base field being 
algebraically closed, since we use the fact 
$\cO_{W,P}/m_{W,P} \cong k$ and the natural isomorphism 
$G = \bigoplus_{n \geq 0}m_{W,P}^{n+1}/m_{W,P}^n 
\cong k[x_1, \dotsc, x_d]$ 
with respect to a fixed regular system of parameters 
$(x_1, \dotsc, x_d)$, as well as 
the fact that we can take the $p$-th root of any element within $k$ 
(when $\on{char}(k) = p > 0$). 
We briefly mention below what happens if we loosen the assumption on the 
base field. 
\begin{subsection}{Perfect case.}\label{0.4.1} Suppose that the base field 
$k$ is perfect, but not necessarily 
algebraically closed.  Upon completion, the algorithm constructed according 
to the program 
should be equivariant under any group action (cf. Part IV).  Therefore, as 
long as the 
base field $k$ is perfect, we see 
that the algorithm established over its algebraic closure 
$\overline{k}$ descends to the one over the original base field $k$, 
utilizing 
the equivariance under the action 
of the Galois group 
$\on{Gal}(\overline{k}/k)$. 
\end{subsection} 
\begin{subsection}{Non-perfect case.}\label{0.4.2} 
Over a non-perfect field $k$, 
we even have to start distinguishing the notion of being regular 
and that of being smooth over $k$.  The discussions, 
including the one on how we may try to reduce the 
non-perfect case to the perfect case 
using the Lefschetz Principle type 
argument, will be given in Part IV. 
\end{subsection} 
\end{section} 
\begin{section}{Other methods and approaches.}\label{0.5} 
We only mention a few of the other methods and approaches than the 
algorithmic approach we follow toward the problem 
of resolution of singularities in positive characteristic.  We refer the 
reader to 
\cite{MR0389901}\cite{MR1395176}\cite{MR1748614} for 
a more detailed account. 
 
Resolution of singularities for curves is a classical result, with many of 
its ideas and methods leading to the 
higher dimensional cases even to this day. 
 
Among several results for surfaces, the most general one seems to be given 
by \cite{MR0276239}\cite{MR0491722}, which 
establish resolution of singularities of an arbitrary excellent scheme in 
dimension $2$. 
 
It is \cite{MR0002864} that initiated the strategy to establish local 
uniformizations first, with the theory of 
valuations as its central tool, and then by patching them to establish 
resolution of singularities globally.  The 
theory of local uniformization has been further developed by many people 
\cite{MR0217069}\cite{MR1748622}\cite{Kuhlmann97}\cite{MR1748629}. 
We should mention the approaches by 
\cite{MR2018565}\cite{S_BIRS} toward local uniformization 
in higher dimensions. 
 
Jung's idea of taking the (generic) projection provides many useful 
approaches toward the problem of resolution of 
singularities.  \cite{MR0217069} uses the method of Albanese projecting from 
a singular point, combined with the 
theory of local uniformization, to resolve singularities of a threefold $X$ 
when $\on{char}(k)$ is not greater than $(\dim X)! = 6$. 
A simplified proof has been recently given by \cite{AG0606530}, which 
also discusses the potential and 
problems if one tries to extend the method to higher dimensions.  There are 
attempts to study the problem in the 
remaining characteristic $\on{char}(k) = 2,3,5$ by 
\cite{MR907903}\cite{MR1395176}\cite{C_BIRS}\cite{P_BIRS} 
in dimension $3$. 
 
Without any restriction on the dimension 
of a variety or on the base field $k$, the most remarkable development in 
the vicinity of the 
problem of resolution of singularities is arguably the method of alteration 
initiated by de Jong \cite{MR1423020}.  Given a variety 
$X$, it constructs a proper and {\it generically finite} morphism 
$f\colon Y \rightarrow X$ 
from a regular variety $Y$.  (In characteristic zero, one can 
refine the method of alteration to realize $f$ as a birational map.  
See \cite{MR1487237}\cite{MR1397679}\cite{MR1714830} for details.) 
The structure of $f$ is rather obscure, 
though its existence follows nicely and 
simply by regarding 
$X$ as a family of curves fibered over a variety of dimension one less and 
hence by paving a way to apply induction. 
The method of alteration even works in mixed characteristics or 
with integral schemes over $\bZ$, and hence 
it allows a wide range of applications for arithmetic purposes. 
\end{section} 
\begin{section}{Origin of our program.}\label{0.6} 
This series of papers is a joint work of H. Kawanoue and 
K. Matsuki as a whole.  However, the program forming the backbone of the 
series was conceived in its entirety 
by the first author toward his Ph.D.~thesis, and revealed to the second 
author in the summer of 2003 at a 
private seminar held at Purdue University as a blueprint toward constructing 
an algorithm for resolution of 
singularities in positive characteristic.  As such all the essential ideas 
are due to the first author. 
Accordingly it should be called the Kawanoue program, which we use as the 
subtitle starting from Part II.  Only the name of the 
first author appears on the cover of Part I, which represents the main 
portion of his Ph.D. thesis. 
 
The only contribution of 
the second author was to help the first author and jointly bring these ideas 
together converging into a coherent algorithm. 
\end{section} 
\begin{section} 
{Acknowledgement.}\label{0.7} 
Our entire project could only be 
possible through the guidance 
and encouragement of Professor Shigefumi Mori both at the personal level and 
in the mathematical context.  He 
not only shared his insight generously with us, but also on several 
occasions in the development of the 
Kawanoue program showed us directly some key arguments to bring us forward. 
Professor Masaki Kashiwara 
also gave us an invaluable and enthusiastic support, without which the 
project would have dissipated into the 
air. 
 
We thank Professors Edward Bierstone, Pierre Milman, Orlando Villamayor, and 
Herwig Hauser, from whom we 
learned most on the subject of resolution of singularities, where the 
tutoring was given in the form of publications and 
personal correspondences.  Only through their teaching, we started 
understanding the greatest ideas of \cite{MR0199184}.  Many of the 
ideas of our project, therefore, find their origins in \cite{MR0199184} 
as well as in the papers of our teachers cited above.  Our 
indebtedness to Professor Heisuke Hironaka, whose influence 
was decisive for us to enter the subject, is immeasurable. 
 
It is a pleasure to acknowledge the helpful comments and suggestions we 
received from Professors Donu Arapura, Johan de Jong, 
Joseph Lipman, Tsuong-Tsieng Moh, Tadao Oda, 
Bernd Ulrich, Jaros{\l}aw W{\l}odarczyk. 
 
Special thanks go to Hidehisa Alikawa, Takeshi Nozawa, and Masahiko 
Yoshinaga, who were both good friends 
and patient listeners in Room 120 for the graduate students of Research 
Institute for Mathematical Sciences 
in Kyoto at the dawn of the Kawanoue program. 
\end{section} 
\begin{section}{Outline of Part I.}\label{0.8} 
Following the itemized table 
of contents at the beginning, we 
describe the outline of the structure of Part I below. 
 
At the end of the introduction in Chapter 0, we give a brief description of 
the preliminaries to read 
Part I and the subsequent series of papers.  In Chapter 1, we recall some 
basic facts on the differential 
operators, especially those in positive characteristic.  Both in the 
description of the preliminaries and in Chapter 1, our 
purpose is not to exhaustively cover all the material, but only to minimally 
summarize what is needed to present our program 
and to fix our notation.  For example, an elementary characterization, in 
terms of the differential operators, of an ideal 
generated by the $p$-th power elements in characteristic 
$p = \on{char}(k)> 0$ is included only due to the lack of an 
appropriate reference.  We should emphasize here that the use of the 
logarithmic differential operators is indispensable 
in our setting in the language of the idealistic filtration 
(See Remark \ref{1.2.2.3}). 
 
Chapter 2 is devoted to establishing the notion of an idealistic filtration, 
and its fundamental properties. 
The most important ingredient of Chapter 2 is the analysis of the 
$\fD$-saturation and $\fR$-saturation 
and that of their interaction.  In our algorithm, given an 
idealistic filtration, we always 
look for its bi-saturation, called the 
${\mathfrak B}$-saturation, which is both 
$\fD$-saturated and $\fR$-saturated 
and which is minimal among such containing the original 
idealistic filtration.  The existence of 
the ${\mathfrak B}$-saturation is theoretically clear. 
However, we do not know a priori whether we can reach 
the ${\mathfrak B}$-saturation by a repetition of 
$\fD$-saturations and $\fR$-saturations 
starting from the given idealistic filtration, even after 
infinitely many times.  The main result here is that the 
${\mathfrak B}$-saturation is actually realized if we 
take the $\fD$-saturation and then 
$\fR$-saturation of the given one, each just once in this 
order.  In our algorithm, we do not deal with an arbitrary idealistic 
filtration, but only with 
those which are generated by finitely many elements with rational levels. 
We say they are of r.f.g.~type 
(short for ``rationally and finitely generated'').  It is then a natural and 
crucial question if the property 
of being of r.f.g.~type is stable under $\fD$-saturation and 
$\fR$-saturation.  We find 
somewhat unexpectedly that the argument of M.~Nagata 
(cf.~\cite{MR0089836}), 
which was originally developed to answer some 
questions posed by P.~Samuel regarding 
the asymptotic behavior of ideals, is 
tailor-made to establish the 
stability under 
$\fR$-saturation (while the stability under 
$\fD$-saturation is elementary). 
 
In Chapter 3, through the analysis of the leading terms of an idealistic 
filtration (which is 
$\fD$-saturated), we define the notion of 
a leading generator system, which, as discussed in \ref{0.2.3.2.1}, 
plays the role of a collective substitute for the notion of a hypersurface 
of maximal contact. 
 
Chapter 4 is the culmination of Part I, establishing the new nonsingularity 
principle of the center for an 
idealistic filtration which is ${\mathfrak B}$-saturated. 
Its proof is given via the three somewhat technical 
but important lemmas, which we will use again later 
in the series of papers. 
 
Our theory in Part 
I is mainly local, dealing almost exclusively with an idealistic filtration 
over the local ring of a closed 
point on a nonsingular ambient variety.  The global theory toward 
constructing an algorithm will be 
discussed in the subsequent papers. 
 
Of course the main purpose of Part I is to establish the foundation of our 
program toward constructing an 
algorithm for resolution of singularities.  However, we believe that the 
results on the idealistic filtration 
we discuss here in Part I, notably the analysis leading to the explicit 
description of the ${\mathfrak B}$-saturation, 
stability of r.f.g.~type, and the nonsingularity principle, 
are of interest on their own in the subject of the ideal theory in 
commutative algebra. 
 
\bigskip
 
This finishes the discussion of the outline of Part I. 
\end{section} 
\begin{section}{Preliminaries.}\label{0.9} 
We summarize a few of the 
preliminaries in order to read Part I and 
the subsequent series of papers. 
\begin{subsection}{The language of schemes.}\label{0.9.1} 
Our entire argument is carried out in the language of schemes.  For 
example, a variety is an integral separated scheme of finite type over $k$. 
Accordingly, when we say ``points'', we refer to 
the scheme-theoretic points and do not confine ourselves to the closed 
points, which correspond to the geometric ones in the 
classical setting.  Thus the invariants that we construct will be defined 
over all the scheme-theoretic points, and not 
confined to the closed points.  However, some of the key notions of our 
program, notably that of a leading generator system, 
are only defined at the level of the closed points, and the values of the 
invariants over the non-closed points are given only 
indirectly through their upper or lower semi-continuity. 
 
Our program is not conceived in the language of schemes originally.  Rather, 
it has its origin in the concrete analysis 
and computation in terms of the coordinates at the closed points.  As such, 
it can be applied to many other ``spaces'' than 
algebraic varieties over $k$, where the same analysis and computation can be 
applied to the coordinates at its closed points. 
The task of presenting a set of axiomatic conditions for the Kawanoue 
program to function, and that of listing explicitly the 
spaces within its applicability will be dealt with elsewhere. 
\end{subsection} 
\begin{subsection}{Basic facts from commutative algebra.}\label{0.9.2} 
For the basic facts in commutative algebra, we try to 
use \cite{MR879273} as the main source of reference. 
\end{subsection} 
\begin{subsection}{Multi-index notation.}\label{0.9.3} 
When we have the multivariables, either as the indeterminates 
in the polynomial ring or as a regular system of parameters, we 
often use the following multi-index notations: 
\par\bigskip
$ 
\left\{\begin{array}{llll} 
\displaystyle{ 
X 
}&\displaystyle{ 
= (x_1, \dotsc , x_d), 
}& 
\qquad 
\displaystyle{ 
I 
}&\displaystyle{ 
= (i_1, \dotsc , i_d) \in \bZ_{\geq0}^d, 
}\\ 
\displaystyle{ 
|I| 
}&\displaystyle{ 
= \sum_{\alpha = 1}^di_{\alpha}, 
}& 
\qquad 
\displaystyle{ 
X^I 
}&\displaystyle{ 
= \prod_{\alpha = 1}^d x_{\alpha}^{i_{\alpha}}, 
}\\ \displaystyle{ 
\binom{I}{J} 
}&\displaystyle{ 
= \prod_{\alpha = 1}^d\binom{i_{\alpha}}{j_{\alpha}} 
}&\text{ for} 
\quad 
\displaystyle{J} 
&= 
\displaystyle{ 
(j_1, \dotsc, j_d) \in \bZ_{\geq0}^d 
}\\ 
&\lefteqn{\displaystyle{ 
\phantom{=}\text{where}\quad 
\binom{i}{j} 
= \frac{i !}{(i-j)!j!} \in \bZ_{\geq0}} 
\text{\ denotes\ the\ binomial\ coefficient},} 
\\ 
\lefteqn{\displaystyle{ 
\phantom{=} 
\text{(We also use the convention that whenever } 
i_{\alpha} < j_{\alpha} 
\text{ we set } 
\binom{i_{\alpha}}{j_{\alpha}}= 0).}} 
\\ 
\displaystyle{ 
\partial_{X^J} 
}& 
\lefteqn{ 
\displaystyle{ 
=\frac{\partial^{|J|}}{\partial_{x_1}^{i_1} \dotsm 
\partial_{x_d}^{i_d}} 
\quad\text{(expressed by }\partial_{J}}\text{ for short).} 
}\\ 
{\mathbf e}_{\alpha} 
&
\lefteqn{ 
= (0,\dotsc,\overset{\underset{\vee}{\alpha}}{1},\dotsc, 0).} 
\end{array}\right.$ 
\end{subsection} 
\end{section} 
\end{chapter} 
\begin{chapter}{Basics on differential operators} 
The purpose of this chapter is to give 
a brief account of the differential operators, 
which play a key role in the Kawanoue program. 
 
We would like to mention that it is through reading the papers 
\cite{MR0269658}\cite{MR0472824} that our attention 
was first brought to the importance of the higher order 
differential operators in the context of the 
problem of resolution of singularities in positive 
characteristic. 
 
Our main reference is EGA IV \S 16 \cite{MR0238860}, where all 
that we need, especially the properties of the higher order differential 
operators of 
Hasse-Schmidt type in positive characteristic, and much more, is beautifully 
presented.  We only try to extract some basic facts 
and discuss them in the form that suits our limited purposes. 
\begin{section}{Definitions and first properties}\label{1.1} 
\begin{subsection}{Definitions.}\label{1.1.1} 
 
Recall that the base field $k$ is 
assumed to be an algebraically closed field of 
$\on{char}(k) \geq 0$. 
\begin{defn}\label{1.1.1.1} 
Let $R$ be a $k$-algebra.  We use the 
following notation: 
\begin{eqnarray*} 
&\mu\colon R \otimes_kR \rightarrow 
R \quad \text{the\ multiplication\ map}, 
\qquad 
I := \on{ker}(\mu) \quad \text{the\ kernel\ of\ }\mu, \\ 
&P_R^n =R \otimes_kR/I^{n+1}, 
\quad 
q_n \colon R \rightarrow R \otimes_kR 
\rightarrow P_R^n \quad \text{for\ }n \in \bZ_{\geq 0} 
\end{eqnarray*} 
where $q_n$ is the composition of the map to the second factor with the 
projection, i.e., 
$$q_n(r) = (1 \otimes r \bmod{I^{n+1}}) 
\quad\text{for}\quad r \in R.$$ 
A differential operator $d$ of degree $\leq n$ on $R$ (over $k$) for 
$n \in\bZ_{\geq 0}$ is a map $d\colon R \rightarrow R$ of the form 
$$d = u \circ q_n \quad\text{with}\quad u \in \on{Hom}_R(P^n_R,R).$$ 
(We note that the $R$-module structure on $P^n_R$ is 
inherited from the $R$-module structure on $R 
\otimes_kR$ given by the multiplication on the first factor.) 
 
We denote the set of differential operators of degree $\leq n$ on $R$ by 
$\dif^n_R$, i.e., 
$$\dif^n_R := \left\{d = u \circ q_n\mid u \in 
\on{Hom}_R(P^n_R,R)\right\}.$$ 
(Note that $\dif^n_R$ inherits the $R$-module structure from the one 
on $\on{Hom}_R(P^n_R,R)$.) 
 
We call $\dif_R = \bigcup_{n = 0}^{\infty}\dif_R^n$ 
(cf.~Lemma \ref{1.1.2.1}) the set of the differential operators on 
$R$ (over $k$). 
 
For a subset $T \subset R$, we also use the following notation 
$$\dif_R^n(T) = (\left\{d(r)\mid d \in 
\dif_R^n, r \in T\right\}).$$ 
\end{defn} 
\end{subsection} 
\begin{subsection}{First properties.}\label{1.1.2} 
\begin{lem}\label{1.1.2.1} 
Let the situation and notation be the same as in 
Definition \ref{1.1.1.1}. 
\item[(1)] 
Let $d$ be a $k$-linear map $d\colon R \rightarrow R$.  Then 
$d$ is a differential operator of degree 
$\leq n$, i.e., $d \in \dif^n_R$ if and only if $d$ satisfies the 
Leibnitz rule of 
degree $n$: 
$$\sum_{T \subset 
S_{n+1}}(-1)^{|T|}\left(\prod_{s 
\in S_{n+1} \setminus T}r_s\right)d\left(\prod_{s \in T}r_s\right) = 0$$ 
where $S_{n+1} = \{1, 2, \ldots, 
n, n+1\}$ and $r_s \in R \text{\ for\ }s \in S_{n+1}$. 
\item[(2)] 
The natural map 
$$\phi_R\colon\on{Hom}_R(P^n_R,R) \rightarrow \dif^n_R,$$ 
given by 
$d = \phi_R(u) = u \circ q_n$ for $u \in \on{Hom}_R(P^n_R,R)$, is 
bijective (and actually an isomorphism between $R$-modules). 
\item[(3)] 
If $R$ is finitely generated as an algebra over $k$, then $P_R^n$ is 
finitely generated as an $R$-module, and 
so is $\on{Hom}_R(P^n_R,R) 
\overset{\sim}{\rightarrow} \dif^n_R$. 
\item[(4)] 
Let $R'$ be the localization $R_S$ of $R$ with respect to a 
multiplicative set $S \subset R$ or the completion $\widehat{R}$ of 
$R$ with respect to a maximal ideal $\fm \subset R$.  We define the 
map 
$\dif^n_R 
\rightarrow 
\dif^n_{R'}$ so that the following diagram commutes 
$$\CD 
\on{Hom}_R(P^n_R,R) @>\phi_R>> @. \dif^n_R \\ 
\downarrow @. @. \downarrow\\ 
\on{Hom}_R(P^n_R,R) \otimes_RR' @>\phi_R \otimes_RR'>> @. 
\dif^n_R \otimes_RR'\\ 
\downarrow @. @. \\ 
\on{Hom}_{R'}(P^n_R \otimes_RR', R \otimes_RR') @.@. @VVV\\ 
\| @.@.\\ 
\on{Hom}_{R'}(P^n_{R'},R') @>\phi_{R'}>> @. \dif^n_{R'}, \\ 
\endCD$$ 
where the vertical arrows are the natural maps. 
 
Consequently, the bijections are compatible with localization and 
completion. 
 
Moreover, if $R$ is essentially of finite type over $k$, then the second 
vertical arrow on the left is an isomorphism, and hence so 
is the second vertical arrow on the right. 
\item[(5)] 
Let $d \in \dif^n_R$ be a differential operator of degree $\leq 
n$ on $R$. 
Then $d$ is a differential operator of degree $\leq m$ for any $n \leq m$. 
That is to say, 
$$\dif^n_R \subset \dif^m_R \text{\ for\ }n \leq m.$$ 
With respect to these inclusions, $\{\dif^n_R\}_{n \in 
\bZ_{\geq 0}}$ forms a projective system. 
\item[(6)] 
Let $d \in \dif^n_R$ be a differential operator of degree $\leq 
n$ on $R$, and 
$d' \in \dif^{n'}_R$ be a differential operator of degree $\leq n'$ 
on $R$.  Then 
the composition $d \circ d'$ is a differential operator of degree $\leq n + 
n'$ on $R$, 
i.e., $d \circ d' \in \dif^{n + n'}_R$. 
\item[(7)] 
Let $R$ be an algebra essentially of finite type over $k$, 
$I \subset R$ an ideal, and let $R'$ be as in $\on{(4)}$. 
Then we have $$\dif_R^n(I)R' = \dif_{R'}^n(IR').$$ 
\end{lem} 
\begin{proof} 
\item[(1)] 
We refer the reader to Proposition (16.8.8) in 
EGA IV \S 16 \cite{MR0238860} for a proof. 
\item[(2)] 
The isomorphism $\phi_R$ is the one mentioned in 
(16.8.3.1) in EGA IV \S 16 \cite{MR0238860}. 
\item[(3)] 
Suppose $R$ is finitely generated as an algebra over $k$.  Let 
$X =\{x_1,\ldots,x_t\}$ be a set of generators for $R$ over 
$k$.  We see that $P_R^n$ is generated by $\{q_n(X^I) 
\mid I \in \bZ_{\geq 0}^t\}$ as an $R$-module 
(cf.~the first note in Definition \ref{1.1.1.1}). 
We also see, by the relation $\prod_{s\in S_{n+1}} 
(1 \otimes r_s - r_s \otimes 1) = 0$ in 
$P_R^n$, that $q_n(X^I)$ for any $I \in \bZ_{\geq 0}^t$ belongs to the 
$R$-span of $\{q_n(X^I)\mid I \in 
\bZ_{\geq 0}^t, |I| \leq n\}$.  Therefore, we conclude that 
$P^n_R$ is finitely generated as an $R$-module and hence that so is 
$\on{Hom}_R(P^n_R,R) 
\overset{\sim}{\rightarrow} \dif^n_R$. 
\item[(4)] 
Compatibility of the bijections with localization and 
completion follows immediately from the definitions and 
from the fact that $P^n_R \otimes_RR' = P^n_{R'}$. 
 
In order to verify the ``Moreover'' part, it suffices to show the assertion 
assuming that $R$ is finitely generated as an algebra 
over 
$k$.  Then since the extension $R 
\rightarrow R'$ is flat and since $P_R^n$ is finitely generated as an 
$R$-module by (3), the second vertical arrow on the left is an 
isomorphism, and hence so is the second vertical arrow on the right. 
\item[(5)] 
The natural surjection $P^m_R = (R 
\otimes_kR)/I^{m+1}\twoheadrightarrow P^n_R = (R 
\otimes_kR)/I^{n+1}$ for $n \leq m$ induces the 
injection $\on{Hom}_R(P^n_R,R) \hookrightarrow 
\on{Hom}_R(P^{n+1}_R,R)$ and hence the 
inclusion $\dif^n_R \subset \dif^m_R$.  It is clear that 
$\{\dif^n_R\}_{n \in \bZ_{\geq 0}}$ forms a 
projective system with respect to these inclusions. 
\item[(6)] 
We refer the reader to Proposition (16.8.9) in 
EGA IV \S 16 \cite{MR0238860}. 
\item[(7)] 
When $R' = \widehat{R}$, the equality $\dif_R^n(I)R' = 
\dif_{R'}^n(IR')$ follows from the ``Moreover'' 
part of (4) and from the fact that the differential operators 
are continuous with respect the $\fm$-adic topology (the 
latter being a consequence of the Leibnitz rule). 
 
Thus we give a proof of the equality only when 
$R' = R_S$ in the following. 
 
Since the inclusion $\dif_R^n(I)R_S \subset 
\dif_{R_S}^n(IR_S)$ follows easily from the 
``Moreover'' part 
of (4), we have only to show the opposite inclusion 
$$\dif_R^n(I)R_S \supset 
\dif_{R_S}^n(IR_S).$$ 
Take $f = s^{-1}r \in IR_S$ with $r \in R, s \in S$, and take $d \in 
\dif_{R_S}^n$.  We want to show $d(f) \in 
\dif_R^n(I)R_S$.  Set $r_1 = \cdots = r_n = s, r_{n+1} = f$. 
Applying the Leibnitz rule of degree $n$ 
for $d 
\in \dif^n_{R_S}$, we have 
$$- s^nd(f) + \sum_{\{n+1\} \subsetneqq T} 
(-1)^{|T|}s^{n+1-|T|}d\left(s^{|T|-2}r\right) + \sum_{n+1 \not\in 
T}(-1)^{|T|}fs^{n - 
|T|}d\left(s^{|T|}\right) = 0$$ 
where the first term in the left hand side corresponds to the range $T = 
\{n+1\}$.  Since $d \in \dif_{R_S}^n 
= \dif_R^nR_S$ by the ``Moreover'' 
part of (4), the second term and the third term of the left hand side belong 
to $\dif_R^n(I)R_S$.  This implies $d(f) \in 
\dif_R^n(I)R_S$. 
 
This completes the proof for Lemma \ref{1.1.2.1}. 
\end{proof} 
\begin{cor}\label{1.1.2.2} 
Let $X$ be a variety over $k$.  Then there 
exists a coherent sheaf 
$\on{{\mathcal H}\!{\it om}}_{\cO_X}({\mathcal P}^n_X,\cO_X) 
\overset{\sim}{\rightarrow} 
{\mathcal D\!\on{\it iff}}^n_X$ 
of the differential operators of degree 
$\leq n$ for 
$n 
\in \bZ_{\geq 0}$ such that for any affine open subset $U = 
\on{Spec}R \subset X$ 
we have 
$$\on{{\mathcal H}\!{\it om}}_{\cO_X}({\mathcal P}^n_X,\cO_X)(U) = 
\on{Hom}_R(P^n_R,R) \overset{\sim}{\rightarrow} \dif^n_R = 
{\mathcal D\!\on{\it iff}}^n_X(U)$$ 
and that for any point $x \in X$ we have a description of the stalk as 
$$\left\{\on{{\mathcal H}\!{\it om}}_{\cO_X}({\mathcal P}^n_X, 
\cO_X)\right\}_x = 
\on{Hom}_{\cO_{X,x}}(P^n_{\cO_{X,x}},\cO_{X,x}) 
\overset{\sim}{\rightarrow} 
\dif^n_{\cO_{X,x}} = \left\{ 
{\mathcal D\!\on{\it iff}}^n_X\right\}_x.$$ 
Moreover, for any closed point $x \in X$ we have a description 
of the completion of the stalk as 
$$\CD 
\left\{\on{{\mathcal H}\!{\it om}}_{\cO_X}({\mathcal P}^n_X, 
\cO_X)\right\}_x \otimes_{\cO_{X,x}} \widehat{\cO_{X,x}} \ 
@. @. \ \left\{{\mathcal D\!\on{\it iff}}^n_X\right\}_x 
\otimes_{\cO_{X,x}} \widehat{\cO_{X,x}} \\ 
\| @.@. \| \\ 
\on{Hom}_{\widehat{\cO_{X,x}}} 
(P^n_{\widehat{\cO_{X,x}}},\widehat{\cO_{X,x}}) @. 
\overset{\sim}{\rightarrow} @. 
\dif^n_{\widehat{\cO_{X,x}}}\\ 
\endCD$$ 
\end{cor} 
\begin{proof} 
This follows immediately from Lemma \ref{1.1.2.1}. 
\end{proof} 
\end{subsection} 
\end{section} 
\begin{section}{Basic properties of differential operators 
on a variety smooth over $k$.}\label{1.2} 
\markright{\ref{1.2}.\ 
BASIC PROPERTIES OF DIFFERENTIAL OPERATORS 
ON A SMOOTH VARIETY} 
The purpose of this section is to discuss some basic properties of 
differential operators on a variety $W$ smooth 
over $k$. 
 
Accordingly, we denote by $R$ the coordinate ring of an affine open subset 
$\on{Spec}R\subset W$, or its localization by some multiplicative set. 
\begin{subsection}{Explicit description of differential operators 
with respect to a regular system of parameters.}\label{1.2.1} 
\begin{defn}\label {1.2.1.1} 
We say $(x_1, \ldots, x_d)$ with $d = \dim W$ is 
a regular system of parameters for $R$ if $\{dx_{\alpha} 
= (1\otimes x_{\alpha} - x_{\alpha} \otimes 1 \bmod{I}) 
\mid \alpha = 1, \ldots , d\}$ forms a basis for the module 
of differentials $\Omega^1_{R/k}$ as an $R$-module, 
i.e., 
$$\Omega^1_{R/k} = (R \otimes_kR)/I 
= \bigoplus_{\alpha = 1}^dRdx_{\alpha} 
\cong R^d,$$ 
where $I \subset R 
\otimes_kR$ is the kernel of the multiplication map 
$\mu\colon R \otimes_kR \rightarrow R$ 
(cf.~Definition \ref{1.1.1.1}). 
 
(Note that in the case where $R$ is the local 
ring associated to a closed point $P \in W$ such a regular system of 
parameters always exists, and that in 
the case where $R$ represents the coordinate ring of an affine open subset 
$\on{Spec}R \subset W$ such a 
regular system of parameters exists by ``shrinking'' 
$\on{Spec}R$ if necessary.) 
\end{defn} 
\begin{lem}\label{1.2.1.2} 
Suppose we have a regular system of 
parameters $(x_1, \ldots, x_d)$ for $R$ with 
$d =\dim W$.  Then we have the following: 
\item[(1)] 
We have a family of maps $\{\partial_{X^J}\colon R \rightarrow R 
\mid J \in \bZ_{\geq 0}^d\}$ such that 
\begin{enumerate} 
\def\labelenumi{(\roman{enumi})} 
\item 
$\partial_{X^J}(X^I) = \binom{I}{J}X^{I - J} 
\text{\ for\ any\ }I \in \bZ_{\geq 0}^d$, and that 
\item 
$\{\partial_{X^J} \mid |J| \leq n\}$ 
forms a basis of $\dif_R^n$ for any $n \in 
\bZ_{\geq 0}$, i.e., 
$$\dif^n_R = \bigoplus_{|J| \leq n}R\partial_{X^J} \cong 
R^{\binom{n+d}{n}}.$$ 
\end{enumerate} 
\item[(2)] 
Let $\widehat{R}$ be the completion of $R$ with respect to a maximal 
ideal $\fm$ (corresponding 
to a closed point $P \in W$).  Then the $\widehat{R}$-module 
$\dif^n_{\widehat{R}} \overset{\sim}{\rightarrow} 
\dif_R^n \otimes_R\widehat{R}$ is free of rank 
$\binom{n+d}{n}$, having a basis 
$\{\partial_{X^J}\mid |J| 
\leq n\}$ of the differential operators of degree 
$\leq n$.  The differential operators are continuous with respect to the 
$\fm$-adic topology. 
 
Set $y_i = x_i - \alpha_i$ for $1 \leq i \leq d$, where $\alpha_i \in k$, so 
that 
$Y = (y_1, \ldots, y_d)$ is a regular system of parameters for $R_\fm$. 
Then for any $f = \sum c_IY^I \in k[[y_1, \ldots, y_d]] = 
\widehat{R}$, we have 
$$\partial_J(f) = \partial_J\left(\sum c_IY^I\right) = 
\sum c_I\partial_J(Y^I) = \sum c_I\binom{I}{J}Y^{I - J},$$ 
where $\partial_J$ is the abbreviated notation for $\partial_{X^J}$. 
\item[(3)] 
We have the generalized product rule 
$$\partial_J(fg) = \sum_{K + L = J}\partial_K(f)\partial_L(g) 
\quad\text{\ for\ }f,g \in R \ (\text{or\ }\widehat{R}).$$ 
\end{lem} 
\begin{proof} 
\item[(1)] 
We refer the reader to Theorem (16.11.2) in EGA IV 
[Gro67]. 
\item[(2)] 
Observe that a differential operator (of degree $\leq n$) is 
continuous with respect to the 
$\fm$-adic topology, a fact which easily follows, e.g., from the 
Leibnitz rule (of degree $n$).  Note that 
$\partial_{Y^J}(f) = \partial_{X^J}(f)$ for any $J \in \bZ_{\geq 0}^d$ 
and $f \in \widehat{R}$ by definition of $Y = (y_1, 
\ldots, y_d)$.  The rest is a direct consequence of (1). 
\item[(3)] 
In order to check the generalized product rule, it suffices to 
check it for the localization $R_\fm$ for any maximal ideals of $R$. 
In order to check it for the localization $R_\fm$, it suffices to 
check it for its completion $\widehat{R}$ with respect to $\fm$. 
 
By choosing a regular system of parameters $Y = (y_1, \ldots, y_d)$ for 
$R_\fm$ as in (2), we can identify $\widehat{R}$ with 
the power series ring $k[[y_1, \ldots, y_d]]$. 
Thus we have only to check (3) for the power series ring 
$k[[y_1, \ldots, y_d]]$.  By (2), it is also clear that 
we have only to check it for the case of one variable, 
i.e., $d = 1$ with 
$y_1 = y$ and that we may even assume $f$ and $g$ are powers of $y$, 
i.e., $f = y^a$ and $g = y^b$.  Then we have 
\begin{eqnarray*} 
\partial_{X^J}(fg) &=& \partial_{Y^J}(fg) = \partial_{y^n}(y^ay^b) = 
\partial_{y^n}(y^{a + b})\\ 
&=& \binom{a + b}{n}y^{a + b - n} = \left(\sum_{l+m = 
n}\binom{a}{l}\binom{b}{m}\right)y^{a+b-n}\\ 
&=& \sum_{l + m = 
n}\binom{a}{l}x^{a-l}\binom{b}{m}x^{b-m} = \sum_{K + L = 
J}\partial_K(f)\partial_L(g), 
\end{eqnarray*} 
which verifies the generalized product rule. 
 
This completes the proof of Lemma \ref{1.2.1.2}. 
 
\end{proof} 
\begin{rem}\label{1.2.1.3} 
\item[(1)] 
It is easy to see that we have a relation 
$$(\partial_{x_1})^{j_1} \circ (\partial_{x_2})^{j_2} 
\circ\cdots\circ 
(\partial_{x_d})^{j_d} = J! \cdot \partial_{X^J}$$ 
where $J! = \prod_{\alpha = 1}^d j_{\alpha}!$ in the multi-index notation. 
 
In characteristic zero, since $J! \neq 0$, the above relation implies that 
all the differential 
operators are expressed as (the linear combinations over $R$ of) the 
composites of the differential operators of 
degree $\leq 1$, e.g., $R$-homomorphisms and $\partial_{x_1}, \ldots , 
\partial_{x_d}$. 
 
In positive characteristic $\on{char}(k) = p > 0$, however, $J!$ could 
well be equal to 0 and hence we start seeing the 
differential operators of higher order which cannot be expressed as (the 
linear combinations over $R$ of) the composites of 
differential operators of lower degrees, e.g., 
$$\partial_{x_{\alpha}^{p^1}}, \partial_{x_{\alpha}^{p^2}}, \ldots , 
\partial_{x_{\alpha}^{p^e}}, \ldots 
\text{\ for\ }\alpha = 1, \ldots , d \text{\ and\ }e \in \bZ_{> 0}.$$ 
It is these operators which play a crucial role in 
positive characteristic. 
\item[(2)] 
The following observation comes in handy when we compute the binomial 
coefficients in positive characteristic $\on{char}(k) = p > 0$: 
 
Let $i = \sum_ea_ep^e$ and $j = \sum_eb_ep^e$ be the expressions of the 
integers $i, j \in \bZ_{\geq 0}$ as $p$-adic numbers with 
$0 \leq a_e, b_e < p$.  Then we have 
$$\binom{i}{j} = \prod_e \binom{a_e}{b_e} \mod{p}.$$ 
The identity follows immediately from the observation that, in 
$(\bZ/p\bZ)[x]$, the number $\binom{i}{j}$ is the 
coefficient of 
$x^j =\prod_ex^{b_ep^e}$ in the polynomial 
$(1 + x)^i = \prod_e(1 + x)^{a_ep^e} = 
\prod_e(1 + x^{p^e})^{a_e}$, 
which can be computed as the product of the coefficients 
$\binom{a_e}{b_e}$ 
of $x^{b_ep^e}$ in $(1 +x^{p^e})^{a_e}$. 
\end{rem} 
\end{subsection} 
\begin{subsection}{Logarithmic differential operators.}\label{1.2.2} 
\begin{defn}\label{1.2.2.1} 
Let $E$ be a simple normal crossing divisor on 
$\on{Spec}R$, and $I_E \subset R$ its defining 
ideal. We define the set $\dif_{R,E}^n$ of the 
logarithmic differential opearators of degree 
$\leq n$ on $R$ with respect to $E$ by 
$$\dif^n_{R,E} = \{d \in \dif_R^n\mid 
d(I_E^t) \subset I_E^t \ \forall t \in \bZ_{\geq 0}\}.$$ 
\end{defn} 
\begin{lem}\label{1.2.2.2} 
Suppose we have a regular system of parameters 
$(x_1, \ldots, x_d)$ for $R$ with 
$d = 
\dim W$, and a simple normal crossing divisor $E$ defined by 
$I_E = (\prod_{i = 1}^m x_i)$ for some $1 \leq m \leq d$. 
Then we have the following: 
\item[(1)] 
The $R$-module $\dif^n_{R,E}$ is free of rank $\binom{n+d}{n}$. 
It has a basis 
$\{X^{J_E}\partial_{X^J}\mid |J| 
\leq n\}$ (cf.~Lemma \ref{1.2.1.2} (1)), 
where $J_E = (j_1, \ldots, j_m, 0, \ldots, 0)$ for 
$J = (j_1, \ldots,j_m,$
$ j_{m+1}, \ldots, j_d)$.  Thus we have 
$$\dif^n_{R,E} = \bigoplus_{|J| \leq n}RX^{J_E}\partial_{X^J} \cong 
R^{\binom{n+d}{n}}.$$ 
\item[(2)] 
We have the logarithmic version of the generalized product formula 
$$X^{J_E}\partial_J(fg) = 
\sum_{K+L=J}X^{K_E}\partial_K(f)X^{L_E}\partial_L(g) 
\text{\ for\ }f,g \in R \quad (\text{or\ }\widehat{R}).$$ 
\end{lem} 
\begin{proof} 
This follows immediately from Lemma \ref{1.2.1.2} and 
Definition \ref{1.2.2.1}. 
\end{proof} 
\begin{rem}\label{1.2.2.3} 
We first learned the explicit use of 
the logarithmic differential operators 
in the context of resolution of singularities from 
\cite{MR907903} and \cite{MR1440306}. 
It is worthwhile noting that even when we look at the 
existing algorithms which only use the usual differential 
operators on the surface 
(e.g.~\cite{MR1748620}\cite{MR1949115}\cite{MR2163383}), 
one could implicitly observe the use of logarithmic ones 
in the proof of Giraud's lemma (cf.~\cite{MR0460712}) 
they depend upon. We invite the reader to look at 
\cite{B_BIRS} \cite{B_Harvard} and \cite{MR2078560} for the 
discussions on how the use of the logarithmic differential 
operators, in contrast to the use of the usual ones, affects 
the functorial properties of the algorithm, and even the 
formulation of the problem of resolution of singularities. 
 
The use of the logarithmic differential operators is a ``must'' for our 
algorithm to function, as we will 
see in Parts III and IV, and is recognized as one of the key ingredients of 
the Kawanoue program from the very beginning of its 
conception. 
\end{rem} 
\end{subsection} 
\begin{subsection}{Relation with multiplicity.}\label{1.2.3} 
We end this section 
by pointing out a basic 
relation between the multiplicity (order) and the differential operators in 
the form of a lemma.  It is because 
of this basic relation that the differential operators play a key role in 
constructing an algorithm for 
resolution of singularities, where the order function constitutes a 
fundamental invariant. 
\begin{lem}\label{1.2.3.1} 
Let $I \subset R$ be an ideal.  Let $P \in \on{Spec}R$ 
be a point.  Then 
$$\on{ord}_P(I) \geq n \Longleftrightarrow P \in 
V(\dif_R^{n-1}(I)).$$ 
In particular, the order function 
$\on{ord}_{*}(I)\colon\on{Spec}R \rightarrow \bZ_{\geq 0}$ is 
upper semi-continuous. 
\end{lem} 
\begin{proof} 
First we show the equivalence in the case when $P$ is a 
closed point.  Let $\fm \subset R$ be the maximal 
ideal corresponding to the closed point 
$P$.  Let $\widehat{R}$ be the completion of $R$ with respect to 
$\fm$.  Note that $\on{ord}_P(I) = 
\on{ord}_P(\widehat{I})$, where $\widehat{I} = I\widehat{R}$.  On the 
other hand, since 
$\dif_{\widehat{R}}^{n-1}(\widehat{I}) = 
\dif_R^n(I)\widehat{R}$ by Lemma \ref{1.1.2.1} (7) 
and since $\widehat{R}$ 
is faithfully flat over $R$, we have 
$\dif_{\widehat{R}}^{n-1}(\widehat{I}) \cap R = 
\dif_R^n(I)$.  Thus we have only to show the equivalence 
at the level of completion.  Choose a regular system of parameters 
$(x_1, \ldots , x_d)$ for $R_\fm$. 
Identify $\widehat{R}$ with the power series ring 
$k[[x_1, \ldots , x_d]]$.  By 
definition, $\on{ord}_P(\widehat{I}) \geq n$ if and 
only if, given 
$f = \sum_J c_JX^J \in \widehat{I} \subset 
k[[x_1, \ldots,x_d]]$ with $c_J \in k$, we have 
$c_J = 0$ for any $J$ with $|J| < n$. 
By Lemma \ref{1.2.1.2} (2), the last condition is 
equivalent to saying 
$\partial_{X^K}(f) \subset \widehat\fm$ for any 
$f\in\widehat{I}$ and $K$ with $|K| < n$.  Since 
$\{\partial_{X^K}\mid |K| < n\}$ generates 
$\dif_{\widehat{R}}^{n-1}$ as an 
$\widehat{R}$-module (cf.~Lemma \ref{1.2.1.2} (2)), 
this condition is equivalent 
to $\dif_{\widehat{R}}^{n-1}(\widehat{I}) \subset 
\widehat\fm$, i.e., $P \in 
V(\dif_{\widehat{R}}^{n-1}(\widehat{I}))$.  Therefore, we conclude 
$$\on{ord}_P(\widehat{I}) \geq n \Longleftrightarrow P \in 
V(\dif_{\widehat{R}}^{n-1}(\widehat{I})).$$ 
 
{From} the above argument it follows that the equivalence asserted in the 
lemma holds for a closed point and that the order function 
is upper semi-continuous if we restrict ourselves to the space of the 
maximal ideals $\fm$-$\on{Spec}R$. 
 
It is then straightforward to see that the same equivalence holds for an 
arbitrary point in $\on{Spec}R$ and that the order 
function is upper semi-continuous over $\on{Spec}R$. 
 
This completes the proof of Lemma \ref{1.2.3.1}. 
\end{proof} 
\end{subsection} 
\end{section} 
\begin{section} 
{Ideals generated by the $p^e$-th power elements.}\label{1.3} 
In this section, we denote by $k$ an algebraically closed field of 
$\on{char}(k) = p > 0$. 
 
The purpose of this section is to give a 
characterization of the ideals 
generated by $p^e$-th power elements, 
fixing $e\in \bZ_{\geq 0}$, 
as the ideals invariant under the action of the set 
of differential operators of degree $\leq p^e - 1$. 
 
We denote by $R$ the coordinate ring of an affine open subset 
$\on{Spec}R$ of a variety $W$ smooth over $k$, or its localization 
at a maximal ideal.  We 
denote by $\widehat{R}$ the completion 
of $R$ with respect to a maximal ideal of $R$. 
\begin{subsection} 
{Characterization in terms of the 
differential operators.}\label{1.3.1} 
\begin{defn}\label{1.3.1.1} 
Fix a nonnegative integer 
$e \in \bZ_{\geq 0}$. 
We denote the $e$-th power of the Frobenius map 
by $$F^e\colon R \rightarrow R$$ 
i.e., $F^e(r) = r^{p^e}$ for $r \in R$. 
We use the same symbol $F^e$ for 
the $e$-th power of the Frobenius map of the localization 
$R_S$ or the completion $\widehat{R}$ by abuse of notation if there is no 
chance of confusion. 
\end{defn} 
\begin{prop}\label{1.3.1.2} 
Let $I \subset R$ be an ideal.  Fix a 
nonnegetive integer $e \in \bZ_{\geq 0}$.  Then the 
following conditions are equivalent: 
\begin{enumerate} 
\item 
The ideal $I$ is generated by the $p^e$-th power elements, i.e., 
$I = (I\cap F^e(R))$. 
\item 
The ideal $I$ is invariant under the action of the set of the 
differential operators of degree $\leq p^e - 1$, i.e., $I = 
\dif_R^{p^e-1}(I)$. 
\end{enumerate} 
Moreover, the equivalence of conditions 
$\on{(1)}$ and $\on{(2)}$ also holds over the 
completion $\widehat{R}$. 
\end{prop} 
 
Before beginning the proof of Proposition \ref{1.3.1.2}, 
we remark a couple of facts in the form of a lemma. 
 
\begin{lem}\label{1.3.1.3} 
Let $R'$ denote the localization $R_S$ with respect 
to a multiplicative set $S \subset R$ or the 
completion $\widehat{R}$ with respect to a maximal ideal $\fm \subset 
R$.  Then we have 
\item[(1)] 
$R \otimes_{F^e(R)}F^e(R') = R'$, 
\item[(2)] 
$\{I \cap F^e(R)\}R' = \{IR' \cap F^e(R')\}R'$. 
\end{lem} 
\begin{proof} 
\item[(1)] 
When $R' = R_S$, the assertion is clear since $R 
\otimes_{F^e(R)}F^e(R_S) = RF^e(R_S) = R_S$.  When $R' = 
\widehat{R}$, we see that $R \otimes_{F^e(R)}F^e(\widehat{R})$ and 
$\widehat{R}$ are the completions of $R$ with respect to the 
topologies defined by $\{F^e(\fm^n)R\}_{n \in \bZ_{> 0}}$ and 
$\{\fm^n\}_{n \in \bZ_{> 0}}$ respectively.  It 
is easy to see that these two topologies coincide. 
\item[(2)] 
Since $I \cap F^e(R)\subset IR' \cap F^e(R')$, we have the inclusion 
$\{I \cap F^e(R)\}R' \subset \{IR' \cap F^e(R')\}R'$.  In 
order to see the opposite inclusion, using the fact that $F^e(R')$ 
is flat over $F^e(R)$, we observe 
\begin{eqnarray*} 
\{I \cap F^e(R)\}R' &\supset& \{I \cap F^e(R)\}F^e(R') = \{I \cap F^e(R)\} 
\otimes_{F^e(R)}F^e(R') \\ 
&=& \{I \otimes_{F^e(R)}F^e(R')\} \cap \{F^e(R) \otimes_{F^e(R)}F^e(R')\} \\ 
&=& \{I \otimes_RR \otimes_{F^e(R)}F^e(R')\} \cap F^e(R') \\ 
&=& \{I \otimes_RR'\} \cap F^e(R') = IR' \cap F^e(R'),\\ 
\end{eqnarray*} 
which implies the desired inclusion. 
 
This completes the proof of Lemma \ref{1.3.1.3}. 
\end{proof} 
\begin{proof}[Proof of Proposition \ref{1.3.1.2}] 
\begin{step}{Reduction to the case over the completion $\widehat{R}$.} 
Firstly note that two ideals of $R$ coincide if their localizations or even 
completions coincide at {\it any} maximal ideal 
$\fm$ of $R$.  Thus it suffices to show the conditions 
$$I\widehat{R_\fm} = (I \cap F^e(R))\widehat{R_\fm} 
\quad\text{and}\quad 
I\widehat{R_\fm} =\dif_R^{p^e-1}(I)\widehat{R_\fm}$$ 
are equivalent for any maximal ideal $\fm \subset R$.  Secondly note 
that 
$$\begin{array}{llll} 
(I \cap F^e(R))\widehat{R_\fm} 
&=&\{I\widehat{R_\fm} \cap 
F^e(\widehat{R_\fm})\}\widehat{R_\fm} 
\quad&(\text{by\ Lemma\ \ref{1.3.1.3}\ (2)}), 
\\ 
\dif_R^{p^e-1}(I)\widehat{R_\fm} 
&=&\dif_{\widehat{R_\fm}}^{p^e-1}(I\widehat{R_\fm}) 
\quad&(\text{by\ Lemma\ \ref{1.1.2.1}\ (7)}). 
\end{array}$$ 
Therefore, it suffices to show the equivalence 
of the conditions in the case 
over $\widehat{R}=\widehat{R_\fm}$. 
In the following consideration, we identify 
$\widehat{R}$ with the power series ring 
$k[[x_1, \ldots, x_d]]$ (by choosing a regular 
system of parameters $(x_1, \ldots, x_d)$ for $R_\fm$). 
\end{step} 
\begin{step} 
{Verification of the implication $\on{(i)}\Longrightarrow\on{(ii)}$.} 
We obviously have $\widehat{I} \subset 
\dif^{p^e - 1}_{\widehat{R}}(\widehat{I})$.  Thus we have only to 
show $\widehat{I} \supset 
\dif^{p^e - 1}_{\widehat{R}}(\widehat{I})$ assuming condition (i). 
By Lemma \ref{1.2.1.2} (2), the set 
$\{\partial_{X^J}\mid |J| \leq p^e - 1\}$ generates 
$\dif^{p^e - 1}_{\widehat{R}}$ as an 
$\widehat{R}$-module.  Therefore, it suffices to check $\partial_{X^J}(f) 
\in \widehat{I}$ for any $f \in \widehat{I}$ and 
$\partial_{X^J}$ with $|J| \leq p^e - 1$.  By assuming condition (i), we may 
assume 
$\widehat{I} = 
(\{r_{\lambda}^{p^e}\mid r_{\lambda} 
\in \widehat{R}\}_{\lambda \in \Lambda})$ 
so that we can write 
$f = \sum_{\lambda \in \Lambda} a_{\lambda}r_{\lambda}^{p^e}$ 
with $a_{\lambda} \in \widehat{R}$. 
We compute via the generalized product rule 
\begin{eqnarray*} 
\partial_{X^J}(f) &=& \partial_{X^J}\left(\sum_{\lambda \in \Lambda} 
a_{\lambda}r_{\lambda}^{p^e}\right) = \sum_{\lambda \in \Lambda} 
\partial_{X^J}\left(a_{\lambda}r_{\lambda}^{p^e}\right) \\ 
&=& \sum_{\lambda \in \Lambda} 
\left\{\sum_{K + L = 
J}\partial_{X^K}\left(a_{\lambda}\right)\partial_{X^L}\left(r_{\lambda}^{p^e 
}\right)\right\} = \sum_{\lambda 
\in 
\Lambda}\partial_{X^J}\left(a_{\lambda}\right)r_{\lambda}^{p^e} 
\in I. 
\end{eqnarray*} 
Note that, in order to obtain the last equality, we use the fact that 
$\partial_{X^L}(r_{\lambda}^{p^e}) = 0$ unless $L = 0$. 
In fact, if $r_{\lambda} = \sum_Jc_JX^J \in k[[x_1, \ldots, x_d]]$, 
then, by Lemma \ref{1.2.1.2} (2), we have 
$$\partial_{X^L}(r_{\lambda}^{p^e}) = \partial_{X^L} 
(\sum_Jc_J^{p^e}X^{p^eJ}) = \sum_Jc_J^{p^e}\partial_{X^L}(X^{p^eJ}) 
= \sum_Jc_J^{p^e}\binom{p^eJ}{L}X^{p^eJ - L}.$$ 
Since 
$\displaystyle \binom{p^eJ}{L} = \prod_{i = 1}^d \binom{p^ej_i}{l_i}$, 
and since $\displaystyle \binom{p^ej_i}{l_i} = 0 \bmod{p}$ unless 
$l_i = 0$ because 
$l_i\leq |L| \leq |J| \leq p^e - 1$, we conclude 
$\partial_{X^L}(r_{\lambda}^{p^e}) = 0$ 
unless $L =0$. 
 
This completes the verification of the implication (i) $\Longrightarrow$ 
(ii). 
\end{step} 
\begin{step}{Verification of the implication 
$\on{(ii)}\Longrightarrow\on{(i)}$.} 
We obviously have $\widehat{I} \supset (\widehat{I} \cap F^e(\widehat{R}))$. 
Thus we have only to show $\widehat{I} \subset 
(\widehat{I} \cap F^e(\widehat{R}))$ assuming condition (ii).  First note 
that, setting $\Gamma = \{0, 1, \ldots, p^e-1\}^d$, we can 
express any $f \in \widehat{R} = k[[x_1, \ldots, x_d]]$ in the form 
$$f = \sum_{M \in \Gamma}a_M^{p^e}X^M,$$ 
where the set of coefficients $\{a_M^{p^e}\mid a_M \in 
\widehat{R}\}_{M \in \Gamma}$ is uniquely 
determined. 
 
It suffices to show that, given $f \in \widehat{I}$ 
and its expression as above, we have 
$\{a_M^{p^e}\mid a_M\in\widehat{R}\}_{M\in \Gamma} 
\subset\widehat{I}$, which implies 
$f \in (\widehat{I} \cap F^e(\widehat{R}))$. 
 
We derive a contradiction assuming 
$\{a_M^{p^e}\mid a_M\in \widehat{R}
\}_{M \in \Gamma}\not\subset\widehat{I}$. 
Set 
$$N = \max\left\{M\in\Gamma\mid a_M^{p^e} 
\not\in\widehat{I}\right\},$$ 
where the maximum is taken with respect 
to the lexicographical order on $\Gamma$. 
We compute via the generalized product rule 
\begin{eqnarray*}
\widehat{I} &=& \dif^{p^e - 1}_{\widehat{R}}(\widehat{I}) 
\ni\partial_{X^N}(f - \sum_{M > N}a_M^{p^e}X^M) 
=\sum_{M \leq N}\partial_{X^N}(a_M^{p^e}X^M) 
\\
&=& \sum_{M \leq N}\sum_{K+L=N}\partial_{X^K} 
(a_M^{p^e})\partial_{X^L}(X^M) 
= \sum_{M \leq N}a_M^{P^e}\partial_{X^N}(X^M) = a_N^{p^e}.
\end{eqnarray*} 
Note that, by the same argument as in Step 2 of this proof, 
we see $\partial_{X^K}(a_M^{p^e}) = 0$ unless $K = 0$. 
This is used to obtain the second last equality. 
Note also that 
$
\binom{M}N= 0$ if $M < N$. 
Indeed, if $M<N$, there exists 
$1 \leq i_o \leq d$ such that $m_{i_o} < n_{i_o}$, 
which implies 
$
\binom{m_{i_o}}{n_{i_o}} = 0$ (cf.~\ref{0.9.3}) 
and hence 
$
\binom{M}N=\prod_{i = 1}^d \binom{m_i}{n_i}=0$. 
Thus $\partial_{X^N}(X^M) = 0$ if $M<N$. 
This is used to obtain the last equality. 
Therefore, we have 
$a_N^{p^e} \in \widehat{I}$, contradicting the choice of $N$. 
This completes the verification of the implication 
$\on{(ii)}\Longrightarrow\on{(i)}$. 
 
This completes the proof of Proposition \ref{1.3.1.2}. 
\end{step}\end{proof} 
We end this section by stating a lemma, which is proved 
in the same spirit as the proof of Proposition \ref{1.3.1.2} 
and is of interest on its own. 
\begin{lem}\label{1.3.1.4} 
Let $\widehat{R}$ be the completion of $R$ with 
respect to a maximal ideal $\fm$ of $R$, and $e \in 
\bZ_{\geq 0}$ a nonnegative integer.  Then 
$$R \cap F^e(\widehat{R}) = F^e(R).$$ 
In other words, if $r \in R$ has its $p^e$-th root $f$ within $\widehat{R}$, 
then $f$ actually belongs to $R$. 
\end{lem} 
\begin{proof} 
Since $\widehat{R}$ is faithfully flat over $R$, so is 
$F^e(\widehat{R})$ over $F^e(R)$.  Applying Theorem 7.5 
in \cite{MR879273} to an $F^e(R)$-module $R/F^e(R)$, 
we see that the natural map 
$$R/F^e(R) \rightarrow (R/F^e(R)) \otimes_{F^e(R)} F^e(\widehat{R})$$ 
is injective.  On the other hand, using Lemma \ref{1.3.1.3}, we analyze 
the target of the above map to be 
$$(R/F^e(R)) \otimes_{F^e(R)} F^e(\widehat{R}) = \{R \otimes_{F^e(R)} 
F^e(\widehat{R})\}/\{F^e(R) \otimes_{F^e(R)} 
F^e(\widehat{R})\} = \widehat{R}/F^e(\widehat{R}).$$ 
That is to say, we conclude that the map $R/F^e(R) \rightarrow 
\widehat{R}/F^e(\widehat{R})$ is 
injective, and hence that $R \cap F^e(\widehat{R}) = F^e(R)$. 
 
This completes the proof of Lemma \ref{1.3.1.4}. 
\end{proof} 
\end{subsection} 
\end{section} 
\end{chapter} 
\begin{chapter}{Idealistic Filtration} 
The purpose of this chapter is to introduce the notion of an idealistic 
filtration, which is the main language 
to describe our program toward constructing an algorithm for resolution of 
singularities, and establish its 
fundamental properties. 
 
We develop our argument over a ring $R$, which is assumed to be the 
coordinate ring of an affine open 
subset of a nonsingular variety $W$ over $k$, or its localization, or its 
completion with respect to a maximal ideal. 
That is to say, more geometrically speaking, we carry out 
our analysis over an affine open 
subset of a nonsingular variety $W$, or over a stalk, or over the analytic 
structure 
at a closed point.  Since the main operations on an idealistic filtration, 
such as the operations of taking 
the $\fD$-saturation and $\fR$-saturation, are compatible with 
localization and 
completion (for an idealistic filtration of r.f.g.~type), it is immediate to 
extend the 
(analytically) local analysis of this chapter to the global argument, which 
we will develop in the 
subsequent papers. 
\begin{section}{Idealistic filtration over a ring.} 
\label{2.1} 
Let $R$ be the coordinate ring of an affine open 
subset of a nonsingular variety $W$ over $k$, or its localization, or its 
completion with respect to a maximal ideal. 
\begin{subsection}{Definitions.}\label{2.1.1} 
\begin{defn}\label{2.1.1.1} 
\item[(1)] 
Let $T \subset R \times \bR$ be a 
subset.  For $a \in \bR$, we set 
$T_a = \{f \in R\mid (f,a)\in T\}$. 
\item[(2)] 
We call a subset $\bI \subset R \times\bR$ 
an idealistic filtration if it satisfies the following conditions: 
$$\left\{\begin{array}{cl} 
(\on{o})& 
\bI_0 = R, 
\\ 
(\on{i})& 
\bI_a \text{\ is\ an\ ideal\ of\ }R 
\text{\ for\ any\ }a \in \bR 
\\ 
&\quad 
(\bI_a \text{\ is\ called\ the\ ideal\ of\ }\bI \text{\ at\ level\ }a) 
\\ 
(\on{ii})& 
\bI_a\bI_b \subset \bI_{a+b} 
\text{\ for\ any\ }a, b \in \bR, 
\\ 
(\on{iii})& 
\bI_b \supset \bI_a \text{\ if\ }b \leq a. 
\end{array}\right. 
$$ 
\item[(3)] 
Let $T \subset R \times \bR$ be a subset.  We call 
the minimal idealistic filtration 
containing 
$T$ the idealistic filtration generated by $T$ and denote it by $G(T)$.  If 
$\bI = G(T)$, we call $T$ a set of 
generators for $\bI$ (cf.~Lemma \ref{2.2.1.1} (2)). 
 
When we want to emphasize the base ring $R$ over which $T$ generates 
the idealistic filtration, 
we write $G_R(T)$ inserting $R$ as a subscript. 
\item[(4)] 
We say an idealistic filtration $\bI$ is of r.f.g.~type (short for 
rationally and finitely 
generated) if there exists a finite set 
$T \subset R \times \bQ \subset R \times \bR$ 
such that $\bI = G(T)$. 
\item[(5)] 
Let $T \subset R \times \bR_{\geq 0}$ be a subset.  Let $P 
\in 
\on{Spec}R$ be a point.  We define the multiplicity 
$\mu_P(T)$ of $T$ at $P$ to be 
$$\mu_P(T) := \inf\left\{\mu_P(f,a) := 
\frac{\on{ord}_P(f)}{a}\mid (f,a) \in T, a > 
0\right\}.$$ 
Note that we set $\mu_P(f,0) = \infty$ 
for any $f\in R$ by definition, 
while $\on{ord}_P(0) = \infty$. 
\item[(6)] 
Let $T \subset R \times \bR$ be a subset. 
We define the support 
$\on{Supp}(T)$ of $T$ to be 
$$\on{Supp}(T) = \{P \in \on{Spec}R\mid 
\mu_P(T)\geq 1\}.$$ 
\end{defn} 
\begin{rem}\label{2.1.1.2} 
\item[(1)] 
It is straightforward to see that a subset $\bI \subset R \times\bR$ 
is an idealistic filtration 
\linebreak 
if and only if it satisfies the following conditions: 
$$\left\{\begin{array}{cl} 
(\on{o})& (f,0) \in \bI \quad \forall f \in R, 
\quad (0,a) \in \bI \quad \forall a \in \bR,\\ 
(\on{i})&(f,a), (g,a) \in \bI 
\Longrightarrow (f + g,a)\in \bI,\\ 
&r \in R, (f,a) \in \bI \Longrightarrow (rf,a) \in \bI,\\ 
(\on{ii})&(f,a), (h,b) \in \bI \Longrightarrow (fh,a + b) 
\in \bI,\\ 
(\on{iii})&(f,a) \in \bI, b \leq a \Longrightarrow (f,b) 
\in \bI. 
\end{array}\right.$$ 
We invite the reader to look at \ref{0.2.3.1} in 
Chapter 0 for the motivation 
behind introducing the notion of an idealistic filtration. 
\item[(2)] 
When $T = \bI \subset R \times \bR$ is an idealistic 
filtration, we define its multiplicity $\mu_P(\bI)$ 
at a point $P \in\on{Spec}R$, 
and its support $\on{Supp}(\bI)$ acccording to 
Definition \ref{2.1.1.1} (5) and (6). 
\end{rem} 
\end{subsection} 
\begin{subsection}{$\fD$-saturation.}\label{2.1.2} 
We define the notion of the differential saturation 
\linebreak 
(which we call the $\fD$-saturation for short) 
of an idealistic filtration. 
Budding of an idea leading to the notion of 
$\fD$-saturation can be observed in the work of 
Giraud and Villamayor, where they discuss the enlargement, 
called the extension, of an ideal obtained by adding the partial 
derivatives of the elements in the ideal. 
\begin{defn}\label{2.1.2.1} 
Let $\bI \subset R \times \bR$ be an 
idealistic filtration.  We say $\bI$ 
is $\fD$-saturated if it satisfies the following condition 
(differential): 
$$(\on{differential}) \quad (f,a) \in \bI, d \in 
\dif^t_R \Longrightarrow (d(f), a-t) \in \bI.$$ 
(We refer the reader to Chapter 1 for the meaning of the notation 
$\dif^t_R$.) 
 
Let $\bI$ be an idealistic filtration.  We call the minimal 
$\fD$-saturated idealistic filtration containing 
$\bI$ the differential saturation (or $\fD$-saturation 
for short) of $\bI$, and denote it by $\fD(\bI)$ 
(cf.~Lemma \ref{2.2.1.1}). 
 
Let $E$ be a simple normal crossing divisor on $W$.  Then using the 
logarithmic differential operators with respect to $E$ instead of the 
usual differential operators 
(cf.~Definition \ref{1.2.2.1}),  we consider the following 
condition $(\on{differential})_E$: 
$$(\on{differential})_E \quad (f,a) \in \bI, d \in 
\dif^t_{R,E} \Longrightarrow (d(f), a - t) \in \bI.$$ 
Replacing condition (differential) with condition 
$(\on{differential})_E$, 
we obtain the notion of an idealistic filtration being 
$\fD_E$-saturated and that of the 
$\fD_E$-saturation. 
\end{defn} 
\end{subsection} 
\begin{subsection}{$\fR$-saturation.} 
\label{2.1.3} 
We define the notion 
of the radical saturation (which we call the 
$\fR$-saturation for short) of an idealistic filtration. 
Note that, for an 
$\fR$-saturated idealistic filtration, we not only require 
that we can take the $n$-th root (radical) of an element within the 
idealistic filtration (if it exists within $R \times \bR$) for any 
$n\in \bZ_{> 0}$, but also require the continuity by definition. 
\begin{defn} 
\label{2.1.3.1} 
Let $\bI \subset R \times \bR$ 
be an idealistic filtration. 
We say $\bI$ is $\fR$-saturated if 
it satisfies the following 
conditions (radical) and (continuity): 
$$ 
\begin{array}{lll} 
(\on{radical}) 
&(f^n,na) \in \bI, f \in R, n \in \bZ_{> 0} 
&\Longrightarrow (f,a) \in \bI\\ 
(\on{continuity}) 
&\{(f,a_l)\} \subset \bI \text{\ with\ } 
\lim_{l \rightarrow \infty}a_l =a 
&\Longrightarrow (f,a) \in \bI. 
\end{array}$$ 
 
Let $\bI$ be an idealistic filtration.  We call the minimal 
$\fR$-saturated idealistic filtration containing 
$\bI$ the radical saturation (or $\fR$-saturation for short) of 
$\bI$, and denote 
it by 
$\fR(\bI)$ (cf.~Lemma \ref{2.2.1.1}). 
\end{defn} 
\begin{rem}\label{2.1.3.2} 
\item[(1)] 
We remark that, in positive 
characteristic $p = \on{char}(k)> 0$, if an idealistic filtration 
$\bI$ satisfies the following condition (Frobenius), which is a priori 
slightly weaker than condition (radical), and 
condition (continuity), then it actually satisfies conditions (radical) and 
(continuity).  Therefore, instead of checking 
conditions (radical) and (continuity) in order to show that a given 
idealistic filtration is 
$\fR$-saturated in positive characteristic, we could check conditions 
(Frobenius) and (continuity): 
$$(\on{Frobenius}) \quad (f^p,pa) \in \bI, f \in R 
\Longrightarrow (f,a) \in \bI.$$ 
 
In fact, suppose we have $(f^n,na) \in \bI$, $f\in R$ 
and $n \in \bZ_{>0}$.  Take $e \in \bZ_{> 0}$ 
so that $p^e > n$, and take $r\in \bZ_{\geq 0}$ 
with $0 \leq r < n$ so that $r\equiv p^e\mod{n}$.  Then 
$$ 
\begin{array}{ll} 
(f^n,na)\in\bI \Longrightarrow (f^{p^e-r},a \cdot (p^e-r))\in \bI 
& 
\text{by\ condition\ (ii)\ in\ Remark\ \ref{2.1.1.2} (1)} 
\\ 
\phantom{(f^n,na) \in \bI}\Longrightarrow (f^{p^e},a \cdot (p^e - r)) 
\in \bI 
&\text{by\ condition\ (i)\ in\ Remark\ \ref{2.1.1.2} (1)} 
\\ 
\phantom{(f^n,na) \in \bI} 
\Longrightarrow (f,a \cdot \left(1 - p^{-e}r\right)) \in \bI 
&\text{by\ condition\ (Frobenius)} 
\\ 
\phantom{(f^n,na) \in \bI} 
\Longrightarrow (f,a) \in \bI 
&\text{by\ condition\ (continuity)\ 
with\ }e \rightarrow \infty. 
\end{array} 
$$ 
\item[(2)] 
In view of condition (iii) in Remark \ref{2.1.1.2}(1), 
requiring condition (continuity) is equivalent to requiring 
the following (left continuity): 
$$(\text{left continuity}) \quad \{(f,a_l)\} \subset \bI 
\text{\ with\ }\{a_l\} \text{\ increasing\ and\ } 
\lim_{l\rightarrow\infty}a_l = a \Longrightarrow 
(f,a) \in \bI.$$ 
In terms of the ideals of an idealistic filtration 
associated to the levels, condition (left continuity) 
translates into the condition 
$$\bI_a = \bigcap_{b <a}\bI_b \quad 
\forall a \in \bR.$$ 
When an idealistic filtration is of r.f.g.~type, 
this condition can be checked rather easily. 
Therefore, we see that condition 
(continuity) is always satisfied for an idealistic 
filtration of r.f.g.~type. 
See Corollary \ref{2.3.2.2} for detail. 
\end{rem} 
\end{subsection} 
\begin{subsection}{Integral closure.}\label{2.1.4} 
We define the notion of the 
integral closure of an idealistic filtration, which is closely 
related to the notion of the $\fR$-saturation.  In general, if an 
idealistic filtration is $\fR$-saturated, then it is integrally 
closed.  In particular, for an idealistic filtration of r.f.g.~type, where 
condition (continuity) is automatic, it is $\fR$-saturated 
if and only if it is integrally closed. 
 
We also conclude in Corollary \ref{2.3.2.6}, through 
the argument showing the stability of 
r.f.g.~type 
under $\fR$-saturation, that, for an idealistic 
filtration of r.f.g.~type, the $\fR$-saturation and 
the integral closure coincide. 
\begin{defn}\label{2.1.4.1} 
Let $\bI \subset R \times \bR$ be an 
idealistic filtration. 
\item[(1)] 
We say an element $(f,a) \in R \times \bR$ is integral over 
$\bI$ if $f$ satisfies a monic equation of the form 
$$f^n + c_1f^{n-1} +\cdots+ c_n = 0 
\text{\ with\ }(c_i,ia) \in\bI 
\text{\ for\ }i = 1, \dotsc , n.$$ 
\item[(2)] 
We say $\bI$ is integrally closed 
if it satisfies the following condition (ic): 
$$(\on{ic}) \quad (f,a) \in R \times \bR \text{\ is\ integral\ 
over\ } \bI \Longrightarrow (f,a) \in \bI.$$ 
Let $\bI$ be an idealistic filtration.  We call the minimal integrally 
closed idealistic filtration containing 
$\bI$ the integral closure of $\bI$, and denote 
it by $\on{IC}(\bI)$ (cf.~Lemma \ref{2.2.1.1}). 
\end{defn} 
\begin{rem}\label{2.1.4.2} 
The notion of the integral closure is important in our program. 
However, since the $\fR$-saturation and the integral closure 
coincide for an idealistic filtration of r.f.g.~type, and 
since almost all the idealistic filtrations we consider are of 
r.f.g.~type, we seldom use the symbol $\on{IC}(\bI)$ or the 
notion of the integral closure explicitly, and almost always 
use the notion of the $\fR$-saturation, 
which is denoted by $\fR(\bI)$. 
\end{rem} 
\end{subsection} 
\begin{subsection}{$\fB$-saturation.}\label{2.1.5} 
We define the notion of the bi-saturation (which we call the 
$\fB$-saturation).  For the purpose of extracting 
the intrinsic information toward a solution of 
the problem of resolution of singularities, 
we take various saturations of a given idealistic filtration 
(cf.~\ref{0.2.3.2.3}).  It would be best if we could take an 
``optimal'' 
one among such.  In our algorithm, the $\fB$-saturation (or 
$\fB_E$-saturation) plays the role of the optimal saturation. 
\begin{defn}\label{2.1.5.1} 
Let $\bI \subset R \times \bR$ be an 
idealistic filtration.  We say $\bI$ is $\fB$-saturated 
(resp. $\fB_E$-saturated) if it is both $\fD$-saturated 
(resp. $\fD_E$-saturated) and $\fR$-saturated. 
Given an idealistic filtration 
$\bI$, we call the minimal 
$\fB$-saturated ($\fB_E$-saturated) idealistic filtration 
containing $\bI$ the $\fB$-saturation (resp. 
$\fB_E$-saturation) of $\bI$, and denote it by $\fB(\bI)$ 
(resp. $\fB_E(\bI)$) (cf.~Lemma \ref{2.2.1.1}). 
\end{defn} 
\begin{rem}\label{2.1.5.2} 
While the existence of the $\fB$-saturation is as straightforward 
as the existence of the other saturations and integral closure, 
its explicit construction is quite remarkable, 
which we will see in Corollary \ref{2.4.2.3}. 
We describe the 
explicit construction of the other saturations and 
integral closure in Lemma \ref{Construction}. 
\end{rem} 
\end{subsection} 
\end{section} 
\begin{section}{Basic properties of an idealistic filtration.} 
\label{2.2} 
In this section, we discuss some basic properties of an idealistic 
filtration over a 
ring.  We use the same notation as in \ref{2.1}. 
\begin{subsection}{On generation, $\fD$-saturation, 
$\fR$-saturation, integral closure, and $\fB$-saturation.} 
\label{2.2.1} 
The next two lemmas 
discuss the existence and explicit construction of the idealistic 
filtration generated by a subset $T \subset R \times \bR$, the 
$\fD$-saturation, $\fR$-saturation, integral closure, and 
$\fB$-saturation. 
\begin{lem}\label{2.2.1.1} 
\item[(1)] 
The intersection $\bigcap_{\lambda \in 
\Lambda}\bI_{\lambda} \subset R \times \bR$ of a non-empty 
collection $\{\bI_{\lambda}\}_{\lambda \in \Lambda}$ of idealistic 
filtrations is again an idealistic filtration.  Moreover, if each 
$\bI_{\lambda}$ is 
$\fD$-saturated (resp. $\fD_E$-saturated, $\fR$-saturated, 
integrally closed, $\fB$-saturated, 
$\fB_E$-saturated), then so is the intersection 
$\bigcap_{\lambda 
\in \Lambda}\bI_{\lambda}$. 
\item[(2)] 
Let $T \subset R \times \bR$ be a subset.  Then $G(T)$ exists 
(cf.~Definition \ref{2.1.1.1} (3)). 
\item[(3)] 
Let $\bI$ be an idealistic filtration.  Then $\fD(\bI)$ 
(resp. $\fD_E(\bI)$, $\fR(\bI)$, 
$\on{IC}(\bI)$, $\fB(\bI)$, $\fB_E(\bI)$) 
exists (cf.~\ref{2.1.2}, \ref{2.1.3}, \ref{2.1.4}, \ref{2.1.5}). 
\end{lem} 
\begin{proof} 
\item[(1)] 
It 
is clear from the definitions. 
\item[(2)] Let 
${\mathcal S} = \{\bI_{\lambda}\mid 
\bI_{\lambda} \supset T\}$ be the collection of all the 
idealistic filtrations containing $T$.  Note that 
${\mathcal S}$ is non-empty, since 
$R \times \bR \in {\mathcal S}$. Now it is clear 
that the intersection $\bigcap_{\bI_{\lambda} \in {\mathcal S}} 
\bI_{\lambda}$ is the minimal idealistic filtration containing $T$. 
\item[(3)] Let 
${\mathcal S} = \{\bI_{\lambda}\mid 
{\bI}_{\lambda} \supset \bI\}$ be the collection of all the 
$\fD$-saturated (resp. $\fD_E$-saturated, $\fR$-saturated, 
integrally closed, $\fB$-saturated, 
$\fB_E$-saturated) idealistic filtrations containing 
$\bI$.  Note that 
${\mathcal S}$ is non-empty, since $R \times \bR \in {\mathcal S}$. 
Now it is 
clear that the intersection $\bigcap_{\bI_{\lambda} \in 
{\mathcal S}} \bI_{\lambda}$ is the minimal $\fD$-saturated 
(resp. $\fD_E$-saturated, $\fR$-saturated, integrally closed, 
$\fB$-saturated, $\fB_E$-saturated) idealistic filtration 
containing $\bI$. 
 
This completes the proof of Lemma \ref{2.2.1.1}. 
\end{proof} 
\begin{lem}\label{Construction} 
Let $\bI$ be an idealistic filtration 
generated by $T = \{(f_{\lambda},a_{\lambda})\} 
\subset R \times \bR$, i.e., $\bI = G(T)$. 
\item[(1)] 
Define a subset $\bI' \subset R \times \bR$ by setting 
$$\bI'_a = 
\left(\prod f_{\lambda}^{n_{\lambda}}\mid n_{\lambda} 
\in \bZ_{\geq 0}, \sum n_{\lambda}a_{\lambda} \geq a\right) 
\quad a\in\bR. 
$$ 
Then $\bI'$ is an idealistic filtration, and $\bI' = \bI$. 
Note that, when $T = \emptyset$, we use the convention 
$$G(\emptyset) = (\{0\} \times \bR) 
\cup (R 
\times \bR_{\leq 0}).$$ 
\item[(2)] 
Let $(x_1, \dotsc, x_d)$ be a regular system of parameters for $R$.  Set 
$$T' = \left\{(\partial_{X^J}f_{\lambda},a_{\lambda} - 
|J|)\mid J \in \bZ_{\geq 0}^d, 
(f_{\lambda},a_{\lambda}) \in 
T\right\}.$$ 
Then we have $\fD(\bI) = G(T')$. 
 
Let $E$ be a simple normal crossing divisor, and say, $\{x_1 \cdots x_m = 
0\}$ defines $E$ for some $1 \leq m \leq d$. 
Set 
$$T'_E = \left\{(X^{J_E}\partial_{X^J}f_{\lambda},a_{\lambda} - 
|J|)\mid J \in \bZ_{\geq 
0}^d, (f_{\lambda},a_{\lambda}) 
\in T\right\}.$$ 
Then we have $\fD_E(\bI) = G(T'_E)$. 
(We refer the reader to \ref{1.2.2} for the notation.) 
\item[(3)] 
Define subsets $\bK, \overline\bK \subset R \times \bR$ 
by 
$$\bK_a = \{f \in R\mid f^n \in \bI_{na} 
\text{\ for\ some\ }n \in \bZ_{> 0}\}, \quad 
\overline\bK_a = 
\bigcap_{b < a}\bK_b \quad (a \in \bR).$$ 
Then $\overline\bK$ is an idealistic filtration, and 
$\fR(\bI) = \overline\bK$. 
\item[(4)] 
Let $\bJ \subset R \times \bR$ be the subset consisting of all 
the elements integral over $\bI$.  Then $\bJ$ is an 
idealistic filtration, and $\on{IC}(\bI) = \bJ$. 
\end{lem} 
\begin{proof} 
\item[(1)] 
It is straightforward to see that $\bI'$ is an 
idealistic filtration, and that any idealistic filtration 
containing 
$T$ necessarily contains $\bI'$.  Therefore, $\bI = \bI$ by 
the definition of $\bI = G(T)$. 
\item[(2)] 
Let $\bI'$ be a $\fD$-saturated idealistic filtration 
containing $T$, or equivalently containing $\bI$.  Then it 
is clear that $\bI' \supset G(T')$.  Therefore, in order to see 
$\fD(\bI) = G(T')$, we have only to show $G(T')$ is 
$\fD$-saturated, which follows from the fact that $\dif_R^t$ is 
generated by $\left\{\partial_{X^J}\mid |J| \leq 
t\right\}$ as an 
$R$-module, and the generalized product rule 
(cf.~Lemma \ref{1.2.1.2}). 
The proof for 
the case of $\fD_E$-saturation is identical to the case 
of $\fD$-saturation, replacing the usual differentials with the 
logarithmic ones. 
\item[(3)] 
Let $\bI'$ be an $\fR$-saturated idealistic filtration 
containing $\bI$.  As $\bI'$ satisfies condition 
(radical), we have $\bK \subset \bI'$.  Therefore, we conclude 
$$\overline\bK_a = \bigcap_{b < a}\bK_a 
\subset \bigcap_{b <a}\bI_b' = \bI_a',$$ 
where the last equality follows since $\bI'$ satisfies condition 
(continuity) (cf.~Remark \ref{2.1.3.2} (2)).  That is to say, 
$\overline\bK \subset \bI'$. 
 
Thus, in order to see $\fR(\bI) = \overline\bK$, we have 
only to show that $\overline\bK$ itself is an idealistic 
filtration containing 
$\bI$, satisfying conditions (radical) and (continuity). 
 
First we show that $\overline\bK$ is 
an idealistic filtration.  We have only to check that 
$\overline\bK_a \quad (a > 0)$ is 
closed under addition (cf.~Definition \ref{2.1.1.1} (1) 
and Remark \ref{2.1.1.2} (1)), 
while the other conditions follow easily.  Take 
$f, g \in \overline\bK_a$.  Then for any $b < a$, 
there exists $n \in\bZ_{> 0}$ such that $f^n, g^n \in \bI_{nb}$. 
Then for any $k \in \bZ_{> 0}$, we have 
$$(f + g)^k = \sum_{i = 0}^k\binom{k}{i}f^ig^{k-i} \in 
\bI_{nb}^{\lfloor \frac{i}{n}\rfloor + \lfloor \frac{k-i}{n}\rfloor} 
\subset \bI_{b(k - 2n)}$$ 
since $\lfloor \frac{i}{n}\rfloor + \lfloor \frac{k-i}{n}\rfloor \geq 
\frac{k}{n} - 2$.  Therefore, we have $f + g \in \bK_{b(1 - 
2nk^{-1})}$.  Since $b < a$ and $k > 0$ are arbitrary (while $n$ depends 
only on $b$), we conclude $f + g \in \bK_c$ for any $c < a$. 
Therefore, we have $f + g \in \bigcap_{c < a}\bK_c 
= \overline\bK_a$. 
 
Secondly we check condition (continuity) for 
$\overline\bK$.  In fact, we have 
$$\overline\bK_a = \bigcap_{b < a} 
\bK_b = \bigcap_{b <a}\bigcap_{c < b}\bK_c 
= \bigcap_{b < a}\overline\bK_b.$$ 
Therefore, $\overline\bK$ satisfies condition 
(continuity) (cf.~Remark \ref{2.1.3.2} (2)). 
 
Finally we check condition (radical) for $\overline\bK$. 
Suppose $f^n\in \overline\bK_a$.  Fix $b < a$.  Then 
$f^n \in\bK_b$ by definition of $\overline\bK$, 
and there exists $m \in\bZ_{> 0}$ such that 
$(f^n)^m \in \bI_{mb}$ by definition of 
$\bK$.  Therefore, we have 
$f \in \bK_{n^{-1}b}$.  Since $b < a$ 
is arbitrary, we have 
$f \in \bigcap_{b < a}\bK_{n^{-1}b} = 
\overline\bK_{n^{-1}a}$.  Therefore, 
$\overline\bK$ satisfies 
condition (radical). 
\item[(4)] 
It is clear that, if $\bI'$ is an idealistic filtration containing 
$\bI$ and satisfying (ic), then $\bJ \subset 
\bI'$.  Thus, in order to see $\on{IC}(\bI) = \bJ$, we 
have only to show that $\bJ$ itself is an idealistic filtration 
containing 
$\bI$, satisfying condition (ic). 
 
It is clear that $\bJ$ 
contains $\bI$.  Consider the graded subring 
${\mathfrak G}{\mathfrak r}(\bI) := \bigoplus_{a \in \bR}\bI_aX^a 
\subset \bigoplus_{a \in \bR}RX^a$, where 
$X$ is a variable transcendental over $R$, and where the structure of the 
graded $R$-algebra is given through multiplication rule $X^aX^b = 
X^{a+b}$. 
 
Observe (cf.~\cite{MR879273}) that 
$$(f,a) \text{\ is\ integral\ over\ }\bI \Longleftrightarrow 
{\mathfrak G}{\mathfrak r}(\bI)[fX^a] \text{\ is\ a\ finite\ } 
{\mathfrak G}{\mathfrak r}(\bI)\text{-}\text{module}.$$ 
Observe also that 
$${\mathfrak G}{\mathfrak r}(G(\bI,(f,a))) 
= {\mathfrak G}{\mathfrak r}(\bI)[fX^a].$$ 
Now from these observations it follows easily that 
$\bJ$ is an idealistic filtration, and that 
$\bJ$ satisfies condition (ic). 
 
This completes the proof for Lemma \ref{Construction}. 
\end{proof} 
\end{subsection} 
\begin{subsection}{$\fR$-saturated implies integrally closed.} 
 \label{2.2.2} 
\begin{prop}\label{2.2.2.1} 
Let $\bI \subset R \times \bR$ be 
an idealistic filtration.  If 
$\bI$ is $\fR$-saturated, then $\bI$ is integrally closed. 
\end{prop} 
\begin{proof} 
Let $\bI$ be an $\fR$-saturated idealistic filtration. 
Suppose $(f,a) \in R \times \bR$ is integral 
over 
$\bI$, i.e., $f$ satisfies a monic equation of the form 
$$(\star) \quad f^n + c_1f^{n-1} +\cdots+ c_n = 0 \text{\ 
with\ }(c_i,ia) 
\in \bI \text{\ for\ } i = 1, \dotsc , n.$$ 
We want to show $(f,a) \in \bI$. 
 
If $a \leq 0$, then obviously $(f,a) \in \bI$ (cf.~conditions (o), 
(iii) in Remark \ref{2.1.1.2} (1)). 
Thus, we may further assume $a >0$.  Let 
$$\beta_l = 1 - \left(\frac{n - 1}{n}\right)^l \quad (l \in 
\bZ_{\geq 0}).$$ 
We show by induction that 
$$(\heartsuit)_l \quad (f,a\beta_l) \in \bI.$$ 
$(\heartsuit)_0$ is clear.  Suppose we have shown $(\heartsuit)_l$. 
Using the monic equation $(\star)$, we have 
$$f^n = - (c_1f^{n-1} +\cdots+ c_n)$$ 
with 
$$(c_if^{n-i},a\{i + (n - 
i)\beta_l\}) \in \bI \quad (1 \leq i \leq n).$$ 
Since $0 \leq \beta_l < 1$, we have 
$$\min_i\{i + (n - i)\beta_l\} = 1 + (n - 1)\beta_l = n\beta_{l+1}.$$ 
Therefore, we conclude 
$$(f^n,an\beta_{l+1}) \in \bI.$$ 
Since $\bI$ is $\fR$-saturated, it follows from condition 
(radical) 
$$(\heartsuit)_{l+1} \quad (f,a\beta_{l+1}) \in \bI.$$ 
Thus $(\heartsuit)_l$ is valid for all $l \in \bZ_{\geq 0}$. 
 
Note that $\lim_{l \rightarrow \infty}a\beta_l = a$.  Therefore, by 
condition (continuity) satisfied by $\bI$, we conclude 
$$(f,a) \in \bI.$$ 
Therefore, $\bI$ is integrally closed. 
 
This completes the proof of Proposition \ref{2.2.2.1}. 
\end{proof} 
\end{subsection} 
\begin{subsection}{Analysis of interaction between $\fD$-saturation 
and $\fR$-saturation.}\label{2.2.3} 
So far, we have studied the $\fD$-saturation and $\fR$-saturation 
separately.  In this subsection, we analyze the interaction of the 
operations of taking $\fD$-saturation and $\fR$-saturation. 
Under the assumption that $R$ has a regular system of parameters, 
our result is stated in the following proposition, which leads to the 
explicit construction of the $\fB$-saturation.  Furthermore, the 
assumption is later removed for an idealistic filtration 
of r.f.g.~type (cf.~Corollary \ref{2.4.2.3}). 
\begin{prop}\label{2.2.3.1} 
Let $\bI$ be an idealistic filtration 
over $R$ which has a regular 
system of parameters $(x_1, \dotsc, x_d)$.  Then 
$\fD\fR(\bI)\subset \fR\fD(\bI)$. 
 
If $E$ is a simple normal crossing divisor defined by 
$\{x_1 \cdots x_m =0\}$ for some $1 \leq m \leq d$, then 
$\fD_E\fR(\bI)\subset \fR\fD_E(\bI)$. 
\end{prop} 
\begin{proof} 
We present a proof of the latter assertion in the 
logarithmic case, as the former is a special case of the latter 
($E=\emptyset$). 
\begin{step}{%
Reduction of the assertion to the statement $(\clubsuit)$.} 
By replacing $\bI$ with $\fD_E(\bI)$ and via the obvious 
inclusion $\fD_E\fR(\bI) \subset \fD_E\fR(\fD_E(\bI))$, 
we see that it suffices to prove the inclusion 
$$\fD_E\fR(\bI)\subset 
\fR(\bI),$$ 
assuming $\bI$ is $\fD_E$-saturated.  In order to show the first 
inclusion above, by Lemma \ref{Construction} (2), we have only to 
show 
$$ 
\bigcup_J\{(D_J(f),a-|J|)\mid 
(f,a)\in\fR(\bI)\}\subset 
\fR(\bI), 
$$ 
where $D_J = X^{J_E}\partial_{X^J}$.  Now since $D_J=D_{j_1{\mathbf 
e}_1}\cdots D_{j_d{\mathbf e}_d}$, this second inclusion then follows if we 
can show, for $1 \leq i \leq d$, the inclusion below 
$$\bigcup_{j\geq 0}\{(D_{j{\mathbf e}_i}(f) 
,a-j)\mid (f,a)\in\fR(\bI)\} 
\subset\fR(\bI).$$ 
Let $\bK, \overline\bK \subset R \times \bR$ be as in 
Lemma \ref{Construction} (3). 
We claim that we may even replace the range 
$\fR(\bI)$ of $(f,a)$ in the left hand side of 
the third inclusion with $\bK$.  That is to say, 
we claim it suffices to show 
$$ 
(\diamondsuit)\quad\bigcup_{j\geq 0}\{(D_{j{\mathbf e}_i}(f) 
,a-j)\mid (f,a)\in\bK\}\subset\fR(\bI).$$ 
In fact, $(\diamondsuit)$ implies 
$$ 
\{(D_{j{\mathbf e}_i}(f),a-j)\mid (f,a)\in\overline\bK\} 
\subset 
\overline{\{(D_{j{\mathbf e}_i}(f),a-j)\mid (f,a)\in\bK\}} 
\subset\overline{\fR(\bI)},$$ 
where, given a subset $T \subset R \times \bR$, the subset 
$\overline{T}$ is defined by $\overline{T}_a = \bigcap_{b < a}T_b$. 
Since 
$\overline\bK=\fR(\bI)=\overline{\fR(\bI)}$, this 
inclusion then would imply the third one. 
 
Finally, we reduce $(\diamondsuit)$ to the following 
general statement: 
$$(\clubsuit) \quad D_{j{\mathbf e}_i}(f) 
\in \fR\fD_E(G\{(f^n,na)\})_{a-j}\quad 
(f\in R, a\in\bR, 1\leq i\leq d, n>0, j\geq0). 
$$ 
Indeed, given $f\in\bK_a$, there exists $n\in\bZ_{>0}$ 
such that $f^n\in\bI_{na}$. 
Thus, $(\clubsuit)$ implies 
$D_{j{\mathbf e}_i}(f) 
\in \fR\fD_E(G\{(f^n,na)\})_{a-j} 
\subset\fR\fD_E(\bI)_{a-j} 
=\fR(\bI)_{a-j}$ 
since $\bI$ is assumed to be $\fD_E$-saturated. 
Thus $(\clubsuit)$ implies $(\diamondsuit)$. 
 
Therefore, we conclude that the assertion of the lemma is reduced to the 
statement $(\clubsuit)$. 
\end{step} 
\begin{step}{Setup for the inductional proof of $(\clubsuit)$.} 
We fix $1 \leq i \leq d$, and omit $i$ from the notation in the following 
argument.  For example, we denote $D_{j{\mathbf e}_i}$ by $D_j$.  We 
also denote $\fR\fD_E(G(\{(f^n,na)\}))$ by $\bJ$ 
to ease the notation. 
 
Set $c=\lfloor a\rfloor$.  We prove the 
statement $(\clubsuit)$ by induction on $c$.  We may assume 
$0\leq j\leq c$, since otherwise we have 
${a-j}<0$ and $\bJ_{a-j}=R$, 
in which case $(\clubsuit)$ clearly holds. 
\begin{case}{$c = 0$.} 
In this case, $j$ must be $0$, and we 
obviously have 
$$(D_j(f),a-j) = (f,a) \in 
\fR(G(\{(f^n,na)\})) 
\subset \bJ.$$ 
Thus $(\clubsuit)$ holds. 
\end{case} 
\begin{case}{$c \geq 1$.} 
In this case, we show $(\clubsuit)$ in Steps 3, 4, and 5 
using the inductional hypothesis. 
 
Observe that 
$$(D_j(f), a - j -1) \in \fR\fD_E 
(G(\{(f^n,na -n)\})) 
\subset \bJ$$ 
for $0 \leq j \leq c - 1 = \lfloor a - 1 \rfloor$, 
from the inductional hypothesis. 
\end{case} 
\end{step} 
\begin{step}{Construction of a sequence $\{b_{u,j}\}$.} 
Our strategy for showing $(\clubsuit)$ is, 
starting from the following initial state 
$$ 
b_{0,j} = 1 \ (1 \leq j \leq c - 1), 
\quad 
b_{0,c} = a-c, 
\quad 
b_{u,0} = 0 \ (u\geq0), 
$$ 
to construct the (double) sequence of numbers $b_{u,j}$ indexed by 
$0 \leq j\leq c$ and $u\in\bZ_{\geq0}$ satisfying the following 
conditions 
$(\spadesuit)$ and $(\heartsuit)$: 
$$ 
(\spadesuit)\quad 
(D_j(f),a-j-b_{u,j}) \in \bJ, 
\qquad 
(\heartsuit)\quad 
\lim_{u \rightarrow \infty}b_{u,j} = 0.$$ 
We construct the numbers $b_{u,j}$ inductively according to the 
lexicographical order on the double index $(u,j)$. 
Suppose we have already constructed all $b_{\alpha,\beta}$ 
with $(\alpha,\beta) < (u,j)$.  Then, we define 
the number $b_{u,j}$ by the following formula 
$$nb_{u,j} = \max\left( 
\{0\}\cup\left\{ 
\sum_{t_l < j}b_{u,t_l} 
+ 
\sum_{j \leq t_l \leq c}b_{u-1,t_l} 
+ 
\sum_{c< t_l }(a-t_l) 
\mid T \in S_{n,j}^* 
\right\}\right).$$ 
where 
$S_{n,j} = \left\{T \in \bZ_{\geq 0}^n\mid |T| = nj\right\}$ 
and 
$S_{n,j}^* = S_{n,j} \setminus \{(j, \dotsc , j)\}$. 
\end{step} 
\begin{step}{Verification of $(\spadesuit)$.} 
By the argument in Step 2, condition $(\spadesuit)$ holds at the initial 
state, 
i.e., if $u=0$ or $j=0$. 
We proceed to check condition $(\clubsuit)$ by induction on 
the pair $(u,j)$ in the lexicographical order. 
Using the logarithmic version of the generalized product rule 
(cf.~Lemma \ref{1.2.2.2} (2)) 
for $f^n$, we compute 
$$(\sharp) \quad D_j(f)^n = D_{nj}(f^n) - 
\sum_{T \in 
S_{n,j}^*}\prod_{l = 1}^nD_{t_l}(f).$$ 
Take $T\in S_{n,j}^*$.  Then, by inductional hypothesis, we have 
$$ 
(\prod_{l = 1}^nD_{t_l}(f), 
\sum_{t_l < j}(a-t_l-b_{t_l,u})+ 
\sum_{j \leq t_l \leq c}(a-t_l-b_{t_l,u-1})) 
\in \bJ. 
$$ 
By definition of $b_{u,j}$ and the fact 
$\sum_{l}(a-t_l)=na-\sum_{l}t_l=na-nj$, we have 
\begin{eqnarray*} 
\lefteqn{ 
\sum_{t_l < j}(a-t_l-b_{t_l,u})+ 
\sum_{j \leq t_l \leq c}(a-t_l-b_{t_l,u-1})} 
\\ 
&=&\displaystyle{ 
\sum_{t_l \leq c}(a-t_l) 
-\sum_{t_l < j}b_{t_l,u} 
-\sum_{j \leq t_l \leq c}b_{t_l,u-1} 
}\\ 
&=&\displaystyle{ 
n(a-j)-\left(\sum_{t_l < j}b_{t_l,u} 
+\sum_{j \leq t_l \leq c}b_{t_l,u-1} 
+\sum_{t_l > c}(a-t_l) 
\right) 
\geq n(a-j)-nb_{u,j}. 
} 
\end{eqnarray*} 
On the other hand, by definition of 
$\bJ$, we have 
$D_{nj}(f^n)\in\bJ_{na-nj}$. 
Since $b_{u,j}\geq0$ by definition, 
we have $(D_{nj}(f^n),n(a-j)-nb_{u,j})\in\bJ$. 
Therfore, by virtue of the formula $(\sharp)$. 
we have $(D_j(f)^n,n(a-j)-nb_{u,j})\in\bJ$. 
As $\bJ$ is $\fR$-saturated, we have 
$(D_j(f),a-j-b_{u,j})\in\bJ$.  Thus 
$(\spadesuit)$ holds for $(u,j)$, 
as desired. 
\end{step} 
\begin{step}{Verification of $(\heartsuit)$.} 
We have only to show the following inequality: 
$$(\flat)\quad b_{u,j} \leq 
\left(1 - n^{-j}\right)\left(1 - n^{-c}\right)^{u-1} 
\quad(u\geq1,\ j\geq0). 
$$ 
In fact, since $b_{u,j}\geq0$ by definition 
and $0<1-n^{-m}<1$, condition $(\heartsuit)$ obviously follows from 
inequality $(\flat)$. 
 
We prove $(\flat)$ by induction on the pair $(u,j)$ in the lexicographical 
order. 
 
Since $b_{u,0}=0$, inequality $(\flat)$ is valid for $j=0$. 
 
By definition of $b_{u,j}$ and from the fact 
$\sum_{c<t_l}(a-t_l)<0$, we have an estimate 
$$nb_{u,j} \leq \max\left\{\sum_{t_l < j}b_{u,t_l} 
+ \sum_{j \leq t_l \leq c}b_{u-1,t_l}\mid T \in 
S_{n,j}^*\right\}.$$ 
By inductional hypothesis, we observe 
the following (i) and (ii): 
\begin{enumerate} 
\def\labelenumi{(\roman{enumi})} 
\item 
For $t_l < j$, we have 
$$b_{u,t_l} \leq \left(1 - n^{- t_l}\right)\left(1 - n^{-c}\right)^{u-1} 
\leq \left(1 - n^{1- j}\right)\left(1 - n^{-c}\right)^{u-1}$$ 
\item 
For $j \leq t_l \leq c$, we have 
$$b_{u-1,t_l} \leq 
\left(1 - n^{- t_l}\right)\left(1 - n^{-c}\right)^{u-2} 
\leq \left(1 - n^{-c}\right)^{u-1}.$$ 
\end{enumerate} 
We also mention that, for any $T = (t_1, \dotsc , t_n) \in S_{n,j}^*$, 
there exists at least one $1\leq l\leq n$ such that $t_l < j$. 
 
By these observations, we obtain the following estimate: 
\begin{eqnarray*} 
nb_{u,j} 
&\leq& \left(1 - n^{1- j}\right)\left(1 - n^{-c}\right)^{u-1} + (n - 
1)\left(1 - n^{-c}\right)^{u-1} \\ 
&=& \left(n - n^{1-j}\right)\left(1 - 
n^{-c}\right)^{u-1}  = n\left(1 - n^{-j}\right)\left(1 - 
n^{-c}\right)^{u-1}, 
\end{eqnarray*} 
which implies inequality$(\flat)$ for $(u,j)$. 
This completes the proof for inequality $(\flat)$, and hence the 
verification of $(\heartsuit)$. 
\end{step} 
\begin{step}{Finishing argument.} 
In the previous Steps, we confirmed conditions $(\spadesuit)$ and 
$(\heartsuit)$. 
Consequently, since $\bJ$ is $\fR$-saturated, 
we have $D_j(f)\in \bJ_{a-j}$ for $0 \leq j \leq c$. 
Namely $(\clubsuit)$ holds for $c = \lfloor a\rfloor$.  This completes the 
inductional proof of $(\clubsuit)$ stated in Step 2. 
\end{step} 
This completes the proof of Proposition \ref{2.2.3.1}. 
\end{proof} 
\begin{cor}\label{2.2.3.2} 
Let $\bI$ be an idealistic filtration 
over $R$ which has a regular system of 
parameters $(x_1, \dotsc, x_d)$.  Then 
$\fB(\bI)=\fR\fD(\bI)$. 
 
If $E$ is a simple normal crossing divisor defined by 
$\{x_1 \cdots x_m =0\}$ for some $1 \leq m \leq d$, then 
$\fB_E(\bI)=\fR\fD_E(\bI)$. 
\end{cor} 
\begin{proof} 
We present a proof of the latter assertion in the 
logarithmic case, as the former is a special case of the latter 
($E=\emptyset$). 
 
Since $\fB_E(\bI)$ is $\fD_E$-saturated, 
we have $\fB_E(\bI)\supset\fD_E(\bI)$. 
Then since $\fB_E(\bI)$ is $\fR$-saturated, 
we have $\fB_E(\bI)\supset\fR\fD_E(\bI)$. 
In order to see the opposite inclusion, we have only to show that 
$\fR\fD_E(\bI)$ 
is $\fD_E$-saturated. 
By Proposition \ref{2.2.3.1}, we see 
$$\fR\fD_E(\bI)\subset 
\fD_E\fR\fD_E(\bI) 
\subset\fR\fD_E\fD_E(\bI) 
=\fR\fD_E(\bI).$$ 
Therefore, we conclude that $\fR\fD_E(\bI)=\fD_E 
\fR\fD_E(\bI)$ is $\fD_E$-saturated. 
 
This completes the proof of Corollary \ref{2.2.3.2}. 
\end{proof} 
\end{subsection} 
\end{section} 
\begin{section}{Idealistic filtration of r.f.g.~type.}\label{2.3} 
In \ref{2.1} and \ref{2.2}, we gave the definition of, 
and carried discussion on the properties of, 
an idealistic filtration in general.  However, 
the idealistic filtrations we deal with 
in our algorithm are all of r.f.g.~type 
(cf.~Definition \ref{2.1.1.1} (4)).  In fact, certain mechanisms in 
our algorithm work only for the idealistic filtrations of r.f.g.~type. 
 
Since the operations of taking the $\fD$-saturation and 
$\fR$-saturation of a given idealistic filtration 
are essential in our algorithm, it is then a natural and important question 
whether the property of being of r.f.g.~type is stable under 
these operations.  The most important result of this section is to give an 
affirmative answer to this question: if an idealistic filtration 
$\bI$ is of r.f.g.~type, then so are 
$\fD(\bI)$ and $\fR(\bI)$.  We remark that some 
related results can be found in \cite{MR1996845}, 
discussing properties of an idealistic exponent. 
 
For $\fD$-saturation, the verification of stability is elementary, 
using compatibility of $\fD$-saturation with 
localization (cf.~Proposition \ref{Compatibility} (2)) 
and using the explicit construction in 
Lemma \ref{Construction}. 
 
For $\fR$-saturation, however, the verification of stability is rather 
subtle.  Our argument presented here is due to 
Professor Shigefumi Mori, who showed us how the contents of 
\cite{MR0089836} 
can be adapted to verify the required stability under $\fR$-saturation. 
The essential point, starting from a given idealistic filtration of 
r.f.g.~type $\bI$, is to show the rationality and boundedness of 
the denominators of the numbers $a$ where $\fR(\bI)_a$ changes. 
Once the crucial rationality and boundedness are shown, 
stability can be reinterpreted as the finite generation 
of the integral closure as an $R$-algebra (in some finite 
extension of the field of fractions) of a certain graded ring, 
which is naturally associated to the 
idealistic filtration $\bI$ of r.f.g.~type. 

In this section, $R$ denotes the coordinate ring of 
an affine open subset of a variety $W$ 
smooth over $k$ of $\on{char}(k) = p \geq 0$, 
or its localization by some multiplicative set. 
\begin{subsection}{Stability of r.f.g.~type under 
$\fD$-saturation.}\label{2.3.1} 
We show that the property of an idealistic 
filtration being of r.f.g.~type is stable under $\fD$-saturation. 
\begin{prop}\label{2.3.1.1} 
Let $\bI \subset R \times \bR$ be 
an idealistic filtration.  If $\bI$ is of r.f.g.~type, then 
so is its $\fD$-saturation $\fD(\bI)$ (or 
$\fD_E$-saturation $\fD_E(\bI)$). 
\end{prop} 
\begin{proof} 
\begin{step}{Reduction to the case where there 
exists a regular system of parameters 
$(x_1, \dotsc ,x_d)$ for $R$, where $d = \dim W$.} 
We take a finite affine cover $\{\on{Spec}R_{g_l}\mid g_l 
\in R\}_{l \in L}$ of $\on{Spec}R$ 
with $\# L < \infty$ so that for each $R_{g_l}$ there exists a regular 
system of parameters for $R_{g_l}$. 
 
Since $\bI$ is of r.f.g.~type, so is 
$\bI_{g_l}$, its localization by $g_l$. 
 
Suppose we have shown that $\fD(\bI_{g_l})$ is of r.f.g.~type, 
i.e., there exists a finite set 
$$T_{\Lambda_l} = 
\{(f_{\lambda_l},a_{\lambda_l})\}_{\lambda_l 
\in \Lambda_l} \subset R_{g_l} \times \bQ$$ 
such that $\fD(\bI_{g_l}) = G_{R_{g_l}}(T_{\Lambda_l})$. 
 
Observe that, since $\fD(\bI_{g_l}) = \fD(\bI)_{g_l}$ 
by compatibility of 
localization with $\fD$-saturation 
(cf.~Proposition \ref{Compatibility} (2)), for each 
$(f_{\lambda_l},a_{\lambda_l})$, there exist $(h_{\lambda},a_{\lambda}) \in 
\fD(\bI)$ and $n_{\lambda} \in \bZ_{> 0}$ such 
that $(f_{\lambda_l},a_{\lambda_l}) = 
(g_l^{-n_{\lambda_l}}h_{\lambda_l},a_{\lambda_l})$. 
 
Then it is easy to see that the finite set 
$$T_{\Lambda} = \{(h_{\lambda_l},a_{\lambda_l})\mid 
\lambda_l \in \Lambda_l, l \in L\} \subset 
\fD(\bI)$$ 
generates $\fD(\bI)$, i.e., $\fD(\bI) = 
G_R(T_{\Lambda})$. 
In fact, by construction, we have 
$$\fD(\bI)_{g_l} \supset G_R(T_{\Lambda})_{g_l} \supset 
G_{R_{g_l}}(T_{{\Lambda}_l}) = 
\fD(\bI_{g_l}) = \fD(\bI)_{g_l},$$ 
i.e., $\fD(\bI)_{g_l} = G_R(T_{\Lambda})_{g_l}$ for any 
$l\in L$, and hence $\fD(\bI) = G_R(T_{\Lambda})$. 
\end{step} 
\begin{step}{%
Proof of the statement in the case where there exists a regular 
system of parameters $(x_1,\dotsc,x_d)$ for $R$, where 
$d = \dim W$.} 
Take a finite set of generators $T_{\Lambda}$ of the form 
$$T_{\Lambda} = 
\{(f_{\lambda},a_{\lambda})\}_{\lambda 
\in \Lambda} \subset R \times \bQ$$ 
such that $\bI = G(T_{\Lambda})$.  We may assume $a_{\lambda} > 0 
\quad \forall \lambda \in \Lambda$ by discarding those with 
$a_{\lambda} \leq 0$.  Let 
$$T_M = 
\{(\partial_{X^J}f_{\lambda},a_{\lambda} - |J|)\mid 
(f_{\lambda}, a_{\lambda}) \in T_{\Lambda}, 0 \leq |J| < 
a_{\lambda}\}.$$ Then clearly we have $\# T_M < \infty$ and $a_{\lambda} - 
|J| \in \bQ \quad \forall \lambda$ and $\forall J 
\text{\ with\ }0 \leq |J| < a_{\lambda}$. 
 
Now it follows from Lemma \ref{Construction} (2) 
that $\fD(\bI) = G(T_M)$. 
Therefore, we conclude that $\fD(\bI)$ is of r.f.g.~type. 
 
The proof for stability under $\fD_E$-saturation is identical. 
\end{step} 
This completes the proof for Proposition \ref{2.3.1.1}. 
\end{proof} 
\end{subsection} 
\begin{subsection}{Stability under $\fR$-saturation.}\label{2.3.2} 
We show that the property of an idealistic filtration being of 
r.f.g.~type is stable under $\fR$-saturation.  We deal with the 
problem of stability in terms of a certain graded ring which is 
naturally associated to an idealistic filtration $\bI$ of r.f.g.~type 
and which ``describes'' $\bI$ in the sense stated below. 
\begin{defn}\label{2.3.2.0} 
Let $A = \bigoplus_{n \in \bZ_{\geq 0}}A_{qn}X^{qn} \subset 
\bigoplus_{n \in \bZ_{\geq 0}}RX^{qn} = R[X^q]$ be a 
graded $R$-subalgebra of the polynomial ring with one variable 
$X^q$ over $R$ for some 
$q \in \bQ_{> 0}$.  Let $\bI \subset R \times \bR$ be an 
idealistic filtration.  We say $A$ describes $\bI$ if it 
satisfies the following condition: 
$$\bI_{qa} = A_{q\lceil a\rceil} 
\text{\ for\ any\ }a \in\bR_{\geq 0}.$$ 
\end{defn} 
\begin{lem}\label{2.3.2.1} 
Let $\bI \subset R \times \bR$ be an 
idealistic filtration.  Then $\bI$ is of r.f.g.~type if and 
only if there exists $A$ which describes $\bI$ (as stated in 
Definition \ref{2.3.2.0}) and which is finitely generated as an $R$-algebra. 
\end{lem} 
\begin{proof} 
Suppose that there exists such $A$ which describes 
$\bI$ and which is generated by a finite set of homogeneous 
elements 
$\{f_{\lambda}X^{qn_{\lambda}}\}_{\lambda \in 
\Lambda}$ as a graded 
$R$-subalgebra in $R[X^q]$.  Then $\bI$ is generated 
by 
the finite set $\{(f_{\lambda},qn_{\lambda})\}_{\lambda 
\in 
\Lambda}$, and hence is of r.f.g.~type. 
 
Conversely, suppose that $\bI$ is an idealistic filtration of 
r.f.g.~type, generated by a finite set $T = 
\{(f_{\lambda},\frac{n_{\lambda}}{\delta})\}_{\lambda 
\in \Lambda} \subset R \times \bQ$ 
for some $\delta \in \bZ_{> 0}$.  It is immediate that, if we take the 
graded $R$-subalgebra $A$ of $R[X^q]$, with 
$q = \delta^{-1}$ and $A_0 = R$, generated by the finite set 
$\{f_{\lambda}X^{\frac{i}{\delta}} 
\mid \lambda \in \Lambda, 0 \leq i \leq n_{\lambda}\}$ over $R$, 
then $A$ describes $\bI$. 
 
This completes the proof of Lemma \ref{2.3.2.1}. 
\end{proof} 
 
We remark that if $\bI$ is an idealistic filtration of 
r.f.g.~type, and if $A \subset R[X^q]$ is a graded 
$R$-subalgebra which describes 
$\bI$ for some $q \in \bQ_{> 0}$, then $A$ is automatically 
finitely generated over $R$. 
\begin{cor}\label{2.3.2.2} 
Let $\bI$ be an idealistic filtration of 
r.f.g.~type.  Then $\bI$ satisfies condition (continuity). 
\end{cor} 
\begin{proof} 
We want to show $\bI_a = \bigcap_{b < a}\bI_b$ 
for any $a \in \bR$ (cf.~Remark \ref{2.1.3.2} (2)). 
 
It is clear when $a \leq 0$ 
(cf.~condition (o) in Definition \ref{2.1.1.1} (2)). 
 
Suppose $a > 0$.  By Lemma \ref{2.3.2.1}, 
there exists a graded $R$-subalgebra $A 
\subset R[X^q]$, for some $q \in \bQ_{> 0}$, 
which describes $\bI$ and which is finitely generated as an 
$R$-algebra.  Then by definition we have 
$\bI_{qa} = A_{q\lceil a\rceil}$.  Since 
$\bI_s \supset \bI_t$ for any $s < t$, we conclude 
$$\bigcap_{b < a}\bI_{qb} 
= \bigcap_{\lceil a\rceil-1<b<a}\bI_{qb} 
= \bigcap_{\lceil a\rceil-1<b<a}A_{q\lceil b\rceil} 
= \bigcap_{\lceil a\rceil-1<b<a}A_{q\lceil a\rceil} 
= A_{q\lceil a\rceil} = \bI_{qa},$$ 
i.e., $\displaystyle \bI_{qa} = \bigcap_{b < a}\bI_{qb}$. 
 
This completes the proof of Corollary \ref{2.3.2.2}. 
\end{proof} 
\begin{prop}\label{2.3.2.3} 
Let $\bI \subset R \times \bR$ be 
an idealistic filtration.  If $\bI$ is of r.f.g.~type, then 
so is its $\fR$-saturation. 
\end{prop} 
Before beginning the proof of Proposition \ref{2.3.2.3}, 
we extract the essence that we need from 
Nagata's paper \cite{MR0089836} 
with some modifications. 

Let $R$ be a noetherian domain, $K=Q(R)$ its field of fractions and 
${\mathfrak a}=(u_1,\dotsc,u_s)\subset R$ a proper ideal of $R$ 
with a finite set of its generators $u_j$.   Set 
$R_j=R\left[\frac{u_1}{u_j},\dotsc,\frac{u_s}{u_j}\right]$, and let 
$\overline{R_j}$ be its normalization 
in $Q(R_j) = K$  for each $j$.  Let 
$\{P_{jk}\}_k\subset\on{Spec}\overline{R_j}$ 
be the set of all the minimal primes of $u_j\overline{R_j}$.  Note that it 
is a finite set and that, by Krull's Hauptidealsatz, the primes 
$P_{jk}$ are of height $1$.   Let $R_{jk}=(\overline{R_j})_{P_{jk}}$ be the 
localization of $\overline{R_j}$ at $P_{jk}$ for each $j,k$. 
Since $R_{jk}$ is a $1$-dimensional noetherian normal ring, 
it is a discrete valuation ring.  We denote the valuation 
of $R_{jk}$ by $v_{jk}$ for each $j,k$.  We consider 
the functions 
$\theta_{\mathfrak a},\overline\theta_{\mathfrak a}\colon 
R\rightarrow\bR_{\geq 0} \cup\{\infty\}$ 
defined by 
\begin{eqnarray*} 
\theta_{\mathfrak a}(r) &=& \sup \left\{\frac{m}{n}\mid  r^n 
\in {\mathfrak a}^m, n, m \in \bZ_{\geq0}, n > 0\right\},\\ 
\overline\theta_{\mathfrak a}(r) &=& \sup 
\left\{\frac{m}{n}\mid r^n \in \overline{{\mathfrak a}^m}, n, 
m \in \bZ_{\geq0}, n > 0\right\}. 
\end{eqnarray*} 
Using the notation as above, we have the following lemmas. 
\begin{lem}\label{2.3.2.4} 
For $n\in\bZ_{>0}$, we have 
$$ 
\overline{{\mathfrak a}^n}=R\cap\bigcap_{j,k}u_j^nR_{jk}. 
$$ 
\end{lem} 
\begin{proof} 
Firstly we show 
$\overline{{\mathfrak a}^n}\subset R\cap\bigcap_{j,k}u_j^nR_{jk}$. 
It suffices to show 
$\overline{{\mathfrak a}^n}\subset u_j^nR_{jk}$ for each $j,k$. 
Fix $j,k$ and take $f\in\overline{{\mathfrak a}^n}$.  Then, there exists a 
monic equation 
$$ 
f^m+a_1f^{m-1}+\cdots+a_m=0\quad(a_i\in{\mathfrak a}^{in}). 
$$ 
Considering the valuation $v_{jk}$ of this equation, 
there exists some $1\leq i\leq m$ such that 
$v_{jk}(f^m)=v_{jk}(a_if^{m-i})$ and hence $v_{jk}(f^i) = v_{jk}(a_i)$. 
Since ${\mathfrak a} R_{jk}=u_j^nR_{jk}$, we have 
$v_{jk}(a_i)\geq in \cdot v_{jk}(u_j)$. 
Consequently $v_{jk}(f)\geq nv_{jk}(u_j)$, 
and hence $f\in u_j^nR_{jk}$.  Thus 
the inclusion $\overline{{\mathfrak a}^n}\subset u_j^nR_{jk}$ holds. 
 
Secondly we show the opposite inclusion $\overline{{\mathfrak a}^n}\supset 
R\cap\bigcap_{j,k}u_j^nR_{jk}$. 
 
Take $g\in R\cap\bigcap_{j,k}u_j^nR_{jk}$.  Set 
$R'=R\left[\frac{u_1^n}g,\dotsc,\frac{u_s^n}g\right]$ and 
${\mathfrak b}=\left( 
\frac{u_1^n}g,\dotsc,\frac{u_s^n}g 
\right)\subset R'$. 
We show $g\in\overline{{\mathfrak a}^n}$ in the following Steps. 
\begin{step}{We show ${\mathfrak b}=R'$.} 
Assume ${\mathfrak b}$ is a proper ideal of $R'$. 
Then there exists a valuation ring $(V,\fm_V)$ of $Q(R')=K$ 
such that $V\supset R'$ and $\fm_V\cap R'\supset{\mathfrak b}$. 
We denote its valuation as $v$. 
Take $j_0$ such that $v(u_{j_0})=\min_{1\leq i\leq s}v(u_i)$. 
Then, as $\frac{u_i}{u_{j_0}}\in V$ for each $i$, 
we have 
$$ 
R_{j_0}\subset V 
\quad\text{and hence }\quad 
\overline{R_{j_0}}\subset V 
$$ 
Since $\frac{u_{j_0}^n}g\in{\mathfrak b}\subset\fm_v$, we have 
$g\not\in u_{j_0}^nV$, and hence 
$g\not\in u_{j_0}^n\overline{R_{j_0}}$. 
Now, since $\overline{R_{j_0}}$ is noetherian normal domain, 
principal ideal $u_{j_0}^n\overline{R_{j_0}}$ is represented as 
$$ 
u_{j_0}^n\overline{R_{j_0}}= 
\overline{R_{j_0}} 
\cap\bigcap_{\on{ht}{\mathfrak p}=1} 
{\mathfrak p}^{v_{\mathfrak p}(u_{j_0}^n)}(\overline{R_{j_0}})_{\mathfrak p} 
=\overline{R_{j_0}}\cap\bigcap_k 
P_{j_0k}^{v_{j_0k}(u_{j_0}^n)}R_{j_0k}. 
$$ 
Therefore there exists some $k$ such that 
$$g\not\in P_{j_0k}^{v_{j_0k} 
(u_{j_0}^n)}R_{j_0k}=u_{j_0}^nR_{j_0k},$$ 
which contradicts to the choice of $g$. 
Thus we have ${\mathfrak b}=R'$. 
\end{step} 
\begin{step}{We show $g\in\overline{{\mathfrak a}^n}$.} 
Since $1\in{\mathfrak b}$ by Step 1, there exists 
$F(X_1,\dotsc,X_s)\in R[X_1,\dotsc,X_s]$ such that $F(0,\dotsc,0)=0$ and 
$F\left(\frac{u_1^n}g,\dotsc, 
\frac{u_s^n}g\right)=1$.  Setting 
$\deg F = n$, we obatin 
$$0 = g^n\left\{1 - F\left(\frac{u_1^n}g,\dotsc, 
\frac{u_s^n}g\right)\right\} = g^n + c_1g^{n-1} +\cdots+ c_n \text{\ with\ 
}c_i \in {\mathfrak a}^i,$$ 
a monic equation which shows $g \in \overline{{\mathfrak a}^n}$. 
\end{step} 
This completes the proof of Lemma \ref{2.3.2.4}. 
\end{proof} 
\begin{lem}\label{2.3.2.5} 
Let $r\in R$.  Then, 
$$\theta_{\mathfrak a}(r)=\overline\theta_{\mathfrak a}(r)=\min_{j,k} 
\left\{\frac{v_{jk}(r)}{v_{jk}(u_j)}\right\} \in \bQ.$$ 
Moreover, for $n,m\in\bZ_{\geq0}$ with $n>0$, 
$r^n\in\overline{{\mathfrak a}^m}$ if and only if 
$\frac{m}n\leq\overline\theta_{\mathfrak a}(r)$. 
\end{lem} 
\begin{proof} 
\begin{step}{We show the first equation 
$\theta_{\mathfrak a}(r)=\overline\theta_{\mathfrak a}(r)$.} 
Since ${\mathfrak a}^m\subset\overline{{\mathfrak a}^m}$, 
it is immediate that 
$\theta_{\mathfrak a}(r)\leq\overline\theta_{\mathfrak a}(r)$. 
We show $\theta_{\mathfrak a}(r)\geq\overline\theta_{\mathfrak a}(r)$. 
Take $n,m\in\bZ_{\geq0}$ with $n>0$ such that 
$r^n\in\overline{{\mathfrak a}^m}$.  By definition, there exists 
a monic equation 
$$ 
(r^n)^{c+1}+a_1(r^n)^c+\cdots+a_{c+1}=0 
\quad\text{with}\quad a_i\in{\mathfrak a}^{im}. 
$$ 
We show 
$$\heartsuit_t\colon\quad 
r^{n(c+t)}\in{\mathfrak a}^{mt}\left( 
r^nR+{\mathfrak a}^m\right)^c 
\qquad(t\in\bZ_{>0})$$ 
by induction on $t$. 
Looking at the monic equation above, we have 
$$ 
r^{n(c+1)}\in 
{\mathfrak a}^{m}(r^n)^c+\cdots+{\mathfrak a}^{(c+1)m}r^{n(c+t)} 
={\mathfrak a}^m\left(r^nR+{\mathfrak a}^m\right)^c, 
$$ 
thus $\heartsuit_1$ holds. 
For the case $t>1$, we have 
\begin{eqnarray*} 
r^{n(c+t)}&=&r^n\cdot r^{n(c+t-1)} 
\in r^n{\mathfrak a}^{m(t-1)}(r^nR+{\mathfrak a}^m)^c 
\qquad(\text{By }\heartsuit_{t-1}) 
\\&\subset& 
{\mathfrak a}^{m(t-1)}\left( 
r^{n(c+1)}R+{\mathfrak a}^m(r^nR+{\mathfrak a}^m)^c 
\right) 
\\&\subset& 
{\mathfrak a}^{m(t-1)}\left( 
{\mathfrak a}^m(r^nR+{\mathfrak a}^m)^c 
\right) 
\qquad(\text{By }\heartsuit_1) 
\\&=& 
{\mathfrak a}^{mt}(r^nR+{\mathfrak a}^m)^c. 
\end{eqnarray*} 
Thus $\heartsuit_t$, 
and hence 
$r^{n(c+t)}\in{\mathfrak a}^{mt}$ 
holds for any $t\in\bZ_{>0}$. 
It follows that 
$$ 
\theta_{\mathfrak a}(r) 
\geq\sup\left\{\frac{mt}{n(c+t)}\mid t\in 
\bZ_{>0}\right\}\geq\frac{m}n. 
$$ 
Since the numbers $n,m \in \bZ_{\geq 0}$ with $n > 0$ such that 
$r^n\in\overline{{\mathfrak a}^m}$ are taken arbitrarily, we have 
$\theta_{\mathfrak a}(r)\geq\overline\theta_{\mathfrak a}(r)$. 
\end{step} 
\begin{step}{We show the second equality.} 
By Lemma \ref{2.3.2.4}, we have 
\begin{eqnarray*} 
r^n\in\overline{{\mathfrak a}^m} 
\quad&\Longleftrightarrow&\quad 
r^n\in u_j^mR_{jk} 
\quad(\forall j,k) 
\quad\Longleftrightarrow\quad 
v_{jk}(r^n)\geq v_{jk}(u_j^m) 
\quad(\forall j,k) 
\\&\Longleftrightarrow&\quad 
nv_{jk}(r)\geq mv_{jk}(u_j) 
\quad(\forall j,k) 
\\&\Longleftrightarrow&\quad 
\frac{m}n\leq\min\left\{\frac{v_{jk}(r)}{v_{jk}(u_j)} 
\mid j,k\right\} 
\end{eqnarray*} 
Therefore 
$\overline\theta_{\mathfrak a}(r)= 
\min\left\{\frac{v_{jk}(r)}{v_{jk}(u_j)} 
\mid j,k\right\} \in \bQ$. 
\end{step} 
The ``Moreover'' part is now obvious. 
 
This completes the proof of Lemma \ref{2.3.2.5}. 
\end{proof} 
We now go back to the proof of Proposition \ref{2.3.2.3}. 
\begin{proof}[Proof of Proposition \ref{2.3.2.3}] 
Take a finite set 
$T = \{(f_{\lambda},a_{\lambda})\}_{\lambda 
\in \Lambda} \subset R \times \bQ$ 
such that $\bI = G(T)$. 
\begin{step}{%
We may assume $T \subset R \times \{L\}$ for some 
$L \in \bZ_{> 0}$.} 
Replacing $T$ with $T \cap R \times \bR_{> 0}$, we may assume $T 
\subset R \times \bQ_{> 0}$.  Set 
$$L = \min\left\{n \in \bZ_{> 0}\mid 
\frac{n}{a_{\lambda}} \in \bZ_{> 
0} \quad \forall \lambda \in \Lambda\right\} 
\quad\text{and}\quad 
T' = \left\{(f_{\lambda}^{\frac{L}{a_{\lambda}}},L) 
\right\}_{\lambda\in\Lambda}.$$ 
Then it is clear that $\fR(G(T)) = \fR(G(T'))$.  Therefore, by 
replacing $T$ with $T'$, we may assume $T \subset R \times \{L\}$. 
\end{step} 
\begin{step}{Description of $\fR(\bI)$ 
in terms of the function $\theta_I$.} 
Let $I=\bI_L$ be the ideal of the idealistic filtration $\bI$ at 
level $L$.  Define $\bJ\subset R\times\bR$ by setting 
$\bJ_{La}=\{f\in R\mid \theta_I(f)\geq a\}$ for $a\in\bR$. 
We show $\fR(\bI)=\bJ$. 
Since $\bI=G(T)=G(I\times\{L\})$, we have 
$\bI_{La}=I^{\lceil a\rceil}$ for any $a\in\bR$ 
by Lemma \ref{Construction} (1). 
(We use the convention that $I^{-n} = R$ for $n \in 
\bZ_{> 0}$). 
Thus, by Lemma \ref{Construction} (3), 
$\fR(\bI)=\overline\bK$ where 
$\overline\bK\subset R\times\bR$ is defined by 
$$ 
\overline\bK_{La} 
=\left\{f\in R\mid \forall b<a,\ \exists n\in\bZ_{>0} 
\ \text{s.t.}\ f^n\in \bI_{nLb} = I^{\lceil nb\rceil}\right\} 
\quad(a\in\bR). 
$$ 
The condition above can be rephrased as follows: 
\begin{eqnarray*} 
\lefteqn{ 
\left(\forall b<a,\ \exists n\in\bZ_{>0}\ 
\text{s.t.}\ f^n\in I^{\lceil nb\rceil}\right) 
}\\ 
&\Leftrightarrow& 
\left(\sup\left\{b\in\bR_{\geq0}\mid 
\exists n\in\bZ_{>0}\ 
\text{s.t.}\ f^n\in I^{\lceil nb\rceil} 
\right\}\geq a\right)\\ 
&\Leftrightarrow& 
\left(\sup\left\{ 
\frac{\lceil nb\rceil}n\mid 
\exists n\in\bZ_{>0},\ 
\exists b\in\bR_{\geq0}\ 
\text{s.t.}\ f^n\in I^{\lceil nb\rceil} 
\right\}\geq a\right)\\ 
&\Leftrightarrow& 
\left(\sup\left\{\frac{m}n\mid 
\exists n, m \in\bZ_{>0} \text{\ with\ }n > 0\ \text{s.t.}\ f^n\in I^m 
\right\}\geq a\right)\\ 
&\Leftrightarrow&\theta_I(f)\geq a 
\end{eqnarray*} 
Thus $\fR(\bI)_{La}=\overline\bK_{La}=\bJ_{La}$ for 
$a\in\bR$, hence $\fR(\bI)=\bJ$. 
\end{step} 
\begin{step}{There exists $\rho\in\bZ_{>0}$ such that 
$\bJ_a=\bJ_{{\lceil\rho a\rceil}/\rho}$ for any $a\in\bR$.} 
We apply Lemma \ref{2.3.2.5} with ${\mathfrak a}=I$ to 
our setting. Let $\rho$ be a common multiple of 
$\{v_{jk}(u_j)\mid j,k\}$. 
Take $f\in\bJ_{La}$. 
Then, we have 
$\theta_I(f)\geq a$. 
Since $\rho\theta_I(f)\in\bZ$ 
by Lemma \ref{2.3.2.5}, we 
have $\rho L\theta_I(f)\geq\lceil\rho La\rceil$. 
Therefore, we have 
$f\in 
\bJ_{{\lceil\rho La\rceil}/\rho}$, 
and hence 
$\bJ_{La}\subset\bJ_{{\lceil\rho(La)\rceil}/\rho}$. 
The opposite inclusion is clear by condition (iii) in 
Definition \ref{2.1.1.1} for the idealistic filtration $\bJ$. 
\end{step} 
\begin{step}{We show $S_1$ describes $\fR(\bI)$ and 
$S_1=\overline{S_0}^{R_1}$ in the following notation:} 
Consider the graded $R$-algebras 
 \begin{alignat*}{5} 
R_0&= 
R[X^L] 
&\quad&\supset&\quad 
S_0&= 
\bigoplus_{n\in\bZ_{\geq0}}I^nX^{Ln} 
\\ 
R_1&= 
R[X^{\frac1\rho}] 
&\quad&\supset&\quad 
S_1&= 
\bigoplus_{n\in\bZ_{\geq0}}\bJ_{\frac{n}\rho}X^{\frac{n}\rho} 
\end{alignat*} 
where $X$ is an indeterminate. 
We denote by $\overline{S_0}^{R_1}$ the 
integral closure of $S_0$ in $R_1$. 
 
It is clear from Step 3 that $S_1$ describes $\bJ = \fR(\bI)$. 
We have only to prove $S_1=\overline{S_0}^{R_1}$. 
 
Firstly we show $S_1 \subset \overline{S_0}^{R_1}$. 
Let $gX^{\frac{n}\rho}\in S_1$ be a homogeneous 
element of $S_1$.  Since $g\in J_{\frac{n}\rho}$, we have 
$\overline\theta_I(g)=\theta_I(g)\geq\frac{n}{\rho L}$. 
Thus, by Lemma \ref{2.3.2.5}, we have $g^{\rho L}\in\overline{I^n}$. 
Therefore there exists a monic equation 
$$ 
\left(g^{\rho L}\right)^m 
+ c_1\left(g^{\rho L}\right)^{m-1} 
+\cdots+ c_m = 0 \text{\ with\ }c_i \in 
(I^n)^i.$$ 
This in turn provides a monic equation of $gX^{\frac{n}\rho}$ 
over $S_0$, i.e., 
$$\left(gX^{\frac{n}\rho}\right)^{\rho Lm} 
+ (c_1X^{Ln})\left(gX^{\frac{n}\rho}\right)^{\rho L(m-1)} 
+\cdots+ c_mX^{Lnm} = 0.$$ 
Therefore, we have $S_1 \subset \overline{S_0}^{R_1}$. 
 
Secondly we show $S_1 \supset \overline{S_0}^{R_1}$. 
Take 
$g=\sum_{n\in\bZ_{\geq0}}g_{\frac{n}\rho} 
X^{\frac{n}\rho}\in\overline{S_0}^{R_1}$. 
Then we have a monic equation of $g$ over $S_0$, i.e., 
$$ 
(\spadesuit)\quad 
g^m+c_1(X^L)g^{m-1}+\cdots+ c_m(X^L) = 0 
\quad\text{with}\quad c_i(X^L)\in S_0.$$ 
Set $G=\sum_{n\in\bZ_{\geq0}}g_{\frac{n}\rho} 
X^{\frac{n}\rho}Y^n\in{R_1}[Y]$ where $Y$ is another indeterminate. 
By replacing $X$ by $XY^\rho$ in $(\spadesuit)$, 
we have a monic equation of $G$ over $S_0[Y]$, i.e., 
$$G^m+c_1(X^LY^{\rho L})G^{m-1} 
+\cdots+ c_m(X^LY^{\rho L}) = 0 
\quad\text{with}\quad c_i(X^LY^{\rho L})\in S_0[Y].$$ 
Since 
$\overline{S_0[Y]}^{R_1[Y]}=\overline{S_0}^{R_1}[Y]$ 
(cf.~{\it Alg. Comm.}, chap.~V, \S1, n$^{\on o}$3, prop.12 
in \cite{MR0194450}), each coefficient of $Y^n$ in $G$ are 
integral over $S_0$, i.e., 
$$g_{\frac{n}\rho}X^{\frac{n}\rho}\in\overline{S_0}^{R_1} 
\quad(n\in\bZ_{\geq0}).$$ 
Thus we may assume $g$ is a homogeneous element in $R_1$, say, 
$g=g_{\frac{l}\rho}X^{\frac{l}\rho}$. 
Looking at the coefficient of $X^{\frac{ml}\rho}$ 
in $(\spadesuit)$, we have 
$$ 
g_{\frac{l}\rho}^m 
+\alpha_1g_{\frac{l}\rho}^{m-1} 
+\cdots+\alpha_m=0 
$$ 
where $\alpha_n$ is the coefficient of $X^{\frac{nl}\rho}$ 
in $c_n\in S_0\subset R_1$.  Note that 
$\alpha_n=0$ if $nl\not\in\rho L\bZ$, 
and $\alpha_n\in I^{\frac{nl}{\rho L}}$ if 
${nl}\in\rho L\bZ$.  Thus, for any 
$1\leq n\leq m$, we have 
$$\alpha_n\in I^{\frac{nl}{\rho L}}=\bI_{\frac{nl}\rho} 
\subset\fR(\bI)_{\frac{nl}\rho}. 
$$ 
 
Since $\fR(\bI)$ is integrally closed 
by Proposition \ref{2.2.2.1}, we have 
$$g_{\frac{l}\rho}\in\fR 
(\bI)_{\frac{l}\rho} 
=\bJ_{\frac{l}\rho} 
\quad\text{and hence}\quad 
g=g_{\frac{l}\rho}X^{\frac{l}\rho}\in 
\bJ_{\frac{l}\rho}X^{\frac{l}\rho} 
\subset S_1.$$ 
Therefore, we have 
$S_1\supset\overline{S_0}^{R_1}$. 
\end{step} 
\begin{step}{We see that $S_1$ is finitely generated over $R$.} 
It is clear when $I=(0)$, since $S_1=R$.  We assume $I\neq(0)$.  

Since $R$ is normal, so is $R_1=R[X^{\frac1\rho}]$. Thus 
$$ 
S_1=\overline{S_0}^{R_1} 
=\overline{S_0}^{Q(R_1)}. 
$$ 
Note that $Q(R_1)$ is a finite extension of 
$Q(S_0)=Q(R[X^L])$.  By \S 33 of \cite{MR879273}, 
it follows that 
$S_1=\overline{S_0}^{Q(R_1)}$ is 
a finite $S_0$-module.  On the other hand, 
$S_0$ is finitely generated over $R$. 
Indeed, taking generators of $I$ as $I=(r_1,\dotsc,r_t)$, 
we have $S_0=R[r_1X^L,\dotsc,r_tX^L]$. 
Thus $S_1$ is also finitely generated over $R$. 
\end{step} 
\begin{step}{Finishing argument.} 
By Steps 2 and 3, we see that 
$S_1$ describes 
the idealistic filtration $\fR(\bI)$. 
Since $S_1$ is 
finitely generated over $R$, 
we conclude that $\fR(\bI)$ is r.f.g.~type 
(cf.~Lemma \ref{2.3.2.1}). 
\end{step} 
This completes the proof of Proposition \ref{2.3.2.3}. 
\end{proof} 
\begin{cor}\label{2.3.2.6} 
Let $\bI \subset R \times \bR$ be an 
idealistic filtration.  Assume $\bI$ is of r.f.g.~type. 
Then its $\fR$-saturation coincides with its integral closure, i.e., 
$$\fR(\bI) = \on{IC}(\bI).$$ 
\end{cor} 
\begin{proof} 
By Proposition \ref{2.2.2.1}, $\fR(\bI)$ is 
integrally closed. 
 
Therefore, we have $\fR(\bI) 
\supset \on{IC}(\bI)$.  Thus we have only to show 
$\fR(\bI) \subset \on{IC}(\bI)$. 
 
By the same argument as in Step 1 of the proof of 
Proposition \ref{2.3.2.3}, we 
may assume that $\bI$ is generated by a finite number of 
elements at level $L$.  In fact, using the same notation, we see 
that $\fR(G(T)) = \fR(G(T'))$ and 
$\on{IC}(G(T)) =\on{IC}(G(T'))$.  Then as shown 
in Step 4 of the proof of Proposition \ref{2.3.2.3}, 
the integral closure 
$S_1 \subset R[X^{\frac{1}{\rho}}]$ of $S_0$ in 
$R[X^{\frac{1}{\rho}}]$ describes 
$\fR(\bI)$, while $S_0$ describes $\bI$. 
 
Take an element $(f,a) \in \fR(\bI)$.  Then we have 
$(f,\frac{\lceil \rho a\rceil}{\rho}) \in \fR(\bI)$ 
(cf.~Lemma \ref{2.3.1.1}), which implies 
$fX^{\frac{\lceil \rho a\rceil}{\rho}} \in S_1$.  Now 
since $fX^{\frac{\lceil \rho a\rceil}{\rho}}$ is 
integral over $S_0$, 
by the same argument as in Step 4 of the proof of 
Proposition \ref{2.3.2.3}, 
we see that $(f,\frac{\lceil \rho a\rceil}{\rho})$ is integral 
over $\bI$, i.e., $(f,\frac{\lceil \rho a\rceil}{\rho}) \in 
\on{IC}(\bI)$.  Finally, since $a \leq \frac{\lceil \rho 
a\rceil}{\rho}$, we conclude $(f,a) \in \on{IC}(\bI)$.  This shows 
the desired inclusion. 
 
This completes the proof of Corollary \ref{2.3.2.6}. 
\end{proof} 
\end{subsection} 
\end{section} 
\begin{section}{Localization and completion of 
an idealistic filtration.}\label{2.4} 
In this section, we discuss the notion of localization and completion of an 
idealistic filtration over $R$, associated to the 
localization and completion of $R$, respectively.  Our main observation here 
is the compatibility of the operations of taking the 
generation, 
$\fD$-saturation, and $\fR$-saturation with localization and 
completion.  The compatibility allows us to reduce the analysis 
of the global properties of these operations to the local or to the analytic 
ones, to which we may apply some explicit computations. 

In this section $R$ denotes the coordinate ring of an affine open subset of 
a nonsingular variety $W$ over $k$. 
\begin{subsection}{Definition.}\label{2.4.1} 
\begin{defn}\label{2.4.1.1} 
Let $\bI \subset R \times \bR$ be an 
idealistic filtration over $R$. 
\item[(1)] 
(Localization)\quad 
Let $S$ be a multiplicative 
set of $R$.  Consider the 
subset $\bI_S \subset R_S \times \bR$ 
defined by 
$$(\bI_S)_a = (\bI_a)_S = \bI_a \otimes_R R_S\quad (a 
\in \bR).$$ 
Then $\bI_S$ is an idealistic filtration, called the localization of 
$\bI$ by $S$. 
 
In case $P \in \on{Spec}R$ is a point corresponding to 
a prime ideal $P \subset R$ 
(we use the same symbol for the point and prime 
ideal by abuse of notation) with $S = R \setminus P$, 
we often denote $\bI_S$ by $\bI_P$. 
\item[(2)] 
(Completion)\quad 
Let $\widehat{R}$ be the completion of $R$ with respect to 
a maximal ideal $\fm \subset R$.  Consider the 
subset $\widehat\bI \subset \widehat{R} \times \bR$ defined by 
$$(\widehat\bI)_a = \widehat{\bI_a} = \bI_a \otimes_R 
\widehat{R} \quad (a \in \bR).$$ 
Then $\widehat\bI$ is an idealistic filtration, called the completion 
of $\bI$ (with respect to $\fm$-adic topology). 
\end{defn} 
\begin{rem}\label{2.4.1.2} 
We remark that, for idealistic filtrations $\bI, \bI' \subset R 
\times \bR$, the following conditions 
are equivalent: 
\begin{enumerate} 
\item $\bI \subset \bI'$, 
\item $\bI_{\fm} \subset \bI'_{\fm}$ for any maximal ideal $\fm 
\subset R$, 
\item $\widehat\bI \subset \widehat{\bI'}$, where the completion 
``$\widehat{\phantom{m}}$'' 
is taken with respect the $\fm$-adic topology, 
for any maximal ideal $\fm \subset R$. 
\end{enumerate} 
In fact, fixing the level $a \in \bR$, we see that the equivalence of 
the conditions on the idealistic filtrations follows from 
the equivalence of the corresponding conditions on the ideals, 
which is a standard result in commutative ring theory. 
\end{rem} 
\end{subsection} 
\begin{subsection}{Compatibility.}\label{2.4.2} 
\begin{prop}\label{Compatibility} 
\item[(1)] 
(Compatibility with generation) Let $T \subset R \times \bR$ be a 
subset.  Then we have 
$$G_R(T)_S = G_{R_S}(T), \quad 
\widehat{G_R(T)} = G_{\widehat{R}}(T).$$ 
In particular, if $\bI = G(T)$ is of r.f.g.~type, then so are 
$\bI_S$ and $\widehat\bI$. 
\item[(2)] 
(Compatibility with $\fD$-saturation) Let $\bI \subset R 
\times \bR$ be an idealistic filtration.  Then we have 
$$\fD(\bI)_S = \fD(\bI_S), \quad 
\widehat{\fD(\bI)} = \fD(\widehat\bI).$$ 
Let $E$ be a simple normal crossing divisor.  Then we have 
$$\fD_E(\bI)_S = \fD_E(\bI_S), \quad 
\widehat{\fD_E(\bI)} = \fD_E(\widehat\bI).$$ 
\item[(3)] 
(Compatibility with $\fR$-saturation) Let $\bI \subset R 
\times \bR$ be an idealistic filtration of r.f.g.~type. 
Then we have 
$$\fR(\bI)_S = \fR(\bI_S), \quad 
\widehat{\fR(\bI)} = \fR(\widehat\bI).$$ 
\end{prop} 
\begin{proof} 
\item[(1)] 
This follows easily from the explicit construction 
of the generation in Lemma \ref{Construction} (1). 
\item[(2)] 
We verify $\fD(\bI)_S = \fD(\bI_S)$.  Firstly we 
show the inclusion $\fD(\bI)_S \subset 
\fD(\bI_S)$.  Note that $\fD(\bI_S) \cap 
\{R\times \bR\}$ is an idealistic filtration over $R$ 
containing $\bI$, and being $\fD$-saturated 
by Lemma \ref{1.1.2.1} (4).  Therefore, we have 
$$\fD(\bI) \subset \fD(\bI_S) \cap \{R \times\bR\} 
\subset \fD(\bI_S).$$ 
At level $a \in \bR$, this implies $\fD(\bI)_a \subset 
\fD(\bI_S)_a$ and hence $(\fD(\bI)_a)_S 
\subset \fD(\bI_S)_a$.  That is to say, we have 
$\fD(\bI)_S \subset \fD(\bI_S)$. 
 
Secondly we show the opposite inclusion $\fD(\bI)_S \supset 
\fD(\bI_S)$.  Note that 
$\fD(\bI)_S$ is an idealistic filtration over $R_S$ containing 
$\bI$, and hence containing 
$\bI_S$.  We claim that 
$\fD(\bI)_S$ is $\fD$-saturated.  In fact, suppose $(f,a) 
\in \fD(\bI)_S$, i.e., $f \in \left\{\fD(\bI)_a\right\}_S$. 
Then, for $d \in \dif_{R_S}^t$, we see by Lemma \ref{1.1.2.1} (7) 
$$d(f) \in \dif_{R_S}^t\left(\left\{\fD(\bI)_a 
\right\}_S\right) = \left\{\dif_R^t\left(\fD(\bI)_a 
\right)\right\}_S \subset\left\{\fD(\bI)_{a-t}\right\}_S.$$ 
That is to say, we have $(d(f),a-t) \in \fD(\bI)_S$, checking 
condition (differential) for $\fD(\bI)_S$.  Thus we 
have $\fD(\bI)_S \supset \fD(\bI_S)$. 
 
This completes the verification for 
$\fD(\bI)_S = \fD(\bI_S)$. 
 
\bigskip
 
The verification for 
$\widehat{\fD(\bI)} = \fD(\widehat\bI)$ 
is identical to the one above using again Lemma \ref{1.1.2.1} 
(7), and left to the reader as an exercise. 
 
\bigskip
 
The verification for the compatibility of localization and completion with 
$\fD_E$-saturation goes almost verbatim to the one above, 
replacing $\fD$ and $\dif_R^t$ with $\fD_E$ and $\dif_{R,E}^t$. 
We leave the verification of the 
statement of Lemma \ref{1.1.2.1} (7) obtained by replacing 
$\dif_R^t$ with $\dif_{R,E}^t$ as an exercise to the reader, 
since it is identical to the one we gave in Chapter 1. 
\item[(3)] 
We use the same notation and argument as in 
Step 1 through Step 4 of the proof of 
Proposition \ref{2.3.2.3} 
(See also Remark \ref{2.4.2.2} (1) below). 
First, since $\bI$, $\bI_S$, and $\widehat\bI$ share the same 
set of generators $T$, 
we may take in Step 1 the common replacement $T'$ at level $L$, 
which keeps the left hand side and right hand side 
of the equation for compatibility intact.  
Therefore, we may assume from the beginning 
that $\bI$ is generated by $T\subset R\times\{L\}$. 
Let $I=\bI_L$ and 
$A = \bigoplus_{n \in \bZ_{\geq0}}I^nX^{Ln} \subset R[X^L]$. 
Note that $A$ describes the idealistic filtration 
$\bI$ (cf.~Definition \ref{2.3.2.0}, Lemma \ref{2.3.2.1}).  
Moreover, 
$$A_S = \bigoplus_{n \in \bZ_{\geq 0}} 
I_S^nX^{Ln} \subset R_S[X^L] 
\quad\text{and}\quad 
\widehat{A}=\bigoplus_{n \in \bZ_{\geq 0}} 
\widehat{I}^nX^{Ln} \subset \widehat{R}[X^L] 
$$ 
describe the localization $\bI_S$ and completion 
$\widehat\bI$, respectively. 
 
Step 2 goes without any change for all $\bI$, $\bI_S$, and 
$\widehat\bI$. 
 
We take $\rho$ in Step 3 so that 
$\rho$ works for all $\bI$, $\bI_S$ and $\widehat\bI$ 
simultaneously.  Set 
$\overline{A}$, $\overline{A_S}$, 
$\overline{\widehat{A}}$ 
as the integral closures 
of $A$ in $R[X^{\frac{1}{\rho}}]$, 
of $A_S$ in $R_S[X^{\frac{1}{\rho}}]$, and 
of $\widehat{A}$ in $\widehat{R}[X^{\frac{1}{\rho}}]$, 
respectively.  Then, in Step 4, we see that 
$\overline{A}$, $\overline{A_S}$, $\overline{\widehat{A}}$ 
describe the idealistic filtrations 
$\fR(\bI)$, $\fR(\bI_S)$, $\fR(\widehat\bI)$, 
respectively. 
 
On the other hand, since $\overline{A}$ describes 
$\fR(\bI)$, it follows by definition that 
$(\overline{A})_S$ and $\widehat{\overline{A}}$ 
describe the localization $\fR(\bI)_S$ 
and completion $\widehat{\fR(\bI)}$, respectively. 
 
Now since the operation of taking the integral closure commutes with 
localization, we have $(\overline{A})_S = \overline{A_S}$.  Thus 
we conclude $\fR(\bI)_S = \fR(\bI_S)$. 
 
As to the question of commutativity of the operation of taking the integral 
closure with completion, recall that $R$ is a finitely generated 
$k$-algebra or its  localization, hence that it is a Grothendieck ring. 
Since $\overline{A}$ is a finitely generated $R$-algebra 
by Step 5 of Proposition \ref{2.3.2.3}, $\overline{A}$ is 
also a Grothendieck ring.  This allows us to conclude that 
$\widehat{\overline{A}}$ is normal, since $\overline{A}$ 
is also normal 
(See Remark 1 to Theorem 32.~6 in \cite{MR879273}). 
Now
$\widehat{\overline{A}}$ is integral over $\widehat{A}$, 
since $\overline{A}$ is integral over $A$. 
Therefore we conclude 
$\widehat{\overline{A}}=\overline{\widehat{A}}$, 
and hence 
$\widehat{\fR(\bI)} = \fR(\widehat\bI)$. 
 
This completes the proof of Proposition \ref{Compatibility}. 
\end{proof} 
\begin{rem}\label{2.4.2.2} 
\item[(1)] 
In \ref{2.3}, the base ring $R$ was assumed to be 
the coordinate ring of an affine open subset of a variety 
$W$ smooth over $k$, or its localization.  
We did {\it not} deal with the case where the 
base ring is the completion $\widehat{R}$.  
Note that the proof of 
Proposition \ref{2.3.2.3} works just as well 
over the base ring being the completion $\widehat{R}$ 
from Step 1 through Step 4, but fails in Step 5, where 
$Q(\widehat{R}[X])$ is not finitely generated over $k$.  
Therefore, we do not claim the stability of 
the idealistic filtrations of r.f.g.~type over 
$\widehat{R}$ under $\fR$-saturation. 
 
Nevertheless, we should emphasize that the 
following assertion is valid: 
\begin{center} 
If an idealistic filtration $\bI$ over $R$ is 
of r.f.g.~type, then so is $\fR(\widehat\bI)$. 
\end{center} 
Indeed, since $\fR(\bI)$ is of r.f.g.~type by 
Proposition \ref{2.3.2.3}, 
the assertion is a direct consequence of compatibility 
$\fR(\widehat\bI)= \widehat{\fR(\bI)}$. 
\item[(2)] 
The assumption of $\bI$ being of r.f.g.~type is indispensable in 
Proposition \ref{Compatibility} (3).  The following gives 
a counterexample to the assertion of compatibility with 
$\fR$-saturation when $\bI$ is not of r.f.g.~type: 
Let $\bI = G(T)$ be an idealistic filtration over $R = k[x,y]$ 
where the set of generators $T$ is an 
infinite set given as below 
$$T = \{(\phi_iy,1-i^{-1}) \mid i \in \bZ_{> 0}\}, 
\ \phi_i = \prod_{j = 1}^i(x - j).$$ 
We claim that, $\fm = (x,y)$ being the maximal ideal 
corresponding to the origin, we have 
$$\fR(\bI_{\fm}) = G_{R_{\fm}}(\{(y,1)\}), \ (y,1) 
\not\in \fR(\bI)_{\fm}.$$ 
This implies that $\fR(\bI_{\fm}) \neq \fR(\bI)_{\fm}$ 
and also that 
$\fR(\widehat\bI) = G_{\widehat{R}}(\{(y,1)\}) 
\neq \widehat{\fR(\bI)} \not\ni (y,1)$ 
where the completion is taken with respect to $\fm$. 
 
Since $\fR(\bI_{\fm}) = G(\{(y,1)\})$ is clear, we only show the 
second part of the claim $(y,1) \not\in 
\fR(\bI)_{\fm}$.  Assume $(y,1) \in \fR(\bI)_{\fm}$. 
Then there exists $f(x,y) \in k[x,y]$ 
such that $f(0,0) \neq 0$ and that $fy \in \fR(\bI)_1$.  Let 
$\bK$ be as in Lemma \ref{Construction} (3). 
Then, for any $l \in \bZ_{> 0}$, we have 
$fy\in\bK_{1-l^{-1}}$ 
and hence $f^{ml}y^{ml} \in \bI_{ml-m}$ 
for some $m \in \bZ_{> 0}$.  Since the generators in $T$ are 
homogeneous with respect to $y$, we see 
that $\bI_{ml-m}$ is a homogeneous ideal with respect to $y$ 
(cf.~Lemma \ref{Construction} (1)).  This implies 
$f(x,0)^{ml}y^{ml} \in \bI_{ml-m}$. 
By Lemma \ref{Construction} (1), we then conclude 
$$f(x,0)^{ml} \in \left(\phi_{i_1}\dotsm\phi_{i_r} 
\mid r \leq ml,\ 
r - \sum_{t=1}^ri_t^{-1}\geq ml - m\right).$$ 
Looking at the range $\{i_1, \dotsc, i_r\}$, we have 
$$1 - r^{-1}\sum_{t=1}^ri_t^{-1} = r^{-1} 
\left(r - \sum_{t=1}^ri_t^{-1}\right) 
\geq (ml)^{-1}(ml-m) = 1 - l^{-1},$$ 
and hence 
$$l^{-1} \geq r^{-1}\sum_{t=1}^ri_t^{-1} 
\geq \left(\max_ti_t\right)^{-1}.$$ 
This implies that each range $\{i_1, \dotsc, i_r\}$ 
contains at least one $i_t$ 
with $i_t \geq l$.  Therefore, we have 
$\phi_l|f(x,0)^{ml}$ and hence $\phi_l|f(x,0)$. 
Since $l \in \bZ_{>0}$ is arbitrary, we conclude 
$f(x,0) = 0$, contradicting the assumption $f(0,0) \neq 0$. 
This contradiction shows $(y,1) \not\in \fR(\bI)_{\fm}$. 
\end{rem} 
We end this section by stating a corollary which says that the results of 
\ref{2.2.3} hold for an idealistic filtration $\bI$ over $R$ 
which is essentially of finite type over $k$, without assuming it has a 
regular system of parameters, if $\bI$ is of r.f.g.~type. 
\begin{cor}\label{2.4.2.3} 
Let $\bI$ be an idealistic filtration 
of r.f.g.~type over $R$ which is essentially of finite type over $k$. 
Then, we have 
$$\fD\fR(\bI)\subset\fR\fD(\bI), 
\quad\fB(\bI)=\fR\fD(\bI).$$ 
Let $E$ be a simple normal crossing divisor. 
Then we have 
$$\fD_E\fR(\bI)\subset\fR\fD_E(\bI), 
\quad\fB_E(\bI) = \fR\fD_E(\bI).$$ 
In particular, 
the operation of taking the $\fB$-saturation or 
$\fB_E$-saturation is compatible with localization or completion 
for an idealistic filtration of r.f.g.~type, 
and the property of being r.f.g.~type is stable under the operation. 
\end{cor} 
\begin{proof} 
Firstly we show the inclusion $\fD\fR(\bI)\subset\fR\fD(\bI)$. 
By Proposition \ref{Compatibility}, 
it suffices to check the inclusion over 
the localization of $R$ at an arbitrary maximal ideal. 
Then, since the localization admits a regular system of parameters, 
we can apply Proposition \ref{2.2.3.1} to verify the inclusion. 
Secondly, in order to prove $\fB(\bI)=\fD\fR(\bI)$, we can repeat 
the argument in Corollary \ref{2.2.3.2}, which is valid 
regardless whether $R$ has a regular 
sytem of parameters or not, once we have the inclusion. 
 
The proof of the logarithmic case is 
almost identical to the one above. 
\end{proof} 
\end{subsection} 
\end{section} 
\end{chapter} 
\begin{chapter}{Leading generator system} 
The purpose of this chapter is to analyze 
the leading terms of the elements of 
an idealistic filtration, i.e., the lowest 
terms of their power series expansions.  Even though our analysis is 
elementary, it leads to the most important notion in the 
Kawanoue program, i.e., that of {\it a leading generator system}.  In this 
chapter, we only give the definition of a leading 
generator system.  However, it could be said that a large portion of our 
entire series of papers, though written with resolution of 
singularities in mind as the principal goal, is a treatise on the 
properties of a leading generator system in its own light. 
 
Our analysis in this chapter is local, or even analytically local. 
Accordingly, we consider an idealistic filtration $\bI$ over $R$ where 
$R$ is taken to be 
the localization at a maximal ideal corresponding to a closed point 
$P \in W$ of the coordinate ring of an affine open subset of a 
variety 
$W$ smooth over an algebraically closed field $k$ of 
$\on{char}(k) = p \geq 0$, or its completion.  We denote by $\fm$ 
the maximal ideal of $R$. 
 
It is worth emphasizing that the main results of this chapter are obtained 
assuming that the idealistic filtration is 
$\fD$-saturated. 
 
\bigskip
 
The main object of our study is the graded $k$-subalgebra 
$L(\bI) = \bigoplus_{n \in \bZ_{\geq 0}}L(\bI)_n$, 
formed by the leading terms of the 
idealistic filtration 
(cf.~\ref{3.1}), 
of the graded ring $G = 
\bigoplus_{\geq 0}\fm^n/\fm^{n+1} = 
\bigoplus_{n \in \bZ_{\geq 0}}G_n$, 
which is isomorphic to a polynomial ring with 
$d (=\dim W)$-variables over 
$k$. 
 
In characteristic zero, if $\bI$ is $\fD$-saturated, $L(\bI)$ is 
generated as a graded algebra over $k$ by its degree one component 
$L(\bI)_1$, i.e., $L(\bI) = k[L(\bI)_1]$.  
Moreover, the hypersurfaces of maximal contact 
correspond exactly to the elements 
of the form $(h,1) \in \bI$ whose leading terms belong to 
$L(\bI)_1 \setminus \{0\}$, i.e., $\overline{h} 
= (h\bmod{\fm}) \in L(\bI)_1 \setminus \{0\}$.  
A fundamental observation of the Kawanoue program is then that the 
invariants we need to build a sequence of blowups for resolution of 
singularities can be 
all constructed from a collection $\{(h_i,1)\} \subset \bI$ with 
$\{\overline{h_i}\}$ forming 
a basis of $L(\bI)_1$ and hence generating the graded algebra 
$L(\bI)$, instead of taking a hypersurface of maximal 
contact one by one. 
 
In characteristic $\on{char}(k) = p > 0$, in contrast, $L(\bI)$ may 
not 
be generated as a graded algebra over $k$ by its degree one component 
$L(\bI)_1$ even if $\bI$ is $ 
\fD$-saturated.  Or worse, $L(\bI)_1$ may be $0$, i.e., there is no 
hypersurface of maximal contact.  However, if 
$\bI$ is $\fD$-saturated, $L(\bI)$ is generated as a graded 
algebra over $k$ by $\bigcup_{e \in 
\bZ_{\geq 0}}L(\bI)_{p^e}^{\on{pure}}$, i.e., 
$$L(\bI) = k[\bigcup_{e \in 
\bZ_{\geq 0}}L(\bI)_{p^e}^{\on{pure}}],$$ 
where $L(\bI)_{p^e}^{\on{pure}} 
\subset L(\bI)_{p^e}$ is the subspace consisting of ``pure'' elements. 
(We call an element $w \in L(\bI)_{p^e}$ ``pure'' 
if $w = v^{p^e}$ for some $v \in G_1$.) 
Observing that there is a sequence of inclusions 
$$L(\bI)_1 = L(\bI)_{p^0}^{\on{pure}}, \quad 
\{L(\bI)_{p^0}^{\on{pure}}\}^p \subset L( 
\bI)_{p^1}^{\on{pure}}, \quad \{L(\bI)_{p^1}^{\on{pure}}\}^p 
\subset L(\bI)_{p^2}^{\on{pure}} \dotsb ,$$ 
which stabilizes to a sequence of equalities after some point, i.e., there 
exists \linebreak $e_N \in \bZ_{\geq 0}$ 
such that for $e > e_N$ the above inclusions become equalities 
$$\{L(\bI)_{p^{e - 1}}^{\on{pure}}\}^p = L( 
\bI)_{p^e}^{\on{pure}},$$ 
we are led to the following notion of {\it a leading generator system}. 
 
\bigskip
 
We call a subset ${\mathbb H} = \{(h_{ij},p^{e_i})\} \subset \bI$ a 
leading generator system if it satisfies the following 
conditions: 
 
(i) $h_{ij} \in \fm^{p^{e_i}}$ and $\overline{h_{ij}} = 
(h_{ij}\bmod{\fm^{p^{e_i} + 1}}) \in L( 
\bI)_{p^{e_i}}^{\on{pure}}$, 
 
(ii) $\{\overline{h_{ij}}^{p^{e-e_i}} 
\mid e_i \leq e\}$ 
consists of $\# \{(i,j)\mid e_i \leq 
e\}$-distinct elements, and forms a basis of 
$L(\bI)_{p^e}^{\on{pure}}$ for any 
$e \in \bZ_{\geq 0}$. 
 
Therefore, if $\bI$ is $\fD$-saturated, the leading terms of the 
elements in the leading generator system generates 
$L(\bI)$ as a graded algebra over $k$, i.e., 
$$L(\bI) = k[\{\overline{h_{ij}}\}].$$ 
(Note that we have $\dim_kL( 
\bI)_{p^e}^{\on{pure}} \leq \dim W$ for any $e \in 
\bZ_{\geq 0}$ and 
hence that condition (ii) implies 
$\# {\mathbb H} \leq \dim W$.) 
 
\bigskip
 
The Kawanoue program in its simplest terms is a program to construct an 
algorithm for resolution of singularities using a leading 
generator system as a collective substitute for a hypersurface of maximal 
contact, which in the existing algorithms in 
characteristic zero is the key for the inductive structure. 
\begin{section}{Analysis of the leading terms 
of an idealistic filtration.}\label{3.1} 
\begin{subsection}{Definitions}\label{3.1.1} 
\begin{defn}\label{3.1.1.1} 
\item[(1)] 
Let $\bI$ be an idealistic filtration over $R$ with its maximal 
ideal $\fm$.  Recall that the maximal ideal $\fm$ 
corresponds to a closed point $P \in W$. 
Set 
$$G = \bigoplus_{n \in \bZ_{\geq 0}}\fm^n/\fm^{n+1} = 
\bigoplus_{n \in \bZ_{\geq 0}}G_n.$$ 
We define the graded $k$-subalgebra 
$$L(\bI) = \bigoplus_{n \in \bZ_{\geq 0}}L(\bI)_n \subset 
G,$$ 
which we call the leading algebra of the idealistic filtration $\bI$ at 
$P$, by setting 
$$L(\bI)_n = \{\overline{f} = (f\bmod{\fm^{n+1}} 
)\mid (f,n) \in \bI, f \in 
\fm^n\}.$$ Note that $L(\bI)_0 = k$ by condition (o) in Definition 
\ref{2.1.1.1} (2). 
\item[(2)] 
Set $p = \on{char}(k)$ when $k$ is of positive characteristic, or 
$p= \infty$ when $k$ is of characteristic zero.  For $e 
\in \bZ_{\geq 0}$ with $p^e \in \bZ_{> 0}$, we define the pure 
part $L(\bI)_{p^e}^{\on{pure}}$ of $L( 
\bI)_{p^e}$ by the formula 
$$L(\bI)_{p^e}^{\on{pure}} = L(\bI)_{p^e} \cap F^e(G_1) \subset 
L(\bI)_{p^e},$$ 
where $F^e$ is the $e$-th power of the Frobenius map of $G$ 
(cf.~Definition \ref{1.3.1.1}). 
 
An element $w \in L(\bI)_{p^e}$ is called pure if $w \in L( 
\bI)_{p^e}^{\on{pure}}$. 
\end{defn} 
\begin{rem}\label{3.1.1.2} 
\item 
If we choose a regular system of parameters $(x_1, \dotsc, x_d)$ for $R$, 
there is a natural isomorphism $G \cong k[x_1, \dotsc, 
x_d]$.  Through this isomorphism, we may identify $G$ with the polynomial 
ring over $k$. 
\item 
We use the convention that $\infty^n = \infty$ for 
$n \in \bZ_{>0}$ and $\infty^0 = 1$ (cf.~\ref{0.2.3.2.1}). 
Therefore, the only pure part we consider in characteristic 
zero is in degree one, where we have 
$$L(\bI)_{\infty^0}^{\on{pure}} = L(\bI)_1^{\on{pure}} = 
L(\bI)_1.$$ 
In other words, in charactersitic zero, all the pure elements are in degree 
one. 
\item 
We see that $L(\bI)_n$ is a $k$-vector subspace of 
$G_n$, which follows from the definition of an idealistic filtration 
$\bI$.  Using the assumption that $k$ is algebraically 
closed, we also see that $L(\bI)_{p^e}^{\on{pure}}$ is a $k$-vector 
subspace of $L(\bI)_{p^e}$. 
\end{rem} 
\end{subsection} 
\begin{subsection}{Heart of our analysis.}\label{3.1.2} 
The following lemma 
sits at the heart of our analysis, though its 
statement and proof are quite elementary. 
\begin{lem}\label{3.1.2.1} 
Let $G = k[x_1, \dotsc, x_d]$ be the polynomial ring 
over $k$ with $d$ variables $X = (x_1, \dotsc, x_d)$. 
We regard 
$G$ as a graded 
$k$-algebra with natural grading defined by the degree. 
We define ``$p$'' as in Definition \ref{3.1.1.1} 
and we use the same convention as in Remark \ref{3.1.1.2}. 
 
Let $L = \bigoplus_{n \in \bZ_{\geq 0}}L_n \subset G$ 
be a graded $k$-subalgebra of $G$ with $L_0 = G_0 = k$. 
Suppose that $L$ is $\fD$-saturated in the sense that it satisfies the 
following condition: 
$$f \in L, \partial \in \dif_G \Longrightarrow \partial(f) \in L.$$ 
Then $L$ is generated as a graded algebra over $k$ by its pure part 
$L^{\on{pure}} := \bigsqcup_{e \in \bZ_{\geq 
0}}L_{p^e}^{\on{pure}}$ where $L_{p^e}^{\on{pure}} = L_{p^e} \cap 
F^e(G_1) \subset L_{p^e}$, i.e., $L = k[L^{\on{pure}}]$. 
 
Moreover, we can choose $\{e_1 < \dotsb < e_N\} \subset \bZ_{\geq 0}$ 
and $V_1 \sqcup \dotsb \sqcup V_N 
\subset G$ with $V_i = \{v_{ij}\}_j$ satisfying 
the following conditions: 
\begin{enumerate} 
\renewcommand{\labelenumi}{(\roman{enumi})} 
\item 
$F^{e_i}(V_i) \subset L_{p^{e_i}}^{\on{pure}}$ for $1 \leq i \leq N$, 
\item 
$\bigsqcup_{e_i \leq e}F^e(V_i)$ is a $k$-basis of 
$L_{p^e}^{\on{pure}}$ for any $e \in \bZ_{\geq 0}$. 
\end{enumerate} 
In particular, we have $L = k[\bigsqcup_{i = 1}^NF^{e_i}(V_i)]$ with 
$\sum_{i = 1}^N \# V_i 
\leq d$. 
\end{lem} 
\begin{proof} 
We prove the following assertion 
$$(\heartsuit)_d \quad L = k[L^{\on{pure}}]$$ 
by induction on the number of variables $d$.  When $d = 0$, we have 
$G = L = k$ and $L^{\on{pure}} = \emptyset$.  Thus we obviously have 
$(\heartsuit)_0$. 
 
Now we prove $(\heartsuit)_d$ assuming $(\heartsuit)_{d-1}$.  Take $f \in 
L$.  It suffices to show $f \in 
k[L^{\on{pure}}]$.  We may assume that 
$f$ is homogeneous of degree 
$r 
\in \bZ_{\geq 0}$, i.e., $f \in L_r$.  Set 
$$s = \max\{t \in \bZ_{\geq 0} \mid f \in F^t(G)\}, 
\quad r = r'p^s,$$ 
and take $g \in G_{r'}$ such that $f = g^{p^s}$.  We write $g = \sum_{|J| = 
r'}c_JX^J \in G_{r'}$ with $c_J \in 
k$. 
 
By the maximality of $s$ (and since $k$ is algebraically closed), we observe 
that there exists $J_o$ with $|J_o| = r'$ such that 
$c_{J_o} \neq 0$ and that 
$p\not| J_o = (j_{o1}, \dotsc , j_{od})$, 
i.e., $p\not|\ j_{o\alpha}$ for some $\alpha$.  By renumbering 
the variables, we may assume $p \not|\ j_{od}$. 
 
We compute 
$$z = \partial_{J_o - {\mathbf{e}_d}}g = j_{od}c_{J_o} \cdot x_d + 
\sum_{\alpha = 1}^{d - 
1}(j_{o\alpha} + 1)c_{J_o - {\mathbf{e}}_d + {\mathbf{e}}_{\alpha}} \cdot 
x_{\alpha}.$$ 
Since $j_{od}c_{J_o} \in k \setminus \{0\}$, we may 
take $(x_1, \dotsc, x_{d-1},z)$ as a new system of variables 
for the polynomial ring $G$.  We set $G' = k[x_1, \dotsc, x_{d-1}]$ 
to be the polynomial ring with $(d-1)$-variables and $L' = L \cap G'$.  Note 
that $L'$ is $\fD$-saturated.  Rewrite 
$g = 
\sum_{i = 0}^{r'}a_iz^i$ with 
$a_i \in G'_{r' - i}$ in terms of the new system of variables. 
 
Note that, for any $h \in G$ and $K \in \bZ_{\geq 0}^d$, we have 
$$\partial_{X^{p^sK}}(h^{p^s}) = \partial_{X^K}(h)^{p^s},$$ 
which follows from the equations (cf.~Remark \ref{1.2.1.3} (2)) 
$$\partial_{p^sK}(X^{p^sJ}) 
= \binom{p^sJ}{p^sK}X^{p^s(J-K)} 
\quad\text{and}\quad\binom{p^sJ}{p^sK}=
\binom{J}K=\binom{J}K^{p^s}$$ 
 
Thus we have 
$$\left\{\begin{array}{rl} 
\partial_{X^{p^sK}}(f) 
&=\partial_{X^{p^sK}}(g^{p^s}) 
 =\partial_{X^K}(g)^{p^s} 
 =z^{p^s}, \quad 
K=J_o - {\mathbf{e}_d},\\ 
  \partial_{z^{p^si}}(f) 
&=\partial_{z^{p^si}}(g^{p^s}) 
 =\partial_{z^i}(g)^{p^s} 
 = a_i^{p^s} + \displaystyle{\sum_{t = i+1}^{r'}
\binom{t}ia_t^{p^s}z^{(t-i)p^s}}\ (0 \leq i \leq r'). 
\end{array}\right.$$ 
Recall that $L$ is $\fD$-saturated.  Thus the first formula implies 
$z^{p^s} \in L$, and the second formula 
implies $a_i^{p^s} \in L$ for $0 \leq i \leq r'$ by descending induction on 
$i$. 
 
On one hand, we have $z^{p^s} \in L^{\on{pure}}$ by definition.  On the 
other hand, since $L' = k[L'^{\on{pure}}]$ 
by inductional hypothesis on the number of variables, we have 
$$\{a_i^{p^s};0 \leq i \leq r'\} \subset L \cap G' = L' = k[L'^{\on{pure}}] 
\subset k[L^{\on{pure}}].$$ 
Therefore, we conclude 
$$f = g^{p^s} = \sum_{i = 0}^{r'}a_i^{p^s}z^{p^si} \in k[L^{\on{pure}}].$$ 
This completes the inductional step and hence the proof for $L = 
k[L^{\on{pure}}]$. 
 
In order to see the ``Moreover'' part of the statement, observe that there 
is 
a sequence of inclusions among the pure parts 
$$L_{p^0}^{\on{pure}} = L_1, \quad \{L_{p^0}^{\on{pure}}\}^p 
\subset L_{p^1}^{\on{pure}}, 
\quad \{L_{p^1}^{\on{pure}}\}^p \subset L_{p^2}^{\on{pure}}, 
\quad \dotsb.$$ 
Let $e_1 < \dotsb < e_i < \dotsb < e_N$ be the integers indicating the 
stages where the above inclusions are not 
equalities, i.e., 
$$\left\{\begin{array}{ll} 
\{L_{p^{e_i - 1}}^{\on{pure}}\}^p &\subsetneqq L_{p^{e_i}}^{\on{pure}} 
\quad(1 \leq i \leq N) \\ 
\{L_{p^{e - 1}}^{\on{pure}}\}^p &= L_{p^e}^{\on{pure}} \quad e \not\in 
\{e_i\mid 1 \leq i \leq 
N\}. 
\end{array}\right.$$ 
Note that the set of such integers is a finite set, since the dimension of 
the pure part is uniformly bounded, 
i.e., $\dim_kL_{p^e}^{\on{pure}}\leq 
\dim G_1 = d$ for any $e \in \bZ_{\geq 0}$. 
 
Now we have only to take $V_i \subset G_1 \quad (i = 1, \dotsc, N)$ 
inductively so that 
$$F^{e_i}(V_i) \cup \bigcup_{j < i}F^{e_i}(V_j) \subset 
L_{p^{e_i}}^{\on{pure}}$$ 
forms a basis of $L_{p^{e_i}}^{\on{pure}}$ for $1 \leq i \leq N$. 
 
This completes the proof of Lemma \ref{3.1.2.1}. 
\end{proof} 
\end{subsection} 
\begin{subsection}{Leading generator system.}\label{3.1.3} 
The statement of Lemma \ref{3.1.2.1} 
leads us to the following notion of a leading 
generator system of a $\fD$-saturated idealistic filtration. 
\begin{defn}\label{3.1.3.1} 
Let $\bI$ be a $\fD$-saturated 
idealistic filtration. 
We call a subset ${\mathbb H} = \{(h_{ij},p^{e_i})\} \subset \bI$ 
a leading generator system if it satisfies the following 
conditions: 
\begin{enumerate} 
\renewcommand{\labelenumi}{(\roman{enumi})} 
\item 
$h_{ij} \in \fm^{p^{e_i}}$ and $\overline{h_{ij}} = 
(h_{ij}\bmod{\fm^{p^{e_i} + 1}}) \in L( 
\bI)_{p^{e_i}}^{\on{pure}}$, 
\item 
$\{\overline{h_{ij}}^{p^{e-e_i}}\mid e_i \leq e\}$ 
consists of $\# \{ij\mid e_i \leq 
e\}$-distinct elements, and forms a basis of 
$L(\bI)_{p^e}^{\on{pure}}$ for any 
$e\in \bZ_{\geq 0}$. 
\end{enumerate} 
\end{defn} 
\begin{prop}\label{3.1.3.2} 
A leading generator system exists for a 
$\fD$-saturated idealistic filtration $\bI$. 
\end{prop} 
\begin{proof} 
Since $\bI$ is $\fD$-saturated, it follows 
that $L(\bI)$ is $\fD$-saturated and hence 
that we can apply Lemma \ref{3.1.2.1} to $L = L(\bI) \subset G$.  
Take a collection ${\mathbb H} = \{(h_{ij},p^{e_i})\}$ 
so that $\overline{h_{ij}} = v_{ij}^{p^{e_i}}$, where 
$\{e_1 < \dotsb < e_N\}$ and 
$V_1 \sqcup \dotsb \sqcup V_N \subset G_1$ with 
$V_i = \{v_{ij}\}_{j}$ are taken as stated in Lemma \ref{3.1.2.1}, 
satisfying conditions (i) and (ii).  Then 
${\mathbb H}$ is a leading generator system for $\bI$. 
\end{proof} 
\begin{rem}\label{3.1.3.3} 
\item[(1)] 
Condition (i) in Definition \ref{3.1.3.1} can be 
rephrased in terms of the 
differential operators as follows: 
\begin{enumerate} 
\renewcommand{\labelenumi}{(\roman{enumi})} 
\item
$D(h_{ij}) \in \fm 
$ for any $
D \in \dif_R^{(p^{e_i})}$, 
\end{enumerate} 
where $\dif_R^{(p^{e_i})}$ 
is defined by the following formula 
$$\dif_R^{(p^{e_i})} := \sum_{a + b = p^{e_i}, 0 < a, b < 
p^{e_i}}\dif_R^a \circ \dif_R^b.$$ 
\item[(2)] 
We often study a subset 
${\mathbb H} = \{(h_{ij},p^{e_i})\} \subset \bI$ 
which satisfies some slightly weaker conditions than 
those for a leading generator system.  Namely, we 
require condition (i), and instead of 
full condition (ii) where 
$\{\overline{h_{ij}}^{p^{e-e_i}}\mid e_i \leq e\}$ should 
form a basis of $L(\bI)_{p^e}^{\on{pure}}$, we only require 
$\{\overline{h_{ij}}^{p^{e-e_i}}\mid e_i \leq e\}$ to be 
$k$-linearly independent. 
 
The class of the subsets described above, which is slightly bigger than the 
class of the leading generator systems, is 
often better suited for the purpose of setting up some inductional proofs. 
We refer the reader to Chapter 4 for the examples of 
such proofs. 
\end{rem} 
\end{subsection} 
\end{section} 
\begin{section}{Invariants $\sigma$ 
and $\widetilde{\mu}$.}\label{3.2} 
In this section, we present the definitions of the two of the most basic 
invariants (at the closed point $P \in W$) that we use in 
our algorithm, $\sigma$ and $\widetilde{\mu}$, in relation to the notion of 
a leading generator system. 
 
We warn the reader, however, that in the actual process of our algorithm, 
the definitions 
of $\sigma$ and $\widetilde{\mu}$ will be slightly modified.  For example, 
in order to determine the 
invariant $\widetilde{\mu}$ in our setting, we also have to take the 
boundary divisor of reference into consideration, just as we 
do to determine the weak order in the classical setting. 
 
The purpose of this presentation is to bring a flavor of how these 
invariants may 
function in our algorithm, while the details, including their fundamental 
properties, will be discussed in Parts II, III, and IV. 
\begin{subsection}{Invariant $\sigma$.}\label{3.2.1} 
\begin{defn}\label{3.2.1.1} 
Let $\bI$ be a $\fD$-saturated idealistic filtration. 
Then the invariant $\sigma$ is defined to be the infinite sequence indexed 
by $e \in \bZ_{\geq 0}$, i.e., 
$$\sigma = (d - l_{p^0}^{\on{pure}}, d - l_{p^1}^{\on{pure}}, \dotsc , d 
- l_{p^e}^{\on{pure}}, \dotsc )$$ 
where 
$$ 
d = \dim W, \quad 
l_{p^e}^{\on{pure}} = \dim_kL(\bI)_{p^e}^{\on{pure}}. 
$$ 
More precisely, $\sigma$ should be considered as a function $\sigma: 
\bZ_{\geq 0} \rightarrow \bZ_{\geq 
0}$ defined by 
$$\sigma(e) = d - l_{p^e}^{\on{pure}}.$$ 
\end{defn} 
\begin{rem}\label{3.2.1.2} 
\item[(1)] 
The reason why we take the (infinite) sequence of $d - 
l_{p^e}^{\on{pure}}$ instead of the (infinite) 
sequence of $l_{p^e}^{\on{pure}}$ is two-fold: 
\begin{enumerate} 
\item[(i)] 
When we consider the invariant $l_{p^e}^{\on{pure}}$, fixing 
$e$ but varying $P \in W$, we see (cf.~Part II) that it is lower 
semi-continuous.  Taking its negative, we have our 
invariant upper semi-continuous as desired. 
 
(In other words, the bigger $l_{p^e}^{\on{pure}}$ is, the better 
the singularities are.  Therefore, as the measure of how bad the 
singularities are, it is natural to take our 
invariant to be its negative.) 
\item[(ii)] 
We reduce the problem of resolution of singularities of a variety $X$ 
to that of embedded resolution.  Therefore, it would be desirable 
or even necessary to come up with an algorithm which would induce 
the ``same'' process of resolution of singularities, no matter what 
ambient variety $W$ we choose for an embedding 
$X\hookrightarrow W$ (locally). 
 
While $l_{p^e}^{\on{pure}}$ is dependent of the choice of $W$, $\dim W - 
l_{p^e}^{\on{pure}}$ is not.  Therefore, the latter 
is more appropriate as an invariant toward constructing such an algorithm. 
\end{enumerate} 
\item[(2)] 
As observed before, the dimension of the pure part is non-decreasing 
(with respect to the power $e 
\in \bZ_{\geq 0}$ of $p^e$), and is  uniformly bounded from above by 
$d =\dim W$, i.e., 
$$0 \leq l_{p^0}^{\on{pure}} \leq l_{p^1}^{\on{pure}} 
\leq\dotsb l_{p^{e - 1}}^{\on{pure}} \leq 
l_{p^e}^{\on{pure}} \leq \dotsb\leq d = \dim W$$ 
and hence stabilizes after some point, i.e., there exists 
$e_N \in \bZ_{\geq 0}$ such that for $e > e_N$ 
the above inequalities become equalities 
$$l_{p^{e - 1}}^{\on{pure}} = l_{p^e}^{\on{pure}}.$$ 
Therefore, though $\sigma$ is an infinite sequence by definition, 
essentially we are only looking at some finite 
part of it. 
\item[(3)] 
In characteristic zero, the invariant $\sigma$ consists of 
only one term $(d - \dim L(\bI)_1)$, where $\dim L(\bI)_1$ 
can be regarded as indicating ``how many linearly independent 
hypersurfaces of maximal contact we can take'' for $\bI$. 
\end{rem} 
\end{subsection} 
\begin{subsection}{Invariant $\widetilde{\mu}$.}\label{3.2.2} 
\begin{defn}\label{3.2.2.1} 
Let $\bI$ be a $\fD$-saturated idealistic filtration. 
Take a leading generator system 
${\mathbb H} = \{(h_{ij},p^{e_i})\}$ of $\bI$. 
Set $\cH = \{h_{ij}\}$ and $(\cH) 
\subset R$ to be the ideal generated by $\cH$. 
 
For $f \in R$, define its multiplicity (or order) modulo 
$(\cH)$ to be 
$$\on{ord}_\cH(f) = \sup\left\{n \in \bZ_{\geq 
0}\mid f \in \fm^n + (\cH)\right\}$$ 
and 
$$\mu_\cH(\bI) := \inf\left\{\mu_\cH(f,a) 
:= \frac{\on{ord}_\cH(f)}{a}\mid (f,a) \in 
\bI, a > 0\right\}.$$ 
(Note that we set $\on{ord}_\cH(0) = \infty$ 
by definition.) 
We define the invariant $\widetilde{\mu}$ by 
$$\widetilde{\mu} = \mu_\cH(\bI).$$ 
\end{defn} 
\begin{rem}\label{3.2.2.2} 
\item[(1)] 
We will see in Part II that $\mu_\cH(\bI)$ is 
independent of the choice of a leading generator system, 
and hence that the invariant $\widetilde{\mu}$ is actually 
an intrinsic one associated to the idealistic filtration $\bI$. 
\item[(2)] 
In characteristic zero, where $\cH$ forms (a part of) a regular 
system of parameters, the upper semi-continuity of the 
invariant $\widetilde{\mu}$ (along the locus where the invariant 
$\sigma$ is constant) follows immediately from the upper semi-continuity of 
the usual multiplicity defined on the nonsingular 
subvariety defined by the ideal $(\cH)$.  In positive 
characteristic, 
however, it is highly non-trivial, and its verification 
is one of the main subjects of Part II. 
\item[(3)] 
In our algorithm, the invariant $\widetilde{\mu}$ is actually computed 
as $\mu_{\cH,E}(\bI)$, using not only the 
information about a leading generator system but also the one about the 
boundary divisor $E$ in reference.  For all the details, 
we refer the reader to Parts II, III, and IV. 
\item[(4)] 
In Part II, we study the power series expansion of $f \in R$ with 
respect to (the elements in $\cH$ associated to) a 
leading generator system.  There the invariant 
$\on{ord}_\cH(f)$ can be computed as the 
multiplicity of the ``constant'' term. 
Again we refer the reader to Part II for its detail. 
\item[(5)] 
In characteristic zero, the invariant $\widetilde{\mu}$ 
corresponds to the multiplicity of what is called 
the coefficient ideal (restricted to a hypersurface 
of maximal contact) in the classical setting. 
\end{rem} 
\end{subsection} 
\end{section} 
\end{chapter} 
\begin{chapter}{Nonsingularity principle.} 
The purpose of this chapter is to establish the nonsingularity principle, 
which guarantees the 
nonsingularity of the center of each blowup in our algorithm 
\linebreak 
(cf.~\ref{0.2.3.2.4} in the introduction). 
 
In \ref{4.1}, we prepare some technical lemmas that 
we use in the proof of the 
nonsingularity principle.  They describe the behavior 
of a leading generator system, which we expect to be parallel to the 
behavior of a collection of hypersurfaces of maximal contact 
forming (a part of) a regular system of parameters.  We will use these 
lemmas again later in our series of papers. 
 
\ref{4.2}, where we present the statement and proof of the nonsingularity 
principle, is literally the culminating 
point of Part I. 
 
In this chapter, $R$ represents the localization at a maximal ideal, 
or its completion, of the coordinate ring of an affine open subset 
of a variety $W$ smooth over 
an algebraically closed field $k$ of 
$\on{char}(k) = p > 0$, or characteristic zero where we 
formally set $p =\infty$ in the arguments below. 
We denote by $\fm$ the maximal ideal of $R$, 
which corresponds to a closed point $P \in W$. 
\begin{section}{Preparation toward the 
nonsingularity principle.}\label{4.1} 
\begin{subsection}{Setting for 
the supporting lemmas.}\label{4.1.1} 
We fix the following setting for the 
three supporting lemmas we present in \ref{4.1.2}: 
 
Let $\cH = \{h_1, \ldots , h_N\} \subset 
R$ be a subset of $R$ consisting of $N$ elements, and let 
$0 \leq e_1 \leq\cdots\leq e_N$ 
be nonnegative integers associated to these 
elements, satisfying the following conditions: 
\begin{enumerate} 
\renewcommand{\labelenumi}{(\roman{enumi})} 
\item 
$h_l \in \fm^{p^{e_l}}$ and $\overline{h_l} = (h_l 
\bmod{\fm^{p^{e_l} + 1}}) \in \gpure{e_l}$ for $l 
= 1, \ldots , N$. 
(See Definition \ref{3.1.1.1}.) 
\item 
$\{\overline{h_l}^{p^{e_s - e_l}}\mid e_l 
\leq e_s\}$ consists of $\# \{l\mid e_l \leq 
e_s\}$-distinct and 
$k$-linearly independent elements 
in the $k$-vector space 
$\gpure{e_s}$ for $s = 1, \ldots, N$. 
\end{enumerate} 
We set 
$$\left\{\begin{array}{ll} 
e &:= e_1 = \min\{e_l\mid l = 1, \ldots , N\}, \\ 
L &:= \max\{l\mid l = 1, \ldots, N, e_l = e\} = \# 
\{l\mid l = 1, \ldots , N, e_l = e\}. \\ 
e' &:= e_{L + 1} \quad (\text{if\ } 
L = N, \text{\ then\ we\ set\ }e' =\infty).\\ 
\end{array}\right.$$ 
Let $(x_1, \ldots , x_d)$ be a regular system of parameters 
for $R$ such that 
$$M_{d,L} = \left[\partial_{x_i^{p^e}} 
(h_l)\right]_{l = 1, \ldots , L}^{i = 1,\ldots , d} 
\in M(d\times L,R)$$ has the invertible $L \times L$ first 
minor, i.e., 
$$M = \left[\partial_{x_i^{p^e}} 
(h_l)\right]_{l = 1, \ldots , L}^{i = 1, \ldots,L} 
\in \on{GL}(L, R).$$ 
 
Let 
$C = \left[c_{ij}\right]_{j =1,\ldots,L}^{i = 1,\ldots,L} 
\in\on{GL}(L,R)$ be the inverse matrix of 
$M$ so that $CM = I_L$. 
 
We introduce the following multi-index notations: 
$$ 
T := (t_1, \ldots , t_L) \in \bZ_{\geq 0}^L, 
\quad 
c^T := \prod_{j = 1}^L\left(c_{L,j}\right)^{t_j}. 
$$ 
\begin{rem}\label{4.1.1.1} 
\item[(1)] 
Condition (i) in Setting \ref{4.1.1} 
can be rephrased in terms of the 
differential operators as follows 
(cf.~Remark \ref{3.1.3.3} (1)): 
\begin{enumerate} 
\item[(i)] 
$D(h_l) \in \fm \quad \forall D \in 
\on{Diff}_R^{(p^{e_l})}$. 
\end{enumerate} 
\item[(2)] 
We are {\it not} assuming in our situation that the subset 
$\cH$ is associated to a leading 
generator system of an idealistic filtration. 
See Remark \ref{3.1.3.3} (2). 
\item[(3)] 
Conditions (i) and (ii) imply that, for any regular system 
of parameters 
\linebreak 
$(y_1, \ldots, y_d)$, the matrix 
$M_y = \left[\partial_{y_i^{p^e}} 
(h_l)\right]_{l = 1, \ldots , L}^{i = 1, \ldots ,d}$ 
is of size $d\times L$ and has the full rank, i.e., 
$\on{rank}\ M_y = L$.  Therefore, 
by a linear change of variables, we may 
always come up with a regular system of parameters 
$(x_1, \ldots , x_d)$ satisfying the 
condition in our situation. 
\end{rem} 
\end{subsection} 
\begin{subsection}{Statements and proofs of 
the supporting lemmas.}\label{4.1.2} 
Given a regular system of parameters $(x_1, \ldots, x_d)$, we have the 
corresponding partial differential operators $\partial_{x_i^u} 
\quad (u\in \bZ_{> 0})$ for $1 \leq i 
\leq d$.  Given a set of elements 
$(h_1, \ldots, h_N)$ as described in the setting 
(e.g.~the set associated to a leading generator system), 
we would like to have their corresponding partial 
differential operators.  The next lemma constructs a 
differential operator 
$D_u$, which behaves like ``$\partial_{h_L^u}$'' in spirit when we look at 
the initial terms of our concern. 
\begin{lem}[Supporting Lemma 1]\label{4.1.2.1} 
Let $u, r$ be integers such that 
$$0 \leq u < p^{e' - e} \quad \text{\ and\ }\quad 0 \leq r.$$ 
Define 
$$D_u := \sum_{|T| = u}c^T\partial_{p^eT} \in \on{Diff}^{up^e}_R \text{\ 
and\ } D_{-1} = 0,$$ 
where we use the abbreviated notation $\partial_{J} = \partial_{X^J}$. 
 
Then for any $\beta \in \fm^r$ and $1 \leq l \leq N$, we have 
$$D_u\left(\beta h_l\right) 
\equiv (D_u \beta) h_l + \delta_{L,l}D_{u - 
1}\beta \mod{\fm^{r + p^{e_l} - up^e + 1}}.$$ 
\end{lem} 
\begin{proof} 
By the generalized product rule, we have 
\begin{eqnarray*} 
D_u\left(\beta h_l\right) 
&=& \sum_{|T| = u}c^T\partial_{p^eT}\left(\beta 
h_l\right) 
= \sum_{|T| = u}c^T \sum_{J \leq p^eT}\left(\partial_{p^eT - 
J}\beta\right)\left(\partial_Jh_l\right). 
\end{eqnarray*} 
We remark here that 
$$\left\{p^{e_l} \not| J\right\} 
\text{\ or\ }\{p^{e_l} < |J|\} 
\Longrightarrow 
\text{ord}_P(\partial_Jh_l) > p^{e_l} - |J|.$$ 
Thus, in the process of continuing the above computation 
$\bmod{\fm^{r + p^{e_l}- up^e + 1}}$, the term 
$\partial_Jh_l$ is relevant only when 
$J = p^{e_l}{\mathbf e}_j \quad(1\leq j \leq L)$ 
or when $J = {\mathbf O}$.  Therefore, we have 
$$ 
D_u\left(\beta h_l\right) \equiv 
\left(\sum_{|T| = u}c^T\partial_{p^eT}\beta\right) 
h_l + \sum_{j = 1}^L \sum_{|T| = u}c^T 
\left(\partial_{p^eT - p^{e_l}{\mathbf e}_j} 
\beta\right)\left(\partial_{p^{e_l}{\mathbf e}_j}h_l\right) 
$$ 
where the first and the second term in the right hand side 
correspond to the case  $J = {\mathbf O}$ and 
$J = p^{e_l}{\mathbf e}_j$ for $1 \leq j \leq L$ 
respectively. 
 
We remark here that in the generalized product rule we only consider the 
case where ${\mathbf O} \leq p^eT - p^{e_l}{\mathbf e}_j$. 
Looking at the $j$-th components, we conclude 
$$0 \leq p^et_j - p^{e_l} \leq p^e|T| - p^{e_l} < 
p^e \cdot p^{e'-e} - p^{e_l} = p^{e'} - p^{e_l}.$$ 
This implies that we only consider the case where 
$$e_l = e \text{\ and\ } t_j \geq 1.$$ 
Therefore, setting $T' = T - {\mathbf e}_j$, we have 
\begin{eqnarray*} 
D_u\left(\beta h_l\right) 
&\equiv& \left(D_u\beta\right) h_l + 
\delta_{e,e_l}\sum_{j = 1}^L \sum_{|T'| = u - 1}c^{T' + {\mathbf 
e}_j}\left(\partial_{p^eT'}\beta\right)\left(\partial_{p^e{\mathbf 
e}_j}h_l\right)\\ 
&=& \left(D_u\beta\right) h_l + \delta_{e,e_l}\left(\sum_{j = 1}^L 
c_{L,j}\partial_{p^e{\mathbf e}_j}h_l\right)\left(\sum_{|T'| = u 
- 1}c^{T'}\partial_{p^eT'}\beta\right) \\ 
&=& \left(D_u\beta\right) h_l + 
\delta_{e,e_l}\left(CM\right)_{L,l}D_{u-1}\beta \\ 
&=& \left(D_u\beta\right) h_l + \delta_{L,l}D_{u-1}\beta.\\ 
\end{eqnarray*} 
Therefore, we conclude 
$$D_u\left(\beta h_l\right) \equiv \left(D_u \beta\right) h_l + 
\delta_{L,l}D_{u - 1}\beta\mod{\fm^{r + p^{e_l} - up^e + 1}}.$$ 
This completes the proof of Lemma \ref{4.1.2.1}. 
\end{proof} 
The next lemma computes the coefficient of $h_L$, using the differential 
operator constructed in the previous lemma, in terms of the coefficients of 
the other elements $h_l \ (l \neq L)$ and terms 
of higher multiplicity. 
\begin{lem}[Supporting Lemma 2]\label{4.1.2.2} 
Let $v, s$ be integers such 
that 
$$1 \leq v < p^{e' - e} \quad \text{\ and\ } \quad 0 \leq s.$$ 
Define 
$$F_v = \sum_{u = 1}^v (-1)^uh_L^{u-1}D_u.$$ 
Suppose that the elements 
$\alpha, \beta_1, \ldots, \beta_N \in R$ satisfy 
the following conditions: 
$$\left\{\aligned 
&\alpha + \sum_{l = 1}^N\beta_l h_l \in \fm^{s + 1} \\ 
&\on{ord}_P(\beta_l) \geq s - p^{e_l} \quad (1 \leq l \leq N).\\ 
\endaligned\right.$$ 
Then we have 
$$\beta_L \equiv F_v\alpha + (- 1)^vh_L^vD_v\beta_L 
+ \sum_{\substack{1 \leq l \leq N,\\ l \neq L}} 
\left(F_v\beta_l\right) h_l\mod{\fm^{s - p^e + 1}}.$$ 
(We use the convention that $\fm^n = R$ when $n \leq 0$.) 
\end{lem} 
\begin{proof} 
From Supporting Lemma 1 it follows that for $1 \leq l 
\leq N$ 
$$D_u\left(\beta_l h_l\right) \equiv \left(D_u\beta_l\right) h_l + 
\delta_{L,l}D_{u-1}\beta_l\mod{\fm^{s - up^e + 1}}.$$ 
Multiplying by $(-1)^uh_L^{u-1}$ and adding them up with $u$ ranging from 
$1$ to $v$, we obtain 
\begin{eqnarray*} 
F_v\left(\beta_l h_l\right) &\equiv& \left(F_v\beta_l\right) h_l + 
\delta_{L,l}\sum_{u = 1}^v (-1)^uh_L^{u-1}D_{u-1}\beta_l 
\mod{\fm^{s - p^e + 1}} \\ 
&=& \left(F_v\beta_l\right) h_l - \delta_{L,l}\sum_{u = 0}^{v-1} 
(-1)^uh_L^uD_u\beta_l. 
\end{eqnarray*} 
Since $\alpha + \sum_{l = 1}^N\beta_l h_l \in \fm^{s + 1}$, we have 
$F_v\left(\alpha + \sum_{l = 1}^N\beta_l h_l\right) \in 
\fm^{s + 1 - p^e}$. 
 
Therefore, we obtain 
\begin{eqnarray*} 
F_v\alpha &\equiv& - F_v\left(\sum_{l = 1}^N\beta_l h_l\right) 
\mod{\fm^{s + 1 -p^e}} 
\\ 
&\equiv& - \sum_{l = 1}^N\left\{\left(F_v\beta_l\right) h_l 
-\delta_{L,l}\sum_{u = 0}^{v-1} (-1)^uh_L^uD_u\beta_l\right\} 
\mod{\fm^{s + 1 - p^e}} \\ 
&=& - \sum_{l = 1}^N\left(F_v\beta_l\right) h_l + \sum_{u = 0}^{v-1} 
(-1)^uh_L^uD_u\beta_L. \\ 
\end{eqnarray*} 
Therefore, we conclude 
\begin{eqnarray*} 
\lefteqn{F_v\alpha + \sum_{\substack{1 \leq l \leq N,\\ l \neq L}} 
\left(F_v\beta_l\right) h_l \equiv - \left(F_v\beta_L\right) h_L 
+ \sum_{u = 0}^{v-1}(-1)^uh_L^uD_u\beta_L \mod{\fm^{s + 1 - p^e}}
}\\ 
&=& - \sum_{u = 1}^v(-1)^uh_L^uD_u\beta_L + \sum_{u = 0}^{v-1} 
(-1)^uh_L^uD_u\beta_L 
= D_0\beta_L - (-1)^vh_L^vD_v\beta_L. 
\end{eqnarray*} 
Since $D_0 = \on{id}_R$, we conclude 
$$\beta_L \equiv F_v\alpha + (- 1)^vh_L^vD_v\beta_L 
+ \sum_{\substack{1 \leq l \leq N,\\ l \neq L}} 
\left(F_v\beta_l\right) h_l \mod{\fm^{s - p^e + 1}}.$$ 
This completes the proof of Lemma \ref{4.1.2.2}. 
\end{proof} 
The next lemma shows that, given a linear combination of 
$(h_1, \ldots, h_L)$, we can retake the coefficients so 
that they have the expected multiplicities. 
This paves the way to the coefficient lemma in the next 
subsection, where we extract more information on the coefficients when 
$(h_l,p^{e_l}) \quad (l = 1, \ldots, N)$ are in a 
($\fD$-saturated) idealistic filtration. 
\begin{lem}[Supporting Lemma 3]\label{4.1.2.3} 
We have 
$$\left(\sum_{l = 1}^N Rh_l\right) \cap \fm^r = \sum_{l = 1}^N 
\fm^{r - p^{e_l}}h_l \text{\ for\ any\ }r \in 
\bZ_{\geq 0}.$$ (We use the convention that $\fm^n = R$ when 
$n \leq0$.) 
\end{lem} 
\begin{proof} 
We have only to show the inclusion 
$$(\diamond) \quad (\sum_{l = 1}^N Rh_l) \cap 
\fm^r \subset \sum_{l = 1}^N \fm^{r - p^{e_l}}h_l,$$ 
while the opposite one is clear. 
 
We prove the inclusion by induction on the triplet 
$(\chi,L,r)$ where 
$$\chi = \{e_l\mid l = 1, \ldots, N\},$$ 
and where the set of the triplets is endowed 
with the lexicographical order. 
\begin{case}{$(\chi, L) = (1,1)$, i.e., $N = 1$.} 
In this case, take 
$\beta h_1 \in (Rh_1) \cap \fm^r$ with 
$\beta \in (\fm^r:h_1)$.  Then  since $h_1 \not\in 
\fm^{p^e + 1}$, we have $\beta \in \fm^{r - p^e}$.  Thus we have 
$$\left(Rh_1\right) \cap \fm^r = (\fm^r:h_1)h_1 \subset 
\fm^{r - p^e}h_1,$$ 
which shows the inclusion $(\diamond)$.  (Note that the inclusion 
$(\diamond)$ holds even when $r < p^e$.) 
\end{case} 
\begin{case}{$(\chi,L) > (1,1)$,\  $r \leq p^{e_N}$.} 
In this case, set $M = \min\{l\mid e_l = 
e_N\}$.  Since $r \leq p^{e_N}$, we observe 
$$(\star) \quad \sum_{l = M}^N Rh_l = \sum_{l = M}^N \fm^{r - 
p^{e_l}}h_l \subset \fm^r.$$ 
 
Assume $\chi = 1$. 
Then we have $M = 1$, and $(\star)$ implies the inclusion $(\diamond)$ 
immediately. 
 
Assume $\chi > 1$. 
Then we we have 
\begin{eqnarray*} 
\left(\sum_{l = 1}^N Rh_l\right) \cap \fm^r &=& 
\left(\sum_{l = 
1}^{M-1} Rh_l + \sum_{l = M}^N Rh_l\right) \cap \fm^r\\ 
&=& \left(\sum_{l = 1}^{M-1} Rh_l\right) \cap \fm^r + \sum_{l = M}^N 
\fm^{r - p^{e_l}}h_l \quad (\text{by\ }(\star))\\ 
&\subset& \sum_{l = 1}^{M-1}\fm^{r - p^{e_l}}h_l + \sum_{l = M}^N 
\fm^{r - p^{e_l}}h_l = \sum_{l = 1}^N \fm^{r - 
p^{e_l}}h_l, 
\end{eqnarray*} 
which implies the inclusion $(\diamond)$.  Note that the inclusion on the 
third line is obtaind by induction on $\chi$, since 
$$\#\left\{e_l\mid 1 \leq l \leq M - 1\right\} = \chi - 
1.$$ 
\end{case} 
\begin{case}{$(\chi,L) > (1,1)$,\ $r > p^{e_N}$.} 
Note that this case happens only when we are 
in positive characteristic 
$0 < p = \on{char}(k) < \infty$. 
In this case, we take an element 
$$g = \sum_{l = 1}^N \beta_lh_l \in \left(\sum_{l = 1}^NRh_l\right) \cap 
\fm^r \subset \left(\sum_{l = 1}^N Rh_l\right) \cap 
\fm^{r - 1}.$$ By induction on $r$, we may assume 
$$\beta_l \in \fm^{r - 1 - p^{e_l}} \quad (1 \leq l \leq N).$$ 
By applying Supporting Lemma 2 with 
$$0 < p^{e' - e} - 1 = v, \quad 0 \leq r - 1 = s, \quad \alpha = 
0,$$ 
as we check the conditions 
$$\alpha + \sum_{l = 1}^N\beta_l h_l\in \fm^{s + 1}, \quad 
\on{ord}_P(\beta_l) \geq s - p^{e_l} \quad 
(1 \leq l \leq N),$$ 
we conclude 
$$\beta_L \equiv F_v\alpha + (-1)^vh_L^vD_v\beta_L 
+ \sum_{\substack{1 \leq l \leq N,\\ l \neq L}} 
\left(F_v\beta_l\right)h_l \mod{\fm^{s - p^e + 1}}.$$ 
Since $F_v\alpha = 0$, we conclude 
$$\beta_L \in Rh_L^{p^{e'-e} - 1} 
+ \sum_{\substack{1 \leq l \leq N,\\ l \neq L}} 
Rh_l +\fm^{r - p^e}.$$ 
Therefore, we have 
\begin{eqnarray*} 
g &=& \sum_{l = 1}^N \beta_lh_l = 
\sum_{\substack{1 \leq l \leq N,\\ l \neq L}} 
\beta_lh_l + \beta_Lh_L \\ 
&\in& 
\left(\sum_{\substack{1 \leq l \leq N,\\ l \neq L}} 
Rh_l + Rh_L^{p^{e'-e}} + 
\fm^{r - p^e}h_L\right) \cap \fm^r \\ 
&=& \left( 
\sum_{\substack{1 \leq l \leq N,\\ l \neq L}} 
Rh_l + Rh_L^{p^{e'-e}}\right) \cap 
\fm^r + \fm^{r - p^e}h_L. 
\end{eqnarray*} 
Now instead of looking at the original 
$$\cH = \{h_1, \ldots, h_{L-1}, h_L, h_{L+1}, \ldots, h_N\} 
\quad\text{with\ }(\chi,L),$$ 
we look at 
$$\cH' = \{h_1, \ldots, h_{L-1}, h_L' = h_L^{p^{e'-e}},h_{L+1}, 
\ldots, h_N\}\quad \text{with\ }(\chi',L').$$ 
If $L = 1$, then we have $\chi' = \chi - 1$.  If $L > 1$, then we have 
$(\chi',L') = (\chi,L-1)$. Hence we always have $(\chi',L') 
< (\chi,L)$.  Therefore, by induction, we conclude 
$$\left( 
\sum_{\substack{1 \leq l \leq N,\\ l \neq L}} 
Rh_l + Rh_L^{p^{e'-e}}\right) \cap 
\fm^r \subset \sum_{\substack{1 \leq l \leq N,\\ 
l \neq L}}\fm^{r - p^{e_l}}h_l 
+ \fm^{r - p^{e'}}h_L^{p^{e'-e}}.$$ 
Plugging in this inclusion for the third line of the 
analysis for $g$, we conclude 
\begin{eqnarray*} 
g &\in& \left( 
\sum_{\substack{1 \leq l \leq N,\\ l \neq L}} 
Rh_l + Rh_L^{p^{e'-e}}\right) 
\cap \fm^r + \fm^{r - p^e}h_L. \\ 
&=& 
\sum_{\substack{1 \leq l \leq N,\\ l \neq L}} 
\fm^{r - p^{e_l}}h_l + \fm^{r - p^{e'}}h_L^{p^{e'-e}} + 
\fm^{r - p^e}h_L \\ 
&=& 
\sum_{\substack{1 \leq l \leq N,\\ l \neq L}} 
\fm^{r - p^{e_l}}h_l + \fm^{r - p^e}h_L = \sum_{l = 1}^N 
\fm^{r - p^{e_l}}h_l. 
\end{eqnarray*} 
Since $g \in \left(\sum_{l = 1}^N Rh_l\right) \cap \fm^r$ is 
arbitrary, we have the inclusion 
$$(\diamond) \quad \left(\sum_{l = 1}^N Rh_l\right) \cap \fm^r 
\subset \sum_{l = 1}^N \fm^{r - 
p^{e_l}}h_l.$$ 
\end{case} 
This completes the proof of Lemma \ref{4.1.2.3}. 
\end{proof} 
\end{subsection} 
\begin{subsection}{Setting for the coefficient lemma.}\label{4.1.3} 
We describe the setting for the coefficient lemma: 
 
Let $\bI$ be a $\fD$-saturated idealistic filtration over 
$R$. 
 
Let $\cH = 
\{h_1, \ldots, h_N\} \subset R$ be a subset of $R$, and let $0 \leq e_1 \leq 
\cdots \leq e_N$ be nonnegative integers associated to 
these elements, satisfying conditions (i) and (ii) as described in 
Setting \ref{4.1.1}, and satisfying one more condition 
\begin{enumerate} 
\item[(iii)] 
$(h_l,p^{e_l}) \in \bI$ for $l = 1, \ldots, N$. 
\end{enumerate} 
We denote by $(\cH) \subset R$ the ideal generated 
by the set $\cH$. 
 
For $f \in R$, set 
$$\on{ord}_\cH(f) = \sup\left\{n \in \bZ_{\geq0} 
\mid f \in \fm^n + (\cH)\right\}$$ 
and 
$$\mu_\cH(\bI) := \inf\left\{\mu_\cH(f,a) 
:= \frac{\on{ord}_\cH(f)}{a}\mid (f,a) \in \bI, 
a > 0\right\}.$$ 
Note that we set $\on{ord}_\cH(0) = \infty$ 
by definition. 
 
We also bring the attention of the reader to the following notation: 
\begin{quote} 
For $B = (b_1, \ldots, b_N) \in \bZ_{\geq 0}^N$, we set 
$[B] = (b_1p^{e_1}, \ldots, b_Np^{e_N})$ and hence 
$\displaystyle |[B]| = \sum_{l =1}^Nb_lp^{e_l}$. 
\end{quote} 
\end{subsection} 
\begin{subsection}{Statement and proof of 
the coefficient lemma.}\label{4.1.4} 
\begin{lem}[Coefficient Lemma]\label{4.1.4.1} 
Let $\mu \in \bR_{\geq 0}$ 
be a nonnegative number such that 
$\mu < \mu_\cH(\bI)$. 
Set 
$$\bI'_t = \bI_t \cap \fm^{\lceil \mu t\rceil},$$ 
where we use the convention that $\fm^n = R$ for $n \leq 0$. 
Then for any $a \in \bR$, we have 
$$\bI_a = \sum_B \ \bI'_{a - |[B]|}H^B.$$ 
\end{lem} 
\begin{proof} 
We set 
$${\mathfrak q}_a = \sum_B \ \bI'_{a - |[B]|}H^B \subset 
\bI_a.$$ 
Our goal is to show $\bI_a = {\mathfrak q}_a$. 
 
When $a \leq 0$, since $R = \bI_a' \subset {\mathfrak q}_a$, we have $ 
\bI_a = R = {\mathfrak q}_a$, the desired equality. 
 
Therefore, in the following, we assume $a > 0$. 
\begin{step}{Proof for the inclusion $(\star)_{c,r}$ defined below.} 
For $c \in \bZ_{> 0}$ and $r \in \bZ_{\geq 0}$, 
we set 
$$J_{c,r} = \fm^{r + 1} + {\mathfrak q}_a + \sum_{|[B]| \geq c} 
\fm^{r - |[B]|}H^B.$$ 
We prove the inclusion 
$$(\star)_{c,r} \quad \bI_a \cap \fm^r \subset J_{c,r} 
\quad (1 \leq c \leq \lceil a \rceil, r \in 
\bZ_{\geq 0})$$ 
by induction on $c$. 
\begin{case}{$c = 1$.} 
In this case, if $\lceil \mu a\rceil \leq r$, then the inclusion 
$(\star)_{1,r}$ holds since 
$$\bI_a \cap \fm^r \subset \bI_a \cap \fm^{\lceil \mu 
a\rceil} = \bI_a' \subset {\mathfrak q}_a \subset 
J_{1,r}.$$ 
If $\lceil \mu a\rceil \geq r + 1$, then we have 
\begin{eqnarray*} 
\bI_a \cap \fm^r 
&\subset& \left(\fm^{\lceil \mu a\rceil} + 
(\cH)\right) \cap \fm^r \quad (\text{by\ 
definition\ of\ }\mu_\cH(\bI) \text{\ and\ }\mu)\\ 
&\subset& \fm^{r+1} + \left(\sum_{l = 1}^NR_lh_l\right) \cap 
\fm^r \quad (\text{since\ }\fm^{\lceil \mu 
a\rceil} 
\subset \fm^{r+1} \subset \fm^r) \\ 
&=& \fm^{r+1} + \sum_{l = 1}^N \fm^{r - p^{e_l}}h_l \quad 
(\text{by\ Supporting\ Lemma\ 3}) \\ 
&=& \fm^{r+1} + \sum_{|B| = 1}\fm^{r - |[B]|}H^B \subset 
\fm^{r+1} + \sum_{|[B]| \geq 1}\fm^{r - |[B]|}H^B 
\subset J_{1,r}, 
\end{eqnarray*} 
and hence the inclusion $(\star)_{1,r}$. 
\end{case} 
\begin{case}{$c \geq 2$ assuming the inclusion 
$(\star)_{c-1, r}$.} 
Using the inclusion $(\star)_{c-1,r}$, we have 
\begin{eqnarray*} 
\bI_a \cap \fm^r 
&\subset& \left(\fm^{r+1} + {\mathfrak q}_a + 
\sum_{|[B]| \geq c - 1}\fm^{r - |[B]|}H^B\right) 
\cap \bI_a \\ 
&=& {\mathfrak q}_a + \left(\fm^{r+1} + \sum_{|[B]| \geq c - 1}\fm^{r 
- |[B]|}H^B\right)\cap \bI_a. 
\end{eqnarray*} 
Since ${\mathfrak q}_a \subset J_{c,r}$, in order to show the inclusion 
$(\star)_{c,r}$, we have only to prove 
$$\left(\fm^{r+1} + \sum_{|[B]| \geq c - 1}\fm^{r - 
|[B]|}H^B\right) 
\cap \bI_a \subset J_{c,r}.$$ 
Let $f$ be an element in the left-hand side of the desired inclusion above, 
so that there exists a finite set 
$$\left\{\alpha_B \in \fm^{r - |[B]|}\mid |[B]| 
\geq c - 1\right\} \subset R$$ 
such that 
$$f - \sum_{|[B]| \geq c - 1}\alpha_BH^B \in \fm^{r+1}.$$ 
Fix a multi-index $B_o$ with $|[B_o]| = c - 1$. 
 
Choose a regular system of parameters $(x_1, \ldots, x_d)$ such that 
$$h_l - x_l^{p^{e_l}} \in \fm^{p^{e_l} + 1} 
\quad (1 \leq l \leq N).$$ 
The partial derivatives in the following computation are 
taken with respect to this regular system of parameters 
$X = (x_1, \ldots,x_d)$.  We use the abbreviation 
$\partial_J = \partial_{X^J}$.  The symbol 
``$\equiv$'' denotes an equality modulo 
$\fm^{r -(c-1) + 1} = \fm^{r - c + 2}$. 
We compute 
\begin{eqnarray*} 
\partial_{[B_o]}f &\equiv& \sum_{|[B]| \geq c -1} 
\partial_{[B_o]}(\alpha_BH^B) \\ 
&=& \sum_{|[B]| \geq c - 1} \sum_{J \leq [B_o]}(\partial_{[B_o] - 
J}\alpha_B)(\partial_JH^B) \quad 
(\text{by\ the\ generalized\ product\ rule})\\ 
&\equiv& \sum_{|[B]| \geq c - 1} \sum_{J \leq 
[B_o]}(\partial_{[B_o] - J}\alpha_B)(\partial_JX^{[B]}) 
\\ 
&=& \sum_{|[B]| \geq c - 1} \sum_{J \leq [B_o]}(\partial_{[B_o] - 
J}\alpha_B)\binom{[B]}{J}X^{[B] - J} 
\end{eqnarray*} 
In the last formula, the binomial coefficient $\binom{[B]}{J}$ is zero 
unless $J = [K]$ for some $K \leq B_o$.  Therefore, we have 
\begin{eqnarray*} 
\partial_{[B_o]}f 
&\equiv& \sum_{|[B]| \geq c - 1} \sum_{K \leq 
B_o}(\partial_{[B_o - K]}\alpha_B)\binom{[B]}{[K]}X^{[B - K]} 
\\ 
&=& \sum_{|[B]| \geq c - 1} \sum_{K \leq B_o}(\partial_{[B_o - 
K]}\alpha_B)\binom{B}{K}X^{[B - K]} \\ 
&\equiv& \sum_{|[B]| \geq c - 1} \sum_{K \leq B_o}(\partial_{[B_o - 
K]}\alpha_B)\binom{B}{K}H^{B - K}. 
\end{eqnarray*} 
In the last formula, the binomial coefficient 
$\binom{B}{K} = 0$ unless $K \leq B$. 
 
If $K < B$, we have $|[B - K]| \geq 1$ and 
$$\partial_{[B_o - K]}\alpha_B \in \fm^{r - |[B]| - |[B_o - K]|} = 
\fm^{r - (c-1) - |[B - K]|}.$$ 
If $K = B$, we have $B = B_o$, since 
$B = K \leq B_o$ and $|[B_o]| = c - 1\leq |[B]|$. 
 
Therefore, we have 
\begin{eqnarray*} 
(*) \quad \partial_{[B_o]}f - \alpha_{B_o} &\in& \sum_{K < B} 
\fm^{r - c + 1 - |[B - K]|}H^{B - K} + \fm^{r - c + 2} 
\\ 
&=& \sum_{|[B]| \geq 1}\fm^{r - c + 1 - |[B]|}H^B + \fm^{r - c 
+ 2}. 
\end{eqnarray*} 
On the other hand, since $f \in \bI_a \cap \fm^r$ and since the 
idealistic filtration $\bI$ is $\fD$-saturated, we have 
$$\partial_{[B_o]}f \in \bI_{a - (c-1)} \cap \fm^{r - (c-1)} 
=\bI_{a - c + 1} \cap \fm^{r - c + 1}.$$ 
Using the inclusion $(\star)_{1,r-c+1}$, we obtain 
$$(**) \quad \partial_{[B_o]}f \in \bI_{a - c + 1} \cap 
\fm^{r - c + 1} \subset 
\fm^{r - c + 2} + {\mathfrak q}_{a-c+1} + \sum_{|[B]| \geq 1}\fm^{r 
- c + 1 -|[B]|}H^B.$$ 
From $(*)$ and $(**)$ it follows that 
\begin{eqnarray*} 
\alpha_{B_o}H^{B_o} &\in& 
\fm^{r - c + 2}H^{B_o} + {\mathfrak q}_{a - c + 
1}H^{B_o} + \sum_{|[B]| \geq 1}\fm^{r - c + 1 - 
|[B]|}H^{B+B_o} \\ 
&\subset& \fm^{r+1} + {\mathfrak q}_a 
+ \sum_{|[B + B_o]| \geq c} 
\fm^{r - c + 1 -|[B]|}H^{B+B_o} \subset J_{c,r}. 
\end{eqnarray*} 
 
Since $B_o$ is arbitrary with $|[B_o]| = c - 1$, we conclude that 
$\alpha_BH^B \in J_{c,r}$ for all $B$ with $|[B]| = c - 1$. 
Therefore, we have 
$$f \in \sum_{|[B]| \geq c-1}\alpha_BH^B + \fm^{r+1} = \sum_{|[B]| = 
c-1}\alpha_BH^B + \left(\fm^{r+1} + \sum_{|[B]| 
\geq c}\alpha_BH^B\right) \subset J_{c,r},$$ 
which implies the desired inclusion $(\star)_{c,r}$. 
 
This completes the proof for the inclusion $(\star)_{c,r}$. 
\end{case} 
\end{step} 
\begin{step}{Finishing argument.} 
We finish the proof of Coefficient Lemma using the result 
of Step 1. 
 
Applying the inclusion 
$(\star)_{\lceil a\rceil,r}$ for $r \in \bZ_{\geq 0}$, 
we have 
$$\bI_a \cap \fm^r \subset \fm^{r + 1} + {\mathfrak q}_a + 
\sum_{|[B]| \geq 
\lceil a \rceil}\fm^{r - |[B]|}H^B = \fm^{r+1} + {\mathfrak q}_a,$$ 
since $\bI_{a - |[B]|}' = R$ for $B$ with $|[B]| \geq \lceil a\rceil$. 
 
Therefore, we have 
$$\bI_a \cap \fm^r \subset \bI_a \cap \left(\fm^{r+1} 
+ {\mathfrak q}_a\right) = \bI_a \cap \fm^{r+1} 
+ {\mathfrak q}_a,$$ 
which implies 
$$\bI_a \cap \fm^r + {\mathfrak q}_a = \bI_a \cap 
\fm^{r+1} + {\mathfrak q}_a,$$ 
for any $r \in \bZ_{\geq 0}$.  In particular, we have 
$$\bI_a = \bI_a \cap \fm^0 + {\mathfrak q}_a = \bI_a \cap 
\fm^{\lceil \mu a\rceil} + {\mathfrak q}_a = 
\bI_a' + {\mathfrak q}_a = {\mathfrak q}_a.$$ 
\end{step} 
This completes the proof of Lemma \ref{4.1.4.1}. 
\end{proof} 
\begin{rem}\label{4.1.4.2} 
\item[(1)] 
The purpose of having a nonnegative number $\mu \in \bR_{\geq 0}$ 
with $\mu < \mu_\cH(\bI)$ involved in 
our statement of Lemma \ref{4.1.4.1} 
is to make it valid even when 
$\mu_\cH(\bI) = \infty$, the case to which 
we often apply Coefficient Lemma. 
When $\mu_\cH(\bI) < \infty$, we may of course 
apply Coefficient Lemma, setting $\mu = \mu_\cH(\bI)$. 
\item[(2)] 
We can restrict the range of $B$ in the expression $\bI_a = 
\sum_B\bI_{a - |[B]|}'H^B$ to a specific finite range, 
e.g., 
$B$ with 
$|[B]|  < a + p^{e_N}$.  In fact, if $|[B]| \geq a + p^{e_N}$, there exists 
$B' < B$ such that $a \leq |[B']| < a + p^{e_N}$.  Then 
we have 
$\bI_{a - |[B]|}'H^B \subset RH^{B'} = \bI_{a - |[B']|}'H^{B'}$. 
Therefore, if $B$ is out of this range, the term $ 
\bI_{a - |[B]|}'H^B$ is redundant, i.e., 
$$\sum_B\bI_{a - |[B]|}'H^B = 
\sum_{|[B]| < a + p^{e_N}}\bI_{a - |[B]|}'H^B.$$ 
\item[(3)] 
In Part II, given an element $(f,a) \in \bI$ of a 
$\fD$-saturated idealistic filtration, we analyze ``the power 
series expansion of $f$ with respect to a set $\cH$ satisfying 
conditions (i), (ii), (iii) (e.g.~a leading generator system of 
$\bI$)''.  This provides a different approach to 
Coefficient Lemma and an alternative proof. 
\end{rem} 
\end{subsection} 
\end{section} 
\begin{section}{Nonsingularity principle.}\label{4.2} 
\begin{subsection}{Statement of the 
nonsingularity principle.}\label{4.2.1} 
\begin{thm}\label{4.2.1.1} 
Let $\bI$ be an idealistic filtration which 
is $\fB$-saturated.  Let $\cH = \{h_1, \ldots, 
h_N\} \subset R$ be a subset of $R$, and let $0 \leq e_1 
\leq\cdots\leq e_N$ be nonnegative integers associated to 
these elements, satisfying 
conditions (i), (ii), (iii) as described in 
Setting \ref{4.1.3}. 
Suppose $\mu_\cH(\bI) = \infty$.  Then 
\item[(1)] 
$\bH = \{(h_l,p^{e_l})\mid l = 1, \ldots, N\}$ 
generates the idealistic filtration $\bI$, i.e., 
$$\bI = G(\bH).$$ 
\item[(2)] 
The elements in $\bH$ are all concentrated at level $p^0 = 1$, 
i.e., 
$$\bH \subset R \times \{1\}.$$ 
(Note that in characteristic zero, where we take $p = \infty$ 
according to our convention, we set $p^0 = \infty^0 =1$. 
cf.~\ref{0.2.3.2.1}.) 

Consequently, we conclude that the support of the 
idealistic filtration is defined by $\cH$, 
i.e., $\on{Supp}(\bI) =V(\cH)$, 
and hence that it is nonsingular. 
\end{thm} 
\begin{rem}\label{4.2.1.2} 
\item[(1)] 
In Theorem \ref{4.2.1.1}, we see from assertion (1) 
that 
$$\{\overline{h_l} =(h_l\bmod{\fm^{p^{e_l}}}) 
\mid l = 1, \ldots, N\}$$ 
generates $L(\bI)$ 
(cf.~Definition \ref{3.1.1.1}), and hence conclude 
that $\bH$ is a leading generator system, even though 
we do not a priori assume so. 
\item[(2)] 
In Part II, we will look at the invariant $\widetilde{\mu}$, which is a 
priori defined to be $\widetilde{\mu} = \mu_ 
\cH(\bI)$ with respect to the set $\cH$ associated to a leading 
generator system.  We will see, however, that $\mu_ 
\cH(\bI)$ is independent of the choice of a leading generator system, 
and hence that $\widetilde{\mu}$ is actually an invariant 
intrinsic to the idealistic filtration $\bI$.  Therefore, the 
nonsingularity principle above can be regarded as the description 
of an idealistic filtration with $\widetilde{\mu} = \infty$, with the 
conclusions holding for {\it any} leading generator system $\bH$. 
\item[(3)] 
Recall that, as we construct the strand of invariants in our algorithm, 
we enlarge the idealistic filtration and construct 
its modifications (cf.~\ref{0.2.3.2.2} and \ref{0.2.3.2.4}). 
At the end of the construction of the strand of invariants, we 
reach the last modification, which is an idealistic filtration 
(which is both $\fR$-saturated and $\fD$-saturated) 
whose leading generator system satisfies the conditions described 
in the above.  The 
maximum locus of the strand of invariants, which we take as 
the center of blowup, coincides with the support of this 
last modification 
(in a neighborhood of each point of the maximum locus), 
and hence is nonsingular according to Theorem \ref{4.2.1.1}. 
This is why it is 
called the nonsingularity principle of the center. 
\item[(4)] 
In order to show $\bI = G(\bH)$, we only need 
$\bI$ to be $\fD$-saturated, while in order to 
show $\bH \subset R\times \{1\}$, 
we need $\bI$ to be $\fB$-saturated. 
\end{rem} 
\end{subsection} 
\begin{subsection}{Proof of the nonsingularity 
principle.}\label{4.2.2} 
\begin{proof}[\underline{\it Proof for assertion $\on{(1)}$}] 
\item
We show that $\bH$ 
generates the idealistic filtration 
$\bI$, i.e., $\bI = G(\bH)$. 
 
Since $\mu_\cH(\bI) = \infty$, 
we can apply Coefficient Lemma with 
an arbitrary non-negative number 
$\bZ_{\geq0}\ni\mu<\mu_\cH(\bI) = \infty$ and obtain 
\begin{eqnarray*} 
\bI_a &=& \sum_B\bI'_{a - |[B]|}H^B = \sum_{|[B]| \geq a} 
\bI'_{a - |[B]|}H^B + \sum_{|[B]| < a}\bI'_{a - 
|[B]|}H^B \\ 
&\subset& \sum_{|[B]| \geq a}RH^B + \bI_{a - (\lceil a\rceil - 1)}' 
\subset \sum_{|[B]| \geq a}RH^B + \fm^{\lceil \mu(a - 
\lceil a \rceil + 1) \rceil}. 
\end{eqnarray*} 
Since $a - \lceil a \rceil + 1 > 0$, this implies by Krull's intersection 
theorem that 
$$\bI_a \subset \bigcap_{r = 0}^{\infty}\left(\sum_{|[B]| \geq a}RH^B + 
\fm^r\right) = \sum_{|[B]| \geq a}RH^B.$$ 
This shows that $\bH$ generates $\bI$, i.e., $\bI = G( 
\bH)$. 
\end{proof} 
\begin{proof}[\underline{\it Proof for assertion $\on{(2)}$}] 
\item
We show that the elements in $\bH$ are concentrated at level $p^0 = 1$, 
i.e., $\bH \subset R \times \{1\}$.  Set 
$$\bH_0 = \{(h_l,p^{e_l}) \in \bH\mid e_l = 0\} 
= \bH \cap (R \times \{1\}).$$ 
 
We will derive a contradiction assuming $\bH_0 \neq \bH$. 
Set $e = \min\{e_l\mid e_l > 0\}$. 
\begin{step} 
{We show that $\bI_a = (\cH)$ for  $0 < a \leq 1$ 
and that $(\cH) = \sqrt{(\cH)}$.} 
In fact, since $\bI = G(\bH)$ and since $\bH \subset R \times 
\bR_{\geq 1}$, Lemma \ref{Construction} (1) implies that 
$$\bI_a = \sum_{l = 1}^NRh_l = (\cH) \quad \text{\ for\ } 
\quad 0 < a \leq 1.$$ 
Suppose $g \in \sqrt{(\cH)}$, i.e., $g^n \in (\cH) = \bI_1$ 
for some $n \in \bZ_{> 0}$.  Since $\bI$ is 
$\fR$-saturated, this implies $g \in \bI_{1/n} = (\cH)$. 
Therefore, we have $(\cH) = \sqrt{(\cH)}$. 
\end{step} 
\begin{step}{We show that 
$(\cH) = ((\cH) \cap R^{p^e}) + (\cH_0)$.} 
Set 
$${\mathcal D} = \{d \in \on{Diff}^{p^e-1}_R\mid 
d((\cH_0)) \subset (\cH_0)\}.$$ 
Observe 
$$(*) \quad {\mathcal D}((\cH)) \subset (\cH).$$ 
In fact, for $d \in {\mathcal D}$, since $\cH \setminus 
\cH_0 \subset \bigcup_{e_l > 0}\bI_{p^{e_l}} = \bI_{p^e}$ 
and since $\bI$ is $\fD$-saturated, we have 
$$d\left((\cH \setminus \cH_0)\right) \subset d\left( 
\bI_{p^e}\right) \subset \bI_{p^e - (p^e - 1)} = \bI_1 
= (\cH).$$ Therefore, we conclude 
$$d((\cH)) = d((\cH_0)) + d\left((\cH \setminus 
\cH_0)\right) \subset (\cH_0) + (\cH) = (\cH).$$ 
Now $(*)$ implies 
$$(**) \quad \on{Diff}^{p^e-1}_{\overline{R}} 
\left(\overline{(\cH)}\right) \subset \overline{( 
\cH)} \text{\ and\ hence\ 
}\on{Diff}^{p^e-1}_{\overline{R}}\left(\overline{(\cH)}\right) = 
\overline{( 
\cH)}$$ 
where 
$$\overline{R} = R/(\cH_0) 
\quad\text{and}\quad 
\overline{(\cH)} = ( 
\cH)/(\cH_0).$$ 
Then, by 
Proposition \ref{1.3.1.2}, 
$(**)$ implies 
$$\overline{(\cH)} = \left(\overline{(\cH)} \cap 
\overline{R}^{p^e}\right).$$ 
Therefore, we have 
$$(\cH) = \left((\cH) 
\cap R^{p^e}\right) + (\cH_0).$$ 
\end{step} 
\begin{step}{Finishing argument.} 
By Step 2, we conclude 
\begin{eqnarray*} 
(\cH) &=& (\cH_0) + \left((\cH) \cap R^{p^e}\right) 
= 
(\cH_0) + \left(\left\{g^{p^e} \in (\cH)\mid 
g \in R\right\}\right) \\ 
&=& (\cH_0) + \left(\left\{g^{p^e}\mid g \in (\cH) 
\right\}\right) \quad (\text{by\ }(\cH) = 
\sqrt{(\cH)})\\ 
&=& (\cH_0) +\left(\left\{g^{p^e}\mid g \in 
(\cH \setminus \cH_0)\right\}\right) 
\subset (\cH_0) + \fm^{p^e + 1}, 
\end{eqnarray*} 
i.e., 
$(\cH) \subset (\cH_0) + \fm^{p^e + 1}$. 
\par 
Choose a regular system of parameters $(x_1, \ldots, x_d)$ so that 
$$\left\{ 
\begin{array}{lllll} 
x_l &= &h_l &\text{for\ }1 \leq l \leq L 
&\text{where\ } L =\#\{l\mid e_l = 0\}\\ 
x_l^{p^{e_l}} &\equiv &h_l &\bmod{\ \fm^{p^{e_l}+1}} 
&\text{for\ } L+1 \leq l \leq N. 
\end{array}\right.$$ 
Then the above inclusion would imply 
$$ 
\left\{(\cH) + \fm^{p^e + 1}\right\}/\fm^{p^e + 1} \subset 
\left\{(\cH_0) + \fm^{p^e + 1}\right\}/ 
\fm^{p^e + 1} \subset R/\fm^{p^e + 1} 
$$ 
and we identify 
$R/\fm^{p^e + 1}\cong k[x_1,\ldots,x_d]/(x_1,\ldots,x_d)^{p^e + 1}$. 
 
On the other hand, however, we have the following element in the first 
quotient 
$$\left(\text{the\ leading\ term\ of\ }h_{L+1}\right) = 
x_{L+1}^{p^{e_{L+1}}} = x_{L+1}^{p^e} \in \left\{(\cH) + 
\fm^{p^e + 1}\right\}/\fm^{p^e + 1},$$ 
which obviously is {\it not} in the middle quotient 
$$(x_1, \ldots, x_L) = 
\left\{(\cH_0) + \fm^{p^e + 1}\right\}/\fm^{p^e + 1},$$ 
a contradiction ! 
 
This contradiction is derived from the assumption 
that $\bH_0 \neq\bH$. 
 
Therefore, we conclude $\bH_0 = \bH$, i.e., 
$$\bH \subset R \times \{1\}.$$ 
This completes the proof of Theorem \ref{4.2.1.1}, 
the nonsingularity principle. 
\end{step} 
\end{proof} 
\end{subsection} 
\end{section} 
\end{chapter} 


\begin{thebibliography}{HLOQ00}
\bibitem[Abh66]{MR0217069}
Shreeram~S. Abhyankar.
\newblock {\em Resolution of singularities of embedded algebraic surfaces}.
\newblock Pure and Applied Mathematics, Vol. 24. Academic Press, New York,
  1966.

\bibitem[Abh77]{MR542446}
Shreeram~S. Abhyankar.
\newblock {\em Lectures on expansion techniques in algebraic geometry},
  volume~57 of {\em Tata Institute of Fundamental Research Lectures on
  Mathematics and Physics}.
\newblock Tata Institute of Fundamental Research, Bombay, 1977.
\newblock Notes by Balwant Singh.

\bibitem[Abh83]{MR713043}
Shreeram~S. Abhyankar.
\newblock Desingularization of plane curves.
\newblock In {\em Singularities, Part 1 (Arcata, Calif., 1981)}, volume~40 of
  {\em Proc. Sympos. Pure Math.}, pages 1--45. Amer. Math. Soc., Providence,
  RI, 1983.

\bibitem[AdJ97]{MR1487237}
Dan Abramovich and A.~J. de~Jong.
\newblock Smoothness, semistability, and toroidal geometry.
\newblock {\em J. Algebraic Geom.}, 6(4):789--801, 1997.

\bibitem[AHV75]{MR0444999}
Jos{\'e}~M. Aroca, Heisuke Hironaka, and Jos{\'e}~L. Vicente.
\newblock {\em The theory of the maximal contact}.
\newblock Instituto ``Jorge Juan'' de Matem\'aticas, Consejo Superior de
  Investigaciones Cientificas, Madrid, 1975.
\newblock Memorias de Matem\'atica del Instituto ``Jorge Juan'', No.~29.
  [Mathematical Memoirs of the ``Jorge Juan'' Institute, No.~29].

\bibitem[AHV77]{MR480502}
Jos{\'e}~M. Aroca, Heisuke Hironaka, and Jos{\'e}~L. Vicente.
\newblock {\em Desingularization theorems}, volume~30 of {\em Memorias de
  Matem\'atica del Instituto ``Jorge Juan'' [Mathematical Memoirs of the Jorge
  Juan Institute]}.
\newblock Consejo Superior de Investigaciones Cient\'\i ficas, Madrid, 1977.

\bibitem[Ben70]{MR0252388}
Bruce~M. Bennett.
\newblock On the characteristic functions of a local ring.
\newblock {\em Ann. of Math. (2)}, 91:25--87, 1970.

\bibitem[BEV05]{MR2174912}
Ana~M. Bravo, Santiago Encinas, and Orlando E.~Villamayor.
\newblock A simplified proof of desingularization and applications.
\newblock {\em Rev. Mat. Iberoamericana}, 21(2):349--458, 2005.

\bibitem[Bie04]{B_BIRS}
Edward Bierstone.
\newblock Resolution of singularities.
\newblock {\em A preprint for series of lectures at the ``Workshop on 
  resolution of singularities, factorization of birational mappings, 
  and toroidal geometry'' at Banff International Research Station, 
  Banff, during 11-16 December}, 2004.

\bibitem[Bie05]{B_Harvard}
Edward Bierstone.
\newblock Functoriality in resolution of singularities.
\newblock {\em Slides for a colloquium at Harvard University on November 10},
  2005.

\bibitem[BM89]{MR1001853}
Edward Bierstone and Pierre~D. Milman.
\newblock Uniformization of analytic spaces.
\newblock {\em J. Amer. Math. Soc.}, 2(4):801--836, 1989.

\bibitem[BM90]{MR1051203}
Edward Bierstone and Pierre~D. Milman.
\newblock Local resolution of singularities.
\newblock In {\em Real analytic and algebraic geometry (Trento, 1988)}, volume
  1420 of {\em Lecture Notes in Math.}, pages 42--64. Springer, Berlin, 1990.

\bibitem[BM91]{MR1106412}
Edward Bierstone and Pierre~D. Milman.
\newblock A simple constructive proof of canonical resolution of singularities.
\newblock In {\em Effective methods in algebraic geometry (Castiglioncello,
  1990)}, volume~94 of {\em Progr. Math.}, pages 11--30. Birkh\"auser Boston,
  Boston, MA, 1991.

\bibitem[BM97]{MR1440306}
Edward Bierstone and Pierre~D. Milman.
\newblock Canonical desingularization in characteristic zero by blowing up the
  maximum strata of a local invariant.
\newblock {\em Invent. Math.}, 128(2):207--302, 1997.

\bibitem[BM99]{MR1748600}
Edward Bierstone and Pierre~D. Milman.
\newblock Resolution of singularities.
\newblock In {\em Several complex variables (Berkeley, CA, 1995--1996)},
  volume~37 of {\em Math. Sci. Res. Inst. Publ.}, pages 43--78. Cambridge Univ.
  Press, Cambridge, 1999.

\bibitem[BM03]{MR2078560}
Edward Bierstone and Pierre~D. Milman.
\newblock Desingularization algorithms. {I}. {R}ole of exceptional divisors.
\newblock {\em Mosc. Math. J.}, 3(3):751--805, 1197, 2003.
\newblock \{Dedicated to Vladimir Igorevich Arnold on the occasion of his 65th
  birthday\}.

\bibitem[Bou64]{MR0194450}
Nicolas Bourbaki.
\newblock {\em \'{E}l\'ements de math\'ematique. {F}asc. {XXX}. {A}lg\`ebre
  commutative. {C}hapitre 5: {E}ntiers. {C}hapitre 6: {V}aluations}.
\newblock Actualit\'es Scientifiques et Industrielles, No. 1308. Hermann,
  Paris, 1964.

\bibitem[BP96]{MR1397679}
Fedor~A. Bogomolov and Tony~G. Pantev.
\newblock Weak {H}ironaka theorem.
\newblock {\em Math. Res. Lett.}, 3(3):299--307, 1996.

\bibitem[BV03]{MR1971154}
Ana~M. Bravo and Orlando E.~Villamayor.
\newblock A strengthening of resolution of singularities in characteristic
  zero.
\newblock {\em Proc. London Math. Soc. (3)}, 86(2):327--357, 2003.

\bibitem[Cos87]{MR907903}
Vincent Cossart.
\newblock Forme normale d'une fonction sur un {$k$}-sch\'ema de dimension {$3$}
  et de caract\'eristique positive.
\newblock In {\em G\'eom\'etrie alg\'ebrique et applications, I (La R\'abida,
  1984)}, volume~22 of {\em Travaux en Cours}, pages 1--21. Hermann, Paris,
  1987.

\bibitem[Cos00]{MR1748622}
Vincent Cossart.
\newblock Uniformisation et d\'esingularisation 
  des surfaces d'apr\`es {Z}ariski.
\newblock In {\em Resolution of singularities (Obergurgl, 1997)}, 
  volume 181 of {\em Progr. Math.}, pages 239--258. Birkh\"auser, Basel, 2000.

\bibitem[Cos04]{C_BIRS}
Vincent Cossart.
\newblock Towards local uniformization along a valuation in Artin-Schreier
  extensions (dimension $3$).
\newblock {\em A talk at the ``Workshop on resolution of singularities,
  factorization of birational mappings, and toroidal geometry'' at Banff
  International Research Station, Banff, during 11-16 December}, 2004.

\bibitem[Cut04]{MR2058431}
Steven~Dale Cutkosky.
\newblock {\em Resolution of singularities}, volume~63 of {\em Graduate Studies
  in Mathematics}.
\newblock American Mathematical Society, Providence, RI, 2004.

\bibitem[Cut06]{AG0606530}
Steven~Dale Cutkosky.
\newblock Resolution of singularities for $3$-folds in positive characteristic.
\newblock {\em {\tt http://arxiv.org/abs/math.AG/0606530}}, 2006.

\bibitem[dJ96]{MR1423020}
A.~J. de~Jong.
\newblock Smoothness, semi-stability and alterations.
\newblock {\em Inst. Hautes \'Etudes Sci. Publ. Math.}, (83):51--93, 1996.

\bibitem[EH02]{MR1949115}
Santiago Encinas and Herwig Hauser.
\newblock Strong resolution of singularities in characteristic zero.
\newblock {\em Comment. Math. Helv.}, 77(4):821--845, 2002.

\bibitem[ENV03]{MR1974392}
Santiago Encinas, A.~Nobile, and Orlando E.~Villamayor.
\newblock On algorithmic equi-resolution and stratification of {H}ilbert
  schemes.
\newblock {\em Proc. London Math. Soc. (3)}, 86(3):607--648, 2003.

\bibitem[EV98]{MR1654779}
Santiago Encinas and Orlando E.~Villamayor.
\newblock Good points and constructive resolution of singularities.
\newblock {\em Acta Math.}, 181(1):109--158, 1998.

\bibitem[EV00]{MR1748620}
Santiago Encinas and Orlando E.~Villamayor.
\newblock A course on constructive desingularization and equivariance.
\newblock In {\em Resolution of singularities (Obergurgl, 1997)}, volume 181 of
  {\em Progr. Math.}, pages 147--227. Birkh\"auser, Basel, 2000.

\bibitem[EV03]{MR2023188}
Santiago Encinas and Orlando E.~Villamayor.
\newblock A new proof of desingularization over fields of characteristic zero.
\newblock In {\em Proceedings of the International Conference on Algebraic
  Geometry and Singularities (Spanish) (Sevilla, 2001)}, volume~19, pages
  339--353, 2003.

\bibitem[Gir74]{MR0460712}
Jean Giraud.
\newblock Sur la th\'eorie du contact maximal.
\newblock {\em Math. Z.}, 137:285--310, 1974.

\bibitem[Gir75]{MR0384799}
Jean Giraud.
\newblock Contact maximal en caract\'eristique positive.
\newblock {\em Ann. Sci. \'Ecole Norm. Sup. (4)}, 8(2):201--234, 1975.

\bibitem[Gir83]{MR734215}
Jean Giraud.
\newblock Forme normale d'une fonction sur une surface de caract\'eristique
  positive.
\newblock {\em Bull. Soc. Math. France}, 111(2):109--124, 1983.

\bibitem[Gro67]{MR0238860}
Alexander Grothendieck.
\newblock \'{E}l\'ements de g\'eom\'etrie alg\'ebrique. {IV}. \'{E}tude locale
  des sch\'emas et des morphismes de sch\'emas {IV}.
\newblock {\em Inst. Hautes \'Etudes Sci. Publ. Math.}, (32):361, 1967.

\bibitem[Hau98]{MR1652479}
Herwig Hauser.
\newblock Seventeen obstacles for resolution of singularities.
\newblock In {\em Singularities (Oberwolfach, 1996)}, volume 162 of {\em Progr.
  Math.}, pages 289--313. Birkh\"auser, Basel, 1998.

\bibitem[Hau03]{MR1978567}
Herwig Hauser.
\newblock The {H}ironaka theorem on resolution of singularities (or: {A} proof
  we always wanted to understand).
\newblock {\em Bull. Amer. Math. Soc. (N.S.)}, 40(3):323--403 (electronic),
  2003.

\bibitem[Hir63]{MR0175898}
Heisuke Hironaka.
\newblock On resolution of singularities (characteristic zero).
\newblock In {\em Proc. Internat. Congr. Mathematicians (Stockholm, 1962)},
  pages 507--521. Inst. Mittag-Leffler, Djursholm, 1963.

\bibitem[Hir64]{MR0199184}
Heisuke Hironaka.
\newblock Resolution of singularities of an algebraic variety over a field of
  characteristic zero. {I}, {II}.
\newblock {\em Ann. of Math. (2) 79 (1964), 109--203; ibid. (2)}, 79:205--326,
  1964.

\bibitem[Hir70]{MR0269658}
Heisuke Hironaka.
\newblock Additive groups associated with points of a projective space.
\newblock {\em Ann. of Math. (2)}, 92:327--334, 1970.

\bibitem[Hir72a]{MR0393555}
Heisuke Hironaka.
\newblock Gardening of infinitely near singularities.
\newblock In {\em Algebraic geometry, Oslo 1970 (Proc. Fifth Nordic Summer
  School in Math.)}, pages 315--332. Wolters-Noordhoff, Groningen, 1972.

\bibitem[Hir72b]{MR0393028}
Heisuke Hironaka.
\newblock Schemes, etc.
\newblock In {\em Algebraic geometry, Oslo 1970 (Proc. Fifth Nordic Summer
  School in Math.)}, pages 291--313. Wolters-Noordhoff, Groningen, 1972.

\bibitem[Hir77]{MR0498562}
Heisuke Hironaka.
\newblock Idealistic exponents of singularity.
\newblock In {\em Algebraic geometry (J. J. Sylvester Sympos., Johns Hopkins
  Univ., Baltimore, Md., 1976)}, pages 52--125. Johns Hopkins Univ. Press,
  Baltimore, Md., 1977.

\bibitem[Hir03]{MR1996845}
Heisuke Hironaka.
\newblock Theory of infinitely near singular points.
\newblock {\em J. Korean Math. Soc.}, 40(5):901--920, 2003.

\bibitem[Hir05]{MR2145950}
Heisuke Hironaka.
\newblock Three key theorems on infinitely near singularities.
\newblock In {\em Singularit\'es Franco-Japonaises}, volume~10 of {\em S\'emin.
  Congr.}, pages 87--126. Soc. Math. France, Paris, 2005.

\bibitem[Hir06]{H-Trieste}
Heisuke Hironaka.
\newblock A program for resolution of singularities, in all characteristics
  $p>0$ and in all dimensions.
\newblock {\em preprint for series of lectures in ``Summer School on Resolution
  of Singularities'' at International Center for Theoretical Physics, Trieste,
  during 12-30 June}, 2006.

\bibitem[HLOQ00]{MR1748614}
Herwig Hauser, Joseph Lipman, Frans Oort, and Adolfo Quir{\'o}s, editors.
\newblock {\em Resolution of singularities}, volume 181 of {\em Progress in
  Mathematics}.
\newblock Birkh\"auser Verlag, Basel, 2000.
\newblock A research textbook in tribute to Oscar Zariski, Papers from the
  Working Week held in Obergurgl, September 7--14, 1997.

\bibitem[Kol05]{AG0508332}
J\'anos Koll\'ar.
\newblock Resolution of singularities -- seattle lecture.
\newblock {\em notes for series of lectures given in the AMS Summer Institute,
  \par {\tt http://arxiv.org/abs/math.AG/0508332}}, 2005.

\bibitem[Kuh97]{Kuhlmann97}
Franz-Viktor Kuhlmann.
\newblock On local uniformization in arbitrary characteristic.
\newblock {\em The Fields Institute Preprint Series}, 1997.

\bibitem[Kuh00]{MR1748629}
Franz-Viktor Kuhlmann.
\newblock Valuation theoretic and model theoretic aspects of local
  uniformization.
\newblock In {\em Resolution of singularities (Obergurgl, 1997)}, volume
 181 of
  {\em Progr. Math.}, pages 381--456. Birkh\"auser, Basel, 2000.

\bibitem[Lip69]{MR0276239}
Joseph Lipman.
\newblock Rational singularities, with applications to algebraic surfaces and
  unique factorization.
\newblock {\em Inst. Hautes \'Etudes Sci. Publ. Math.}, (36):195--279, 1969.

\bibitem[Lip75]{MR0389901}
Joseph Lipman.
\newblock Introduction to resolution of singularities.
\newblock In {\em Algebraic geometry (Proc. Sympos. Pure Math., Vol. 29,
  Humboldt State Univ., Arcata, Calif., 1974)}, pages 187--230. Amer. Math.
  Soc., Providence, R.I., 1975.

\bibitem[Lip78]{MR0491722}
Joseph Lipman.
\newblock Desingularization of two-dimensional schemes.
\newblock {\em Ann. Math. (2)}, 107(1):151--207, 1978.

\bibitem[Lip83]{MR713245}
Joseph Lipman.
\newblock Quasi-ordinary singularities of surfaces in {${\bf C}\sp{3}$}.
\newblock In {\em Singularities, Part 2 (Arcata, Calif., 1981)}, volume~40 of
  {\em Proc. Sympos. Pure Math.}, pages 161--172. Amer. Math. Soc., Providence,
  RI, 1983.

\bibitem[LT74]{LT74}
Monique Lejeune-Jalabert and Bernard Teissier.
\newblock Cl\^oture int\'egrale des id\'eaux et \'equisingularit\'e.
\newblock {\em S\'eminaire sur les singularit\'e \`a l'Ecole Polytechnique},
  pages 1--66, 1974--1975.

\bibitem[Mat86]{MR879273}
Hideyuki Matsumura.
\newblock {\em Commutative ring theory}, volume~8 of {\em Cambridge Studies in
  Advanced Mathematics}.
\newblock Cambridge University Press, Cambridge, 1986.
\newblock Translated from the Japanese by M. Reid.

\bibitem[Mk06]{AG0103120}
Kenji Matsuki.
\newblock Resolution of singularities in characteristic zero with focus on the
  inductive algorithm by {V}illamayor and its simplification by
  {W}{\l}odarczyk.
\newblock {\em preprint, formally a revision of\ \ {\tt math.AG/0103120}, but
  completely rewritten from scratch}, 2006.

\bibitem[Moh92]{MR1186705}
T.~T. Moh.
\newblock Quasi-canonical uniformization of hypersurface singularities of
  characteristic zero.
\newblock {\em Comm. Algebra}, 20(11):3207--3249, 1992.

\bibitem[Moh96]{MR1395176}
T.~T. Moh.
\newblock On a {N}ewton polygon approach to the uniformization of singularities
  of characteristic {$p$}.
\newblock In {\em Algebraic geometry and singularities (La R\'abida, 1991)},
  volume 134 of {\em Progr. Math.}, pages 49--93. Birkh\"auser, Basel, 1996.

\bibitem[Nag57]{MR0089836}
Masayoshi Nagata.
\newblock Note on a paper of {S}amuel concerning asymptotic properties of
  ideals.
\newblock {\em Mem. Coll. Sci. Univ. Kyoto. Ser. A. Math.}, 30:165--175, 1957.

\bibitem[Nar83a]{MR684627}
R.~Narasimhan.
\newblock Hyperplanarity of the equimultiple locus.
\newblock {\em Proc. Amer. Math. Soc.}, 87(3):403--408, 1983.

\bibitem[Nar83b]{MR715853}
R.~Narasimhan.
\newblock Monomial equimultiple curves in positive characteristic.
\newblock {\em Proc. Amer. Math. Soc.}, 89(3):402--406, 1983.

\bibitem[Oda73]{MR0472824}
Tadao Oda.
\newblock Hironaka's additive group scheme.
\newblock In {\em Number theory, algebraic geometry and commutative algebra, in
  honor of Yasuo Akizuki}, pages 181--219. Kinokuniya, Tokyo, 1973.

\bibitem[Oda83]{MR723465}
Tadao Oda.
\newblock Hironaka's additive group scheme. {II}.
\newblock {\em Publ. Res. Inst. Math. Sci.}, 19(3):1163--1179, 1983.

\bibitem[Oda87]{MR894302}
Tadao Oda.
\newblock Infinitely very near singular points.
\newblock In {\em Complex analytic singularities}, volume~8 of {\em Adv. Stud.
  Pure Math.}, pages 363--404. North-Holland, Amsterdam, 1987.

\bibitem[Par99]{MR1714830}
Kapil~H. Paranjape.
\newblock The {B}ogomolov-{P}antev resolution, an expository account.
\newblock In {\em New trends in algebraic geometry (Warwick, 1996)}, volume 264
  of {\em London Math. Soc. Lecture Note Ser.}, pages 347--358. Cambridge Univ.
  Press, Cambridge, 1999.

\bibitem[Pil04]{P_BIRS}
Olivier Piltant.
\newblock Applications of ramification theory to resolution of
  three-dimensional varieties.
\newblock {\em A talk at the ``Workshop on resolution of singularities,
  factorization of birational mappings, and toroidal geometry'' at Banff
  International Research Station, Banff, during 11-16 December}, 2004.


\bibitem[Spi04]{S_BIRS}
Mark Spivakovsky.
\newblock Puiseaux expansions, a local analogue of Nash's space 
  of arcs and the local uniformization theorem.
\newblock {\em A talk at the ``Workshop on resolution of singularities,
  factorization of birational mappings, and toroidal geometry'' at Banff
  International Research Station, Banff, during 11-16 December}, 2004.

\bibitem[Tei03]{MR2018565}
Bernard Teissier.
\newblock Valuations, deformations, and toric geometry.
\newblock In {\em Valuation theory and its applications, Vol. II (Saskatoon,
  SK, 1999)}, volume~33 of {\em Fields Inst. Commun.}, pages 361--459. Amer.
  Math. Soc., Providence, RI, 2003.

\bibitem[Vil89]{MR985852}
Orlando E.~Villamayor.
\newblock Constructiveness of {H}ironaka's resolution.
\newblock {\em Ann. Sci. \'Ecole Norm. Sup. (4)}, 22(1):1--32, 1989.

\bibitem[Vil92]{MR1198092}
Orlando E.~Villamayor.
\newblock Patching local uniformizations.
\newblock {\em Ann. Sci. \'Ecole Norm. Sup. (4)}, 25(6):629--677, 1992.

\bibitem[W{\l}o05]{MR2163383}
Jaros{\l}aw W{\l}odarczyk.
\newblock Simple {H}ironaka resolution in characteristic zero.
\newblock {\em J. Amer. Math. Soc.}, 18(4):779--822 (electronic), 2005.

\bibitem[Zar40]{MR0002864}
Oscar Zariski.
\newblock Local uniformization on algebraic varieties.
\newblock {\em Ann. of Math. (2)}, 41:852--896, 1940.
\end{thebibliography}
\end{document}